\documentclass[DIV=15,abstract=true]{scrartcl}

\usepackage{graphicx,amsmath,amssymb,amsthm,tikz-cd,mathrsfs,comment,adjustbox,epic,etoolbox,dynkin-diagrams,todonotes,quiver,enumitem,multirow,array,makecell,float,xparse,xstring,subfiles,verbatim,listings,nameref,setspace,epigraph,url}
\usepackage[T1]{fontenc}
\usepackage[utf8]{inputenc}
\usepackage{CJKutf8}

\def\row#1/#2!{#1_{\IfStrEq{#2}{}{n}{#2}} & \dynkin{#1}{#2}\\}

\usepackage{hyperref}
\hypersetup{
    colorlinks,
    citecolor=magenta,
    filecolor=red,
    linkcolor=blue,
    urlcolor=black
}

\usepackage[nameinlink]{cleveref}

\newtheorem{theoremalphabetical}{Theorem}[section]
\newtheorem{propositionalphabetical}[theoremalphabetical]{Proposition}
\AddToHook{env/propositionalphabetical/begin}{\crefalias{theoremalphabetical}{propositionalphabetical}}

\newtheorem{conjalphabetical}[theoremalphabetical]{Conjecture}
\AddToHook{env/conjalphabetical/begin}{\crefalias{theoremalphabetical}{conjalphabetical}}

\theoremstyle{plain}
\newtheorem{theorem}{Theorem}[section]

\newtheorem{lemma}[theorem]{Lemma}
\AddToHook{env/lemma/begin}{\crefalias{theorem}{lemma}}

\newtheorem{proposition}[theorem]{Proposition}
\AddToHook{env/proposition/begin}{\crefalias{theorem}{proposition}}

\newtheorem{cor}[theorem]{Corollary}
\AddToHook{env/cor/begin}{\crefalias{theorem}{cor}}

\newtheorem{conj}[theorem]{Conjecture}
\AddToHook{env/conj/begin}{\crefalias{theorem}{conj}}

\newtheorem{cordef}[theorem]{Corollary/Definition}
\AddToHook{env/cordef/begin}{\crefalias{theorem}{cordef}}

\AddToHook{env/propdef/begin}{\crefalias{theorem}{propdef}}

\theoremstyle{definition}

\newtheorem{definition}[theorem]{Definition}\AddToHook{env/definition/begin}{\crefalias{theorem}{definition}}

\newtheorem{exmp}[theorem]{Example}
\AddToHook{env/exmp/begin}{\crefalias{theorem}{exmp}}

\newtheorem{rem}[theorem]{Remark}
\AddToHook{env/rem/begin}{\crefalias{theorem}{rem}}

\newtheorem{algorithm}[theorem]{Algorithm}
\AddToHook{env/algorithm/begin}{\crefalias{theorem}{algorithm}}

\crefname{theorem}{Theorem}{Theorems}
\crefname{lemma}{Lemma}{Lemmata}
\crefname{proposition}{Proposition}{Propositions}
\crefname{cor}{Corollary}{Corollaries}
\crefname{definition}{Definition}{Definitions}
\crefname{rem}{Remark}{Remarks}
\crefname{exmp}{Example}{Examples}
\crefname{conj}{Conjecture}{Conjectures}
\crefname{algorithm}{Algorithm}{Algorithms}
\crefname{cordef}{Corollary/Definition}{Corollaries/Definitions}
\crefname{propdef}{Proposition/Definition}{Propositions/Definitions}
\crefname{section}{Section}{Sections}
\crefname{chapter}{Chapter}{Chapters}
\crefname{table}{Table}{Tables}
\crefname{rqenumeratei}{item}{items}
\crefname{rqenumerateii}{item}{items}
\Crefname{rqenumerateii}{Item}{Items}

\newcounter{proofcntr}
\newcounter{proofcntr2}

\NewDocumentCommand{\prooflabell}{m m}{%
  \refstepcounter{proofcntr2}%
  \label{#2}%
}

\newlist{rqenumerate}{enumerate}{2}
\setlist[rqenumerate,1]{%
    label=\arabic*.,
    ref=\arabic*,
    leftmargin=*, labelindent=0pt, itemsep=0.5em}
\setlist[rqenumerate,2]{%
    label=\bfseries\itshape RQ \arabic{rqenumeratei}.\arabic*.,
    ref=\arabic{rqenumeratei}.\arabic*,
    leftmargin=1.25em, labelindent=0pt, itemsep=0.25em}    

\AtBeginEnvironment{exmp}{\pushQED{\qed}}
\AtEndEnvironment{exmp}{\popQED\endexmp}

\AtBeginEnvironment{rem}{\pushQED{\qed}}
\AtEndEnvironment{rem}{\popQED\endrem}

\AtBeginEnvironment{algorithm}{\pushQED{\qed}}
\AtEndEnvironment{algorithm}{\popQED\endalgorithm}

\DeclareMathOperator{\CH}{CH}
\DeclareMathOperator{\Pic}{Pic}
\DeclareMathOperator{\Ind}{Ind}

\DeclareMathOperator{\cl}{cl}
\DeclareMathOperator{\Mod}{Mod}
\DeclareMathOperator{\Rr}{R}
\DeclareMathOperator{\supp}{supp}
\DeclareMathOperator{\Cyc}{Cyc}

\DeclareMathOperator{\Gl}{GL}

\DeclareMathOperator{\Hom}{Hom}

\DeclareMathOperator{\im}{im}

\DeclareMathOperator{\htt}{ht}

\newcommand{\overbar}[1]{\mkern 2mu\overline{\mkern-4mu#1\mkern-1.5mu}\mkern 2mu}

\newcommand{\Pro}{\mathbb{P}}
\newcommand{\F}{\mathbb{F}}
\newcommand{\Q}{\mathbb{Q}}

\newcommand{\N}{\mathbb{N}}
\newcommand{\A}{\mathbb{A}}
\newcommand{\mO}{\mathcal{O}}

\newcommand{\Z}{\mathbb{Z}}

\newcommand{\minn}{\mathrm{min}}
\newcommand{\1}{\mathbf{1}}
\newcommand{\St}{\mathbf{St}}

\newcommand{\Id}{\mathrm{Id}}

\newcommand{\G}{\mathbb{G}}

\newcommand{\Ql}{\overline{\Q_\ell}}

\newcommand{\et}{\mathrm{\acute{e}t}}
\newcommand{\mot}{\mathrm{mot}}
\newcommand{\Lich}{\mathrm{Lich}}
\newcommand{\pr}{\mathrm{pr}}
\newcommand{\upr}{\mathrm{upr}}
\newcommand{\Sch}{\mathrm{Sch}}

\title{On the Picard group of Deligne--Lusztig varieties}
\author{Erik Barinaga}
\date{\today}

\begin{document}

\maketitle

\begin{abstract}
    We study Picard groups of classical Deligne--Lusztig varieties in any Dynkin type, both in the integral and rational sense (i.e. we tensor the abelian group with an extension of $\Q$), and discuss the implications on their $\ell$-adic cohomology through work of Lusztig and Digne--Michel--Rouquier. 
    More precisely, we give an explicit integral description in the longest element case and discuss some additional cases in untwisted type $A$.
    For the rational Picard groups, we give a recursive formula whose starting point concerns Coxeter Deligne--Lusztig varieties.
    In the case of untwisted type $A$, this covers all Deligne--Lusztig varieties associated to Weyl group elements of full support.  
    This complements Hansen's previous work on Coxeter Deligne--Lusztig varieties in classical type.
\end{abstract}

\tableofcontents

\section{Introduction}\label{Introduction}

In their landmark paper \cite{originalpaper}, Deligne and Lusztig constructed representations of finite groups of Lie type using the $\ell$-adic cohomology of now called Deligne--Lusztig varieties. 
The $\ell$-adic cohomology of such varieties associated to Coxeter elements was computed by Lusztig in \cite{Lusztig1976}. 
\'Etale cohomology with integral and modular coefficients and other cohomology theories entirely have also been explored to gain further insight into these varieties; see, for example, \cite{BonRouq_CoxOrbModReps}, \cite{Dudas_CoxOrbBrauerTreesI}, \cite{Dudas_CoxOrbBrauerTreesII}, \cite{WangCompDL}, \cite{GK_crys_DL}, and \cite{hansen,hansen2002picard}.
The aim of this work is to build on the aforementioned contributions of Hansen and to study the Picard group of Deligne--Lusztig varieties in greater generality.

We start by briefly explaining the contents of this work.
Let $G_0$ be a connected reductive group over a finite field $\F_q$ of characteristic $p\neq \ell$, and let $G$ be its base change to $\overline{\F_q}$ with geometric Frobenius endomorphism $F$. 
Let $T\subset B$ be a maximal torus and a Borel subgroup of $G$, both $F$-stable.
Let $W=N_G(T)/T$ be the Weyl group of $G$, and for $w\in W$ denote by $X(w)\subset G/B$ the corresponding Deligne--Lusztig variety. 
One conjecture\footnote{Although \Cref{conjintroduction} is weaker than the stated version in \cite[Conjecture 1]{hansen}, it became clear in a private conversation with Søren Have (formerly Søren Have Hansen) that the conjecture was meant to be implicitly restricted to $F$-Coxeter elements.} that Hansen stated (\cite[Conjecture 1]{hansen}) is the following. 

\begin{conjalphabetical}[Hansen]\label{conjintroduction}
    Let $w\in W$ be $F$-Coxeter\footnote{A notion of a Coxeter element which takes the Frobenius into account, see \cref{Cox_supp}. If $G^F$ is untwisted, then the $F$-Coxeter elements are precisely the Coxeter elements.} and $i\in \N$ such that $i<l(w)$.
    Then $\CH_i(X(w))$ has rank zero, i.e., $\CH_i(X(w))_{\Q}=0$. 
\end{conjalphabetical} 

Hansen was able to prove his conjecture if $G^F$ is of classical type and the Chow group is the Picard group, i.e., if $i=l(w)-1$ (\cite[Corollary 4]{hansen2002picard}).
If $G^F$ is of untwisted type $A$ and $w$ is the standard Coxeter element, then $X(w)$ is well-known to be the (finite) Drinfeld halfspace (\cite[2.2]{originalpaper}) and hence \Cref{conjintroduction} is even true integrally for any $i<l(w)$ (\cite[Theorem 5]{hansen2002picard}). 
Hansen was also able to show that Deligne--Lusztig cycles form a $\Q$-basis for $\CH_*(G/B)_{\Q}$, which implies that the Chow groups of $X(w_0)$, where $w_0\in W$ is the longest element, have rank zero (\cite[Section 3.3]{hansen}). 
One of our results is an explicit description of $\Pic(X(w_0))$ as a torsion abelian group for every finite group of Lie type except the Suzuki and Ree groups, stated as follows. 
Denote by $S\subset W$ the subset of simple reflections.

\begin{propositionalphabetical}[\cref{Picard_longest_elem}]\label{Picard_longest_elem_introduction}
    The action of $G^F$ on $\Pic(X(w_0))$ is trivial, and an explicit description of the abelian group $\Pic(X(w_0))$ is given in the following table in the case where $G^F$ is not of type $A_1, \prescript{2}{}{B_2}, \prescript{2}{}{G_2}$ or $\prescript{2}{}{F_4}$.
\end{propositionalphabetical}

\begin{table}[H]
\centering
\begin{tabular}{ m{10em} | m{15em}} 
 \makecell{Type of $W$} & \makecell{$\Pic(X(w_0))$} \\ [0.5ex] 
 \hline\hline
 \\[-1em]
 \makecell{$A_{2m}$} & \makecell{$(\Z/(q^2-1)\Z)^{m}$} \\
 \\[-1em]
 \hline 
 \\[-1em]
 \makecell{$A_{2m+1}$, $m>0$} & \makecell{$(\Z/(q^2-1)\Z)^{m}\oplus \Z/(q+1)\Z$}  \\
 \\[-1em]
 \hline 
 \\[-1em]
 \makecell{$\prescript{2}{}{A_{n}}$, $C_n$, $D_n (n\text{ even})$, \\
$\prescript{2}{}{D_n} (n\text{ odd})$, $F_4$, $E_7$,\\ $E_8$, $G_2$, $\prescript{2}{}{E_6}$} & \makecell{$(\Z/(q+1)\Z)^m$ where $m\coloneq |S|$}\\
 \\[-1em]
 \hline 
 \\[-1em]
 \makecell{$D_n (n\text{ odd})$,\\
$\prescript{2}{}{D_n} (n\text{ even})$} & \makecell{$(\Z/(q+1)\Z)^{n-2}\oplus \Z/(q^2-1)\Z$} \\
 \\[-1em]
 \hline 
 \\[-1em]
 \makecell{$\prescript{3}{}{D_4}$} & \makecell{$\Z/(q+1)\Z\oplus \Z/(q^3+1)\Z$} \\
\\[-1em]
 \hline 
 \\[-1em]
 \makecell{$E_6$} & \makecell{$(\Z/(q+1)\Z)^{2}\oplus (\Z/(q^2-1)\Z)^2$} 
\end{tabular}
\label{table:2_introduction}
\end{table} 

The proof comes down to expressing the cycles $[\overline{X(sw_0)}]$ for $s\in S$ as a $\Z$-linear combination of Schubert cycles in $\Pic(G/B)$. 
This can be done using a result of Kim (\cite[Theorem 3.10]{kim2016}), as explained in \cref{subsec_longest}.
Moreover, we can combine his result with Schubert calculus and fine Deligne--Lusztig varieties (these are analogs in partial flag varieties) to study $\Pic(X(sw_0))$ for $s\in S$ in the case where $G^F$ is of type $A_{n}$ for $n\geq 2$ (see \cref{sec:other_weyl_gp_elems}). 
Namely, we show that fine Deligne--Lusztig cycles form a $\Q$-basis of $\Pic(G/P(s))_{\Q}$ which is an analogue of Hansen's result for the full flag variety in codimension $1$.
This implies that $\Pic(X(sw_0))$ is torsion if $s$ is not the middle vertex in the corresponding Dynkin diagram.
We explain how our methods can be used to explicitly compute $\Pic(X(sw_0))$ integrally, and illustrate this with an example.

Notably, Hansen's work does not feature an example of a Deligne--Lusztig variety $X(w)$ whose Chow group $\CH_i(X(w))$ has positive rank for some $i<l(w)$. 
We establish that such examples exist for the Picard group (where $w$ is not $F$-Coxeter, in agreement with \Cref{conjintroduction}).
Moreover, these examples give rise to (nonzero) $G^F$-representations over $\Ql$ if we consider the tensor product $\Pic(X(w))_{\Ql}\coloneq \Pic(X(w))\otimes_{\Z}\Ql$. 
Viewed from this representation-theoretic lens, one can use the result \cite[Corollary 2.7]{Lusztig1976} by Lusztig (also observed by Hansen in \cite[3.1.3]{hansen}) to conclude that \Cref{conjintroduction} is at least true on the principal series part.
In studying these representations beyond the $F$-Coxeter case, we obtain our main result stated as the following \Cref{inductive_step_general_introduction}. 
To prove it, we use, among others as explained later, tools from \cite{digne2006cohomologie} and \cite{orlik2014cohomology}.  
To state it, we explain some notation. 
We denote by $F^+$ the monoid generated by $S$, and one can define for every $w\in F^+$ a Deligne--Lusztig variety $X(w)$ in a similar way as for Weyl group elements (in the sense of \cite{digne2006cohomologie}, see \cref{def_gen_DL_var}). 
If $S_F$ is the set of $F$-orbits of $S$, then we let $\varphi\colon S\to S_F$ be the function mapping a simple reflection $s$ to its $F$-orbit $F\cdot s$. 
We say $w\in F^+$ has \textit{full $F$-support} if every fiber of $\varphi$ contains at least one $s\in S$ appearing in the word $w$. 
For $J\subset S$, we set $i_J\coloneq \Ind^{G^F}_{P_J^F}\1_{P_J^F}$.

\begin{theoremalphabetical}[\cref{inductive_step_general}]\label{inductive_step_general_introduction}
    Let $s\in S$ and $v\in F^+$ such that $w\coloneq svF(s)$ has full $F$-support.
    \begin{rqenumerate}
        \item \label{inductive_step_general_1_introduction}Suppose the tensored higher Chow groups $\CH^1(X(w),1)_{\Ql}$ and $\CH^1(X(vF(s)),1)_{\Ql}$ vanish. 
        Then we have an isomorphism of $G^F$-representations \[\Pic(X(w))_{\Ql}\cong (\Pic(X(vF(s)))_{\Ql}\oplus i_{\varphi^{-1}(\supp_F(v))})/\1_{G^F}.\]
        \item \label{inductive_step_general_2_introduction}Suppose $vF(s)$ is $F$-Coxeter, then we have an isomorphism of $G^F$-representations \[\Pic(X(w))_{\Ql}\cong \Pic(\overline{X(vF(s))})_{\Ql}/(\1_{G^F}\oplus \bigoplus_{J\in S_F\setminus \{F\cdot s\}}i_{\varphi^{-1}(S_F\setminus \{J\})}).\]
    \end{rqenumerate}
\end{theoremalphabetical}

Although \Cref{inductive_step_general_introduction} appears to require $w$ to have a specific form, a result of Geck--Kim--Pfeiffer (\cite[Theorem 2.6]{GECK2000570}) shows that any Weyl group element that is not of minimal length in its conjugacy class can be brought into this form after finitely many cyclic shifts.
Moreover, it turns out that these cyclic shifts induce isomorphisms on the $\ell$-adic cohomology and rational Chow groups of the corresponding Deligne--Lusztig varieties.
This gives one a tactic to recursively compute $\Pic(X(w))_{\Ql}$ for $w\in W$ not of minimal length in its conjugacy class. 
Additionally, we have good control over the vanishing conditions in \cref{inductive_step_general_1} of \Cref{inductive_step_general_introduction}.
More precisely, we use the Beilinson--Lichtenbaum conjecture to show the following.
For the convenience of the reader, we state the result in a simplified manner here.  

\begin{propositionalphabetical}[\cref{cycle_map_lemma_app}]\label{cycle_map_lemma_app_introduction}
    Let $w\in F^+$. 
    If $H^1_{\et}(X(w),\Ql)^{B^F}$ vanishes, then so does $\CH^1(X(w),1)_{\Ql}$. 
\end{propositionalphabetical}

If one restricts to the case where $G^F$ is of untwisted type $A$, then Orlik was able to prove in \cite[Proposition 7.6]{orlik2014cohomology} that $H^1_{\et}(X(w),\Ql)$ vanishes if $w$ is not Coxeter (we noticed that this statement requires an additional support condition, and we correct it in \cref{orlik_codim_11}).
We show that a similar phenomenon also carries over to all other types, in the sense of the following result.
The proof uses ideas of Lusztig in \cite{Lusztig1976}, and in fact we use \Cref{cycle_map_lemma_app_introduction} and the following \Cref{H1_vanishing_other_gps_introduction} to prove \Cref{inductive_step_general_introduction}.
For a $G^F$-representation $V$, we denote by $V_\upr$ and $V_\pr$ its (unipotent) principal series part, and for $J\subset S$ we denote by $W_J\subset W$ the subgroup generated by $J$.

\begin{propositionalphabetical}[\cref{H1_vanishing_other_gps}]\label{H1_vanishing_other_gps_introduction}
        Let $s\in S$ and $v\in W_{S\setminus F\cdot s}$ such that $vF(s)$ is $F$-Coxeter.
    \begin{enumerate}
        \item \label{H1_vanishing_other_gps1_introduction}
        If $G^F$ is not of type $\prescript{2}{}{A_{2n}}, \prescript{2}{}{B_2}, \prescript{2}{}{F_4}$ or $\prescript{2}{}{G_2}$ then $H^{1}_{\et}(X(svF(s)),\Ql)=0$.
        \item \label{H1_vanishing_other_gps2_introduction} 
        If $G^F$ is of type $\prescript{2}{}{A_{2n}}, \prescript{2}{}{B_2}, \prescript{2}{}{F_4}$ or $\prescript{2}{}{G_2}$, then $H^{1}_{\et}(X(svF(s)),\Ql)_\pr=0$.
    \end{enumerate}
\end{propositionalphabetical}

In order to compute $\Pic(X(w))_{\Ql}$ for $w\in W$ recursively using \Cref{inductive_step_general_introduction} and Hansen's results for $F$-Coxeter elements, one needs to understand the combinatorics of $W$ due to the aforementioned \cite[Theorem 2.6]{GECK2000570}.  
If $W$ is of type $A$, then the elements in $W$ of minimal length in their conjugacy class are Coxeter in some parabolic subgroup of $W$ (\cite[Example 3.1.16]{GECK2000570}).
This means that if $G^F$ is of untwisted type $A$, then the preceding results are enough to recursively compute $\Pic(X(w))_{\Ql}$ for any $w\in W$ non-Coxeter and of full support. 
For other types, the main obstruction is that there are such minimal length elements $w$ (henceforth called \textit{elements of height zero}) which are of full $F$-support, not $F$-Coxeter and unequal to $w_0$. 
There appear to be no general results on the cohomology of the corresponding varieties $X(w)$ in the literature.
An application of \Cref{inductive_step_general_introduction} is the following. 
Denote by $\overbar{X}(w)$ the (smooth) Bott--
Samelson--Demazure--Hansen compactification of $X(w)$, where $w\in F^+$ (see \cref{smooth_compactification}).
Using results of Lusztig and Digne--Michel--Rouquier (\cite[Sections 2 and 3]{Lusztig1998} and \cite[Proposition 3.4.2, Remarque 3.4.3]{digne2006cohomologie}), one sees that the cycle map $\cl\colon \CH^k(\overbar X(w))_{\Ql,\upr}\to H^{2k}_{\et}(\overbar X(w),\Ql)_{\upr}$ on the unipotent principal series part is an isomorphism of $G^F$-representations for any $w\in F^+$ and $k\in \N$ (see \Cref{sec:cycle_class_map}).
This implies that the cycle map $\cl\colon \Pic(X(w))_{\Ql,\upr}\to H^2_{\et}(X(w),\Ql)_{\upr}$ is injective for any $w\in F^+$, which means that by computing the rational Picard group one is computing at least part of the $\ell$-adic cohomology group $H^2_{\et}(X(w),\Ql)$. 
In the case that $G^F$ is of untwisted type $A_n$, we find a sufficient requirement for $\cl$ to be an isomorphism. 
Further, we use this to check that $\cl$ is indeed an isomorphism for all $w\in W$ of full support and non-Coxeter if $1\leq n\leq 4$. 

This preprint is organised as follows. 
In \cref{sec_DLvars} we define Deligne--Lusztig varieties associated to full and partial flag varieties, and recall some of their properties from the literature. 
\Cref{sec_integral_picard} and \cref{sec_rat_Picard_group} study $\Pic(X(w))$ integrally and rationally (i.e., we tensor the abelian group with $\Ql$) respectively for the various $w$ specified above.
We consider the case where $w$ is the longest element in \cref{subsec_longest}, and in \cref{sec:other_weyl_gp_elems} we study $\Pic(G/P(s))_{\Q}$ and $\Pic(X(sw_0))$ while restricting to groups of untwisted type $A$ as mentioned before.
For \cref{sec_rat_Picard_group}, we first explain in \cref{sec:motivic_cohom} how we use the Beilinson--Lichtenbaum conjecture to derive \Cref{cycle_map_lemma_app_introduction}. 
We discuss the Geck--Kim--Pfeiffer theorem \cite[Theorem 2.6]{GECK2000570}, cyclic shift classes and their consequences for the cohomology of Deligne--Lusztig varieties in \cref{sec:cyclic_shifts}.
In \cref{sec:cycle_class_map} we explain the consequences of the aforementioned results of Lusztig and Digne--Michel--Rouquier to cycle class maps of Deligne--Lusztig varieties. 
In \cref{sec:the_cox_case}, we compute some rational (higher) Chow groups of Coxeter Deligne--Lusztig varieties. 
The main results of \cref{sec_rat_Picard_group} are proved in \cref{sec:the_general_case}.
We restrict to $G^F$ of type $A$ in \cref{section_type_A} and derive some additional results in that case.
Finally, we specify our Dynkin diagram conventions and compute $\Pic(X(w))_{\Ql}$ for $w\in W$ of full support and positive height in some small rank examples in the \nameref{Appendix}.

\subsection*{Acknowledgements}

As this preprint is based on my PhD thesis, I would like to thank my advisor, Sascha Orlik, for suggesting this interesting topic and for his valuable feedback on earlier versions of my dissertation.
I am very grateful to Olivier Dudas, Jean Michel and Xuhua He for fruitful discussion and for answering my questions on their work.
I want to thank Andreas Bode, Sergio Cía, Nicolas Dupré, Arnaud Eteve, Jan Hennig, Dennis Peters, Dzoara Núñez Ramos, Julian Reichardt, Fei Ren, Yingying Wang, Ziqian Yin, Alexander Ziegler and Hugo Zock for helpful discussions.
This work was realised while I was a PhD student in the research training group \textit{GRK 2240: Algebro-Geometric Methods in Algebra, Arithmetic and Topology}, which is funded by the Deutsche Forschungsgemeinschaft. 
I want to explicitly mention the principal investigators Kay R\"ulling, Britta Späth and Matthias Wendt, with whom I had very useful discussions for which I wish to thank them. 

\subsection*{AI statement}

The author did not use AI in any way in the preparation of this preprint.
This preprint has been written solely by the author and the author bears full responsibility for it.

\subsection{Notation and conventions}\label{not_and_conv}

\begin{itemize}
    \item For us, a variety is a separated and reduced scheme of finite type over an algebraically closed field. 
    \item For $G$ a group and $H\leq G$ a subgroup, we set $\prescript{g}{}H\coloneq gHg^{-1}$ for $g\in G$.
    \item Let $\ell$ and $p$ be unequal primes. 
    If $Y$ is a scheme over $\overline{\F_p}$ containing a closed subscheme $Z$, then we denote by $H^i(Y)$, $H^i_c(Y)$ and $H_Z^i(Y)$ respectively the $\ell$-adic cohomology, the $\ell$-adic cohomology with compact support and the $\ell$-adic cohomology with support in $Z$ of $Y$ in degree $i$, all three with coefficients in $\Ql$. 
    As usual, we denote by $(m)$ the Tate twist of degree $m$ where $m$ is an integer.
    We use the convention $H^2_c(\A^1)=\Ql(-1)$.
    \item If $Y$ is a scheme of dimension $d$ over a field, then for $k\in\N$ we consider its Chow group $\CH_k(Y)$ of $k$-cycles as treated in \cite[Section 1.3]{fulton2012intersection}.
    We sometimes also use the notation $\CH^k(Y)$ to denote its Chow group of $(d-k)$-cycles.
    If $Y$ is smooth, then there is a canonical isomorphism $\Pic(Y)\cong \CH^1(Y)$ by the Auslander--Buchsbaum theorem (see \cite[Theorem 20.3]{Ma} and \cite[Example 2.1.1]{fulton2012intersection}).
    We will mainly be interested in studying the Picard group of (smooth) Deligne--Lusztig variety, and so we will often identify its Picard and codimension $1$ Chow group.
    \item If $W$ is a (finite) Weyl group, then we denote by $w_0\in W$ the unique longest element with respect to the length function $l\colon W\to \N$ (for more information, see for example \cite[Section 2.3]{bjorner2006combinatorics}).
\end{itemize}

\section{Deligne--Lusztig varieties in full and partial flag varieties}\label{sec_DLvars}

Before treating Deligne--Lusztig varieties associated to full and partial flag varieties, we briefly recall the following.

\subsection{Finite groups of Lie type, Weyl groups and Dynkin diagrams}\label{sec:fin_gps_weyl_dynkin}

We briefly recall finite groups of Lie type.
For more information, see for example \cite[Chapter 4]{Digne_Michel_2020}, \cite[Part III]{Malle_Testerman_2011} or \cite[Chapter 1]{Carter1972}.
Let $q$ be a power of a prime $p$\footnote{Here we allow rational exponents.}.
Let $d$ be a positive integer such that $q^d\in \N$, $k_0\coloneq \F_{q^d}$ and $\F$ be an algebraic closure of $\F_p$ containing $k_0$.
Let $G_0$ be a connected reductive group over $k_0$ with Frobenius $F_0\colon G_0\to G_0$ defined by raising the functions on $G_0$ to the $q^d$-th power. 
Let $G$ be the base change of $G_0$ to $\F$ and $F\colon G\to G$ an endomorphism such that $F^d=F_0\otimes_{k_0}\Id$ is a geometric Frobenius (see \cite[Definition 4.1.6]{Digne_Michel_2020}).
We will refer to $F$ as the \textit{Frobenius} for simplicity, it is also called a Frobenius root in the literature (\cite[Definition 4.2.4]{Digne_Michel_2020}). 
In this text, we will have $d=1$ or $d=2$ and $q$ is an odd power of $\sqrt{2}$ or $\sqrt{3}$.
The group of fixed points $G^F$ is called \textit{a finite group of Lie type}. 
Let $T\subset B$ be a maximal torus and a Borel subgroup of $G$, both $F$-stable.
Denote by $\Phi$ the reduced crystallographic root system relative to $T$ and by $\Phi^+$ the positive subsystem determined by $B$ (see \cite[Theorem 2.3.1]{Digne_Michel_2020}).
We assume that $\Phi$ is irreducible.
Let $W\coloneq N_G(T)/T$ be the (finite) Weyl group of $G$, with simple reflections $S$.
Also, denote by $\leq$ the Bruhat order of $W$ and by $l\colon W\to \N$ the length function.
The Frobenius induces a group automorphism $F\colon W\to W$ such that $F(S)=S$ (see \cite[Lemma 4.2.21]{Digne_Michel_2020}).
Hence, it induces an automorphism of the underlying connected Dynkin diagram of $W$.
Denote by $\delta$ the smallest natural number such that $F^\delta$ is the identity on $W$.
The type of the Dynkin diagram together with the graph automorphism induced by $F$ determines the type of $G^F$, which must be one of the following 
\begin{align*}
    &A_n (n\geq 1), B_n (n\geq 2), C_n (n\geq 3),D_n (n\geq 4),E_n (n=6,7,8),F_4, G_2,\\
    &\prescript{2}{}A_n(n\geq 2), \prescript{2}{}B_2,\prescript{2}{}D_n(n\geq 4),\prescript{3}{}D_4,\prescript{2}{}E_6,\prescript{2}{}F_4,\prescript{2}{}G_2.
\end{align*}
The superscript equals $\delta$ and the subscript equals $|S|$, and we follow the Dynkin diagram convention of \cite[Section 1.11]{Carter1972}.
For the convenience of the reader, the diagrams with labelled vertices are also displayed in the \nameref{Appendix}.
If $F$ acts nontrivially on the Dynkin diagram (equivalently, if $\delta>1$), we call $G^F$ \textit{twisted}.
If not, we call it \textit{untwisted}.
If $W$ is of type $A$, $C$ or $D$, we call $W$ and $G^F$ of \textit{classical type}.
If not, we call them of \textit{exceptional type} (note that $F$ may act nontrivially on $W$).
The types $\prescript{2}{}{B_2}$, $\prescript{2}{}{G_2}$ and $\prescript{2}{}{F_4}$ correspond to the \textit{Suzuki} and \textit{Ree} groups, and for these $d=2$ and $q$ is an odd power of $\sqrt{2}$, $\sqrt{3}$ and $\sqrt{2}$ respectively (these are also sometimes called the \textit{very twisted} types).
For all other types (including the untwisted ones), $d=1$.
Recall that the Weyl groups corresponding to types $B_n$ and $C_n$ are isomorphic, making the results obtained in this preprint for one type easily carry over to the other type.
Also recall that the Weyl group $W$ associated to type $A_n$ is isomorphic to the symmetric group $S_{n+1}$ via the mapping $W\to S_{n+1}$ where $s_i\mapsto (i\ i+1)$ for $1\leq i\leq n$ (see for example \cite[Example 1.2.3]{bjorner2006combinatorics}).
Very often, we will identify these groups and write $s_i=(i\ i+1)$. 

\subsubsection{Orbits and support in \texorpdfstring{$S$}{S}}

For $s\in S$, we denote by $F\cdot s\coloneq \{F^r(s):r\in \N\}$ the $F$-orbit of $s$ in $S$ and by $S_F$ the set of $F$-orbits in $S$.
There is a canonical map $\varphi\colon S\to S_F$.  

\begin{definition}\label{Cox_supp}
	Let $w\in W$, and write $w=s_1\cdots s_r$ as a reduced decomposition.
\begin{enumerate}
	\item The \textit{support} $\supp(w)$ of $w$ is defined as the set \begin{align*}
	\supp(w)\coloneq \{s_1,\ldots,s_r\}.
	\end{align*}
	We say that $w$ has \textit{full support in} $W$ if $\supp(w)=S$.
	By taking $F$-orbits of the $s_i$ instead, we can analogously define the \textit{$F$-support} $\supp_F(w)$ of $w$ as 
	\begin{align*}
	\supp_F(w)\coloneq \{F\cdot s_1,\ldots,F\cdot s_r\}.
	\end{align*}
	We say that $w$ has \textit{full $F$-support} if $\supp_F(w)=S_F$.
    For $w_1,\ldots,w_r\in W$, we denote 
	\begin{align*}
		\supp_F(w_1,\ldots,w_r)\coloneq \supp_F(w_1)\cup \cdots \cup \supp_F(w_r).
	\end{align*}
    Let $F^+$ denote the monoid (freely) generated by $S$.
    We can analogously define the notions above for elements of $F^+$ instead of $W$.
	\item We say $w$ is a \textit{Coxeter element} (respectively an \textit{$F$-Coxeter element}) if it has full support (respectively full $F$-support) and the $s_i$ (respectively $F\cdot s_i$) are pairwise distinct. 
\end{enumerate}
\end{definition}

Clearly, if $F$ acts trivially on $W$ then the two Coxeter and support notions defined in \cref{Cox_supp} coincide. 
Note that the definition of an $F$-Coxeter element is equivalent to the one given in \cite[(1.7)]{Lusztig1976} (Lusztig calls them Coxeter elements, but we prefer to have two separate notions if $F$ acts nontrivially).
They are also sometimes called \textit{twisted} Coxeter elements (\cite[Section 7]{Springer_regelemens}).

\subsection{Deligne--Lusztig varieties}\label{sec:DL_vars}

In this section we define the main objects of interest: Deligne--Lusztig varieties. 
Throughout this entire section, let $G^F$ be a finite group of Lie type as defined in \cref{sec:fin_gps_weyl_dynkin}.

\subsubsection{Full flag variety}\label{sec:Full_flag_variety}

We start by discussing classical and generalised Deligne--Lusztig varieties associated to the full flag variety $X$.
In \cref{background_partial_flag}, we discuss the counterparts of the classical varieties in partial flag varieties.

Before that, we first recall some useful notions. 
Set $X\coloneq G/B$, which can be identified as the set of Borel subgroups of $G$.
It is a smooth projective variety over $\F$, and is referred to as the full flag variety of $G$. 
For $w\in W$, we denote by $\Omega_{\Sch}(w)$ the corresponding Schubert cell $BwB/B\subset X$ and by $\Omega_w$ the corresponding Schubert variety (the Zariski closure $\overline{\Omega_\Sch(w)}$ in $X$). 
It is well-known that $\Omega_\Sch(w)\cong \A^{l(w)}_{\F}$.
Let $w\in W$ and define $\mO(w)$ to be the $G$-orbit of $(B,wBw^{-1})$ in $X\times X$. 
Following \cite[Section 2]{digne2006cohomologie}, let $\underline W$ be a set together with a bijection $W\to\underline W$. 
Give $\underline W$ the same length function as $W$, i.e. for any $\underline w\in \underline W$ set $l(\underline w)\coloneq l(w)$. 
For $\underline w\in \underline W$, set $\mO(\underline{w})\coloneq \overline{\mO(w)}$. 
Let $w_1,\ldots,w_r\in W\cup \underline W$, and construct the following fiber product \begin{align*}
\mO(w_1,\ldots,w_r)&\coloneq \mO(w_1)\times_{X} \cdots \times_{X}\mO(w_r)
\end{align*}
inductively by considering the first and second projections.
Explicitly, we may identify \begin{align*}
\mO(w_1,\ldots,w_r)&=\{(B_0,B_1,\ldots,B_r)\in X^{r+1}:\forall\hspace{0.25em}0\leq i\leq r-1, (B_i,B_{i+1})\in\mO(w_{i+1})\}.
\end{align*}
By construction, we have that if $w_1,\ldots, w_r\in W$ and $l(w_1\cdots w_r)=l(w_1)+\cdots +l(w_r)$ then we have a $G$-equivariant isomorphism $\mO(w_1,\ldots,w_r)\cong \mO(w_1\cdots w_r)$.
Denote by $\Gamma_F$ the graph of $F$ in $X\times X$. 

\begin{definition}\label{def_gen_DL_var}
For $w_1,\ldots,w_r\in W\cup \underline W$, we define the \textit{(generalised) Deligne--Lusztig variety} corresponding to $(w_1,\ldots,w_r)$ as the pullback \[\begin{tikzcd}
X(w_1,\ldots,w_r) \arrow[r] \arrow[d] & \Gamma_F \arrow[d,hook] \\
\mO(w_1,\ldots,w_r) \arrow[r]         & X\times X         
\end{tikzcd}\]in the category of $\overline{\F_q}$-schemes  where the bottom horizontal map is given by $(B_0,B_1,\ldots,B_r)\mapsto (B_0,B_r)$. 
If $r=1$ and $w\in W$, we recover the \textit{classical} Deligne--Lusztig variety denoted by $X(w)$.
If additionally $w$ is $F$-Coxeter, we call $X(w)$ a \textit{Coxeter variety}.

If $w_1,\ldots,w_r\in W$, we set $\overbar{X}(w_1,\ldots,w_r)\coloneq X(\underline{w_1},\ldots,\underline{w_r})$ and we write $X_G(w_1,\ldots,w_r)$ and $\overbar{X}_G(w_1,\ldots,w_r)$ if we want to emphasize the reductive group in which we are working.
\end{definition}

Explicitly, we may identify a generalised Deligne--Lusztig variety as \begin{align*}
X(w_1,\ldots,w_r)=\{(B_0,B_1,\ldots,B_r)\in X^{r+1}:\forall \hspace{0.25em}0\leq i\leq r-1, (B_i,B_{i+1})\in\mO(w_{i+1})\text{ and }F(B_0)=B_r\}.
\end{align*}

In \cref{smooth_compactification}, we discuss how to construct a smooth compactification of any Deligne--Lusztig variety $X(w)$.
This was first done in \cite[9.10]{originalpaper}, and later generalised in \cite[Paragraph 2.2.11 -- 2.3.7]{digne2006cohomologie}. 
We closely follow the latter reference, but first recall some relevant results from the literature.
Recall from \cite[Definition A1]{LusztigKazhdan} the definition of rational smoothness.
By purity (\cite[Theorem 16.1, Theorem 16.7]{milneLEC}), smoothness is a stronger property. 

\begin{proposition}\label{gen_iso_normal}
Let $w_1,\ldots,w_r\in W$.
\begin{rqenumerate}\label{rqgen_iso_normal}
\item \label{rqgen_iso_normal1}The variety $X(w_1,\ldots,w_r)$ is smooth and quasi-projective of pure dimension $l(w_1)+\cdots+l(w_r)$.
\item \label{rqgen_iso_normal11}Suppose $r=1$ and set $w\coloneq w_1$. 
The variety $X(w)$ carries a left $G^F$-action and is locally closed in $X$.
Its Zariski closure $\overline{X(w)}$ is given by 
\begin{align*}
    \overline{X(w)}=\displaystyle \dot\bigcup_{v\leq w} X(v) 
\end{align*} and is rationally smooth if and only if the corresponding Schubert variety $\Omega_w$ is. 
If $G$ is of type $A$, $D$ or $E$ then we have the previous equivalence for smoothness instead of rational smoothness.  
\item \label{rqgen_iso_normal2}If $\overline{\mO(w_i)}$ is smooth for every $1\leq i\leq r$, then $\overbar{X}(w_1,\cdots ,w_r)$ is smooth.
\item \label{rqgen_iso_normal3}The variety $\overbar{X}(w_1,\ldots, w_r)$ is projective and contains $X(w_1,\ldots, w_r)$ as an open dense subset.
\item \label{rqgen_iso_normal5}If $l(w_1\cdots w_r)=l(w_1)+\cdots+l(w_r)$ then there is an isomorphism \begin{align*}
X(w_1,\ldots,w_r)\cong X(w_1\cdots w_r). 
\end{align*}
\item \label{rqgen_iso_normal4}
Let $I\coloneq \varphi^{-1}(\supp_F(w_1,\ldots,w_r))$.
There is an $G^F\times F^m$-equivariant isomorphism 
\begin{align*}
	X(w_1,\ldots,w_r)\cong G^F/U_I^F\times_{L_I^F}X_{L_I}(w_1,\ldots,w_r)
\end{align*}
where $m\in\N$ is such that $F^m$ fixes $(w_1,\ldots,w_r)$.
In particular, the scheme $X(w_1,\ldots,w_r)$ has $|G^F/P_{I}^F|$ irreducible components and $G^F$ acts transitively on them. 
Hence, it is irreducible if and only if $I=S$.
The foregoing statements also hold for $\overbar{X}(w_1,\ldots,w_r)$.
\end{rqenumerate}
\end{proposition}
\begin{proof}
\Cref{rqgen_iso_normal1,rqgen_iso_normal2,rqgen_iso_normal3} can be found in \cite[Propositions 2.3.5, 2.3.6]{digne2006cohomologie}, using the fact that each $\overline{\mO(w_i)}$ contains $\mO(w_i)$ as an open dense subset.
\Cref{rqgen_iso_normal11} follows from \cite[Section 1.2, Definition 1.4]{originalpaper}, \cite[Lemme 2.2.9, Proposition 3.2.5]{digne2006cohomologie}, \cite[Theorem 6.0.4]{laksbilley} and \cite[Lemma 9.11]{originalpaper}.
\Cref{rqgen_iso_normal5} follows from the fact that $\mO(w_1,\ldots,w_r)\cong \mO(w_1\cdots w_r)$ (see also \cite[2.3.1]{digne2006cohomologie}).
For \cref{rqgen_iso_normal4}, see \cite[Lemma 2.2.1, Proposition 2.3.3]{WangCompDL} and \cite[Proposition 2.3.8]{digne2006cohomologie}.
The classical Deligne--Lusztig variety case already had proofs in \cite[Theorem 2]{bonnafé2006irreducibility} and later in \cite[Corollary 1.2, Proposition 4.3]{goertz2008connectednessdelignelusztigvarieties} and \cite[Proposition 3.3]{kim2016}. 
\end{proof}

\begin{rem}\label{rem_pattern_avoidance}
    In type $A$, we have the pattern avoidance criterion of Lakshmibai--Sandhya at our disposal characterizing smooth Schubert varieties \cite{Laks_Sand_SmSch_Sln}. 
\end{rem}

The following is \cite[Lemma 9.11]{originalpaper}.

\begin{cordef}[Deligne--Lusztig]\label{smooth_compactification}
For every $s\in S$, the variety $\overline{\mO(s)}$ is smooth.
In particular, if $w=s_{i_1}\cdots s_{i_r}\in W$ is a reduced decomposition then $\overbar{X}(w)\coloneq \overbar{X}(s_{i_1},\ldots, s_{i_r})$ is a smooth compactification of $X(w)$.
We write $\partial \overbar{X}(w)\coloneq \overbar{X}(w)\setminus X(w)$ for the closed complement. 
\end{cordef}
\begin{proof}
The first claim follows from the fact that either projection $\overline{\mO(s)}\to X$ is a $\Pro^1$-bundle (see \cite[1.2]{originalpaper} and \cite[Proposition 2.2.3]{digne2006cohomologie}). 
The second follows from \Cref{gen_iso_normal}. 
\end{proof}

An element $w=s_{i_1}\cdots s_{i_r}\in F^+$ is essentially a sequence $(s_{i_1},s_{i_2},\ldots,s_{i_r})\in S^r$. 
We write $X^{F^+}(w)$ for $X(s_{i_1},s_{i_2},\ldots,s_{i_r})$ and occasionally leave out the superscript if it is clear that $w\in F^+$. 
Using \cite[Section 1.4, Application 2]{Deligne_braids} one can show that the smooth compactification $\overbar{X}(w)$ is independent of the choice of reduced expression for $w\in W$.

\begin{proposition}\label{smooth_cox_comp}
Let $w=s_1\cdots s_r\in W$ be a reduced decomposition. 
If the $s_i$ are pairwise distinct then we have an isomorphism $\overbar{X}(w)\cong \overline{X(w)}$.
In particular, the condition is satisfied by $F$-Coxeter elements.
\end{proposition}
\begin{proof}
	This follows from the isomorphism $\mO(\underline{s_1},\ldots,\underline{s_r})\cong \overline{\mO(w)}$ (see \cite[Proposition 2.2.16]{digne2006cohomologie}).
	See also \cite[Remark 3.18]{orlik2014cohomology}.
\end{proof}

\begin{cordef}\label{partial_compactification}
	Let $v=s_{i_1}\cdots s_{i_r}\in F^+$ and $t\in S$.
	Define \begin{align*}
		X(\underline t v)&\coloneq X(\underline t,s_{i_1},\ldots ,s_{i_r}),\\
		X(v\underline t )&\coloneq X(s_{i_1},\ldots ,s_{i_r},\underline t).
	\end{align*}
	The schemes $X(\underline t v)$ and $X(v\underline t)$ are smooth. 
	Moreover, they respectively contain $X(tv)$ and $X(vt)$ as open subschemes and the closed complement in both cases is $X(v)$. 
\end{cordef}
\begin{proof}
	This is \cite[Proposition 2.3.5]{digne2006cohomologie}, see also \Cref{smooth_compactification}.
	The second claim follows from the open embedding $\mO(t)\subset \overline{\mO(t)}$. 
\end{proof}

The following lemma states that for particular generalized Deligne--Lusztig varieties, their $\ell$-adic cohomology with compact support is determined by the cohomology of a generalized Deligne--Lusztig variety of codimension $1$.
In the most basic case, this uses ideas discussed in \cite[Case 2 of the proof of Theorem 1.6]{originalpaper}.
The more general version is due to Digne--Michel--Rouquier, which we exhibit in \Cref{affine_bundle}.
The stated version is less general than their result, but suffices for our purposes.
As the conclusions on cohomology stem from an $\A^1$-bundle, we can make an analogous statement for Chow groups. 

\begin{lemma}[Digne--Michel--Rouquier]\label{affine_bundle}
Let $v\in F^+$ and $s\in S$.  
By \Cref{partial_compactification} $X(\underline s sv)$, $X(\underline sv)$, $X(vs\underline s)$ and $X(v\underline s)$ are smooth schemes.
There are affine line bundles $X(\underline s s v)\to X(\underline s v)$ and $X(vs\underline s)\to X(v\underline s)$ inducing isomorphisms \begin{align*}
	H^{i}_{c}(X(\underline s sv))&\cong H^{i-2}_c(X(\underline s v))(-1)\\
	H^{i}_c(X(vs\underline s))&\cong H^{i-2}_c(X(v\underline s))(-1)\\
	\CH_{i}(X(\underline s sv))&\cong \CH_{i-1}(X(\underline s v))\\
	\CH_{i}(X(vs\underline s))&\cong \CH_{i-1}(X(v\underline s))
\end{align*}for every $i\in \N$.
These are $G^F\times F^\delta$-equivariant in the first two cases and $G^F$-equivariant in the last two. 
\end{lemma}
\begin{proof}
The $\A^1$-bundles and $\ell$-adic cohomology statements are shown in (the proof of) \cite[Proposition 3.2.3(i)]{digne2006cohomologie} and in the analogous proposition discussed after the proof.
The claim for Chow groups is now given by \cite[Lemma 2.2]{totaro_2014}. 
\end{proof}

Once again, we use a result by Digne--Michel--Rouquier but state it in a less general form. 
See also \cite[Proposition 4.3]{orlik2014cohomology} for a less general version in untwisted type $A$.

\begin{proposition}[Digne--Michel--Rouquier]\label{univ_homeo_DL_general}
Let $w_1,w_2\in F^+$ such that \begin{align*}
l(w_1w_2)=l(w_1)+l(w_2)=l(w_2)+l(F(w_1))=l(w_2F(w_1)).
\end{align*}
There are universal homeomorphisms\begin{align*}
\sigma&\colon X(w_1w_2)\to X(w_2F(w_1)),\\
\tau&\colon X(w_2F(w_1))\to X(F(w_1w_2)),\\
\sigma^{(q)}&\colon X(F(w_1w_2))\to X(F(w_2)F^2(w_1))
\end{align*}
sitting in a commutative diagram 
\[
    \begin{tikzcd}
    X(w_1w_2) \arrow[r, "\sigma"] \arrow[d, "F"'] & X(w_2F(w_1)) \arrow[d, "F"] \arrow[ld, "\tau"] \\
    X(F(w_1w_2)) \arrow[r, "\sigma^{(q)}"]        & X(F(w_2)F^2(w_1)).                             
    \end{tikzcd}
\]
\end{proposition}
\begin{proof}
See \cite[Proposition 3.1.6]{digne2006cohomologie}. 
\end{proof}

\begin{rem}\label{subsec:elsofpropsupp}
    We end this section by explaining why we are only interested in elements $w\in F^+$ of full $F$-support. 
    This is because by \cref{rqgen_iso_normal4} of \cref{gen_iso_normal}, computing the cohomology (rational Chow groups or $\ell$-adic) of $X(w)$ or $\overbar{X}(w)$ comes down to computing the cohomology of a Deligne--Lusztig variety associated to an element of full $F$-support and then taking Harish-Chandra induction.
    More precisely, suppose $I\coloneq \supp_F(w)\subsetneq S_F$. 
    Then 
    \begin{align*}
        H^i_c(X(w))&\cong \Ql[G^F/U_I^F]\otimes_{\Ql[L_I^F]}H^i_c(X_{L_I}(w))\\
        &= \Rr^{G}_{L_I}H^i_c(X_{L_I}(w))
    \end{align*}
    by \cref{rqgen_iso_normal4} of \cref{gen_iso_normal} and the definition of Harish-Chandra induction.
    The same thing happens if we take $\CH^i(-)_{\Ql}$ instead (see \cite[Example 1.7.6]{fulton2012intersection}) and likewise for $\overbar{X}(w)$.
\end{rem}

\subsubsection{Partial flag varieties}\label{background_partial_flag}

Partial flag varieties contain two, in general different, types of Deligne--Lusztig varieties called \textit{parabolic} and \textit{fine} Deligne--Lusztig varieties.
We first recall Schubert varieties in partial flag varieties.

\begin{definition}\label{partial_Sch}
	Let $J\subset S$. 
	\begin{enumerate}
		\item Let $P_J\subset G$ be the standard parabolic subgroup generated by the simple reflections in $J$, and $L_J$ and $U_J$ the corresponding Levi subgroup and unipotent radical such that $P_J=U_J\rtimes L_J$.
		Denote by $\Phi_J$ the corresponding root subsystem of $\Phi$.
		If $J=\{s\}$ is a singleton, we also sometimes write $P(s)$ and $L(s)$ instead of $P_J$ and $L_J$ respectively. 
		We analogously define $W_J\subset W$ as the subgroup generated by $J$.
		\item We call the quotient $G/P_J$ \textit{the partial flag variety} corresponding to $J$. 
        We can identify the partial flag variety $G/P_J$ with the variety of parabolic subgroups of $G$ conjugate to $P_J$.
		If $J=\varnothing$ we recover the full flag variety. 
		\item Define \begin{align*}
			W^J\coloneq \{w\in W:l(ws)>l(w)\text{ for all } s\in J\}.
		\end{align*}It is the subset of minimal length representatives of $W/W_J$ by \cite[Proposition on page 19]{Humphreys_1990}.
        We can define $\prescript{J}{}W$ analogously.
		\item Let $P\coloneq P_J$ and for $w\in W^J$ define $\Omega_{\Sch}^P(w)\coloneq BwP/P\subset G/P$ and $\Omega^P_w\coloneq \overline{\Omega_{\Sch}^P(w)}$. 
	\end{enumerate}
\end{definition}

\begin{lemma}[Kempf]\label{Sch_birat}
Let $s\in S$ and consider the quotient map $\pi_s\colon X\to G/P(s)$.
Let $w\in W$. 
If $w\in W^{\{s\}}$, then $\pi_s$ restricts to a birational morphism $\Omega_w\to \Omega_w^{P(s)}$. 
If not, then $\dim\pi_s(\Omega_w)=l(w)-1<l(w)$. 
\end{lemma}
\begin{proof}
This is \cite[Section 2, Lemma 1]{kempf_sch}. 
\end{proof}

We start by defining parabolic Deligne--Lusztig varieties. 
Let $J\subset S$.
As the Frobenius $F\colon G\to G$ gives us $F(L_J)=L_{F(J)}$ and $F(U_J)=U_{F(J)}$ we get $F(P_J)=P_{F(J)}$.
This means that the induced Frobenius $X\to X$ on the full flag variety induces a Frobenius $G/P_J\to G/P_{F(J)}$ which we will also denote by $F$. 
For $w\in W$, denote by $\mO_J(w)$ the $G$-orbit of $(P_J,wP_{F(J)}w^{-1})$ in $(G/P_J)\times (G/P_{F(J)})$. 
It only depends on the double coset $W_JwW_{F(J)}$.

\begin{definition}
For $w\in W$, define the corresponding \textit{parabolic Deligne--Lusztig variety} to be the pullback \[\begin{tikzcd}
X_J(w) \arrow[r] \arrow[d] & \Gamma_F \arrow[d,hook] \\
\mO_J(w) \arrow[r,hook]         & G/P_J\times G/P_{F(J)}         
\end{tikzcd}\]in the category of $\overline{\F_q}$-schemes. Working out the definition, it may be identified as \begin{align*}
X_J(w)&= \{P\in G/P_J:(P,F(P))\in\mO_J(w)\}. 
\end{align*}
\end{definition}

\begin{rem}\label{rem:parab_covers}
	If $\prescript{w}{}J=F(J)$ then $X_J(w)$ comes with a natural \'etale cover $Y_J(\dot w)\to X_J(w)$ as in \cite[1.8]{originalpaper} where $\dot w \in N_G(T)$ is a lift of $w$. 
	In this case, we get $\prescript{w}{}{L_J}=L_{F(J)}$ and $\prescript{w}{}{U_J}=U_{F(J)}$ and we may hence identify 
	\begin{align*}
		X_J(w)=\{gP_J\in G/P_J\colon g^{-1}F(g)\in P_JwP_{F(J)}\}.
	\end{align*}
	The lift $Y_J(\dot w)$ is defined as the analogue in $G/U_J$.
	More precisely,
	\begin{align*}
		Y_J(\dot w)=\{gU_J\in G/U_J\colon g^{-1}F(g)\in U_J\dot wU_{F(J)}\}.
	\end{align*}
	These carry a left $G^F$-action and a right $L_J^{\dot w F}$-action, where the latter is the finite group of fixed points of $L_J$ under $\dot wF\coloneq \mathrm{ad}(\dot w)\circ F$. 
	These are central in the definition of \textit{Lusztig induction}, a first definition of which appeared in \cite{Lusztig1976IND}.
	These varieties were further studied in \cite{Digne_2014}.
    In this preprint, they will only appear in the proof of \cref{H1_height0_nonCox}.
	Note $X_{\varnothing}(w)=X(w)$ and $Y_{\varnothing}(\dot w)=Y(\dot w)$. 
\end{rem}

We can now define fine Deligne--Lusztig varieties. 

\begin{definition}\cite[2.1]{HE_2019}
For $w\in \prescript{J}{}W$, define the corresponding \textit{fine Deligne--Lusztig variety} as \begin{align*}
X_J\{w\}&=\{gP_J\in G/P_J:g^{-1}F(g)\in P_J\cdot_FBwB\}
\end{align*}
where $x\cdot_F y=xyF(x)^{-1}$. 
If $J=\{s\}$ is a singleton, then we sometimes write $X_s\{w\}$ to shorten notation. 
Note $X_\varnothing \{w\}=X(w)$.
\end{definition}

Although there are multiple definitions of these objects (see \cite[Section 4]{lusztig2006class}, \cite[Section 3.2]{he2008gstable} or \cite[Definition 2.1.6]{hoeve}), we will think of them as the images of classical Deligne--Lusztig varieties under the quotient map $X\to G/P_J$ (see \Cref{sneaky_definition}).
In \cite[Theorem 2.1.7.]{hoeve} we see that these fine Deligne--Lusztig varieties are uniquely determined by the following properties. 

\begin{theorem}[Hoeve]\cite[Theorem 2.1.7.]{hoeve}\label{Hoeven}
The stratification \begin{align*}
G/P_J=\mathop{\dot{\bigcup}}_{w\in \prescript{J}{}W}X_J\{w\}
\end{align*}
satisfies the following properties and is uniquely determined by them.
\begin{enumerate}
\item \label{Hoeven_1}If $\prescript{w}{}{F(J)}=J$ then $X_J\{w\}=X_J(w)$. 
\item \label{Hoeven_2}If $I\subset J$ and $w$ is in $\prescript{J}{}W\subset \prescript{I}{}W$ then the natural morphism $\alpha_{I,J}\colon G/P_I\to G/P_J$ gives a finite \' etale surjective morphism $X_I\{w\}\to X_J\{w\}$.
\item \label{Hoeven_3}If $w$ is in $\prescript{J}{}W$ and $I=J\cap \prescript{w'}{}{F(J)}$ where $w'$ is the representative in $\prescript{I}{}W^{F(I)}$ of $W_I wW_{F(I)}$ then the morphism $X_I\{w\}\to X_J\{w\}$ is an isomorphism. 
\end{enumerate}
\end{theorem}

\Cref{Hoeven} allows us to think of fine Deligne--Lusztig varieties in a more practical way, given in the following \cref{sneaky_definition} (see also \cite[Theorem 3.1]{he2008gstable}). 

\begin{cor}\label{sneaky_definition}
For $w\in \prescript{J}{}W$, the fine Deligne--Lusztig variety $X_J\{w\}$ is the scheme theoretic image of the morphism $X(w)=X_\varnothing\{w\}\to X_J\{w\}$ induced by the quotient map $X\to G/P_J$.  
\end{cor}
\begin{proof}
By \cref{Hoeven_2} of \Cref{Hoeven} we know that the morphism $X(w)\to X_J\{w\}$ is finite and thus quasi-compact \cite[Section 01WG]{stacksproject}.
We also know that it is surjective, and so we are done as $X_J\{w\}$ is a (dense) subset of the scheme theoretic image (see \cite[Lemma 01R8]{stacksproject}). 
\end{proof}

Fine Deligne--Lusztig varieties satisfy some of the same properties as the classical ones, as stated in the following \cref{prop_fine}.

\begin{proposition}\cite[Definition 2.1.6, Corollary 2.1.8]{hoeve}\label{prop_fine}
For each $w\in\prescript{J}{}W$, the scheme $X_J\{w\}$ is locally closed, smooth and of dimension $l(w)$. 
\end{proposition}

\begin{rem}
    We can verify that parabolic and fine Deligne--Lusztig varieties are different in general by considering their dimensions.
    This is exhibited in the following example. 
    Let $G^F$ be of type $A_2$, $J=\{s_1\}$ and $w=s_2$. 
    Then $\dim X_{s_1}(s_2)=2$ by \cite[I,2.1]{hoeve} whereas $\dim X_{s_1}\{s_2\}=1$ by \cref{prop_fine}. 
\end{rem}

\section{Integral Picard group}\label{sec_integral_picard}

In this section we are interested in describing the abelian group $\Pic(X(w))$ for various $w\in W$.
In \cref{subsec_longest} we consider the case $w=w_0$, and in \cref{sec:other_weyl_gp_elems} we restrict to $G^F$ of untwisted type $A$ to show that fine Deligne--Lusztig varieties form a $\Q$-basis for $\Pic(G/P(s))_{\Q}$ for any $s\in S$, which allows us to derive results for $\Pic(X(sw_0))$.

\subsection{The case of the longest Weyl group element}\label{subsec_longest}

Let $G^F$ be a finite group of Lie type as in \cref{sec:fin_gps_weyl_dynkin}, not of type $\prescript{2}{}{B_2}$, $\prescript{2}{}{G_2}$ or $\prescript{2}{}{F_4}$.
That is, $d=1$ and $F\colon G\to G$ is a geometric Frobenius.
Recall the longest element $w_0\in W$ from \cref{not_and_conv}.
In this section, we compute $\Pic(X(w_0))$.
By the following, we know that it is a torsion abelian group. 

\begin{proposition}[Hansen]\cite[Corollary 3.21]{hansen}\label{rat_w0_vanishes}
    The vector space $\CH_i(X(w_0))_{\Q}$ vanishes for all $i<l(w_0)$. 
\end{proposition}

Our results in this chapter will follow from the following beautiful theorem. 
Recall that the intersection product $\cdot$ makes $\CH^*(X)$ into a commutative, graded ring with unit $[X]$ (see \cite[Chapters 6, 8]{fulton2012intersection}). 

\begin{theorem}\cite[Theorem 3.10]{kim2016}\label{kim}
In $\CH_*(X)$ we have \begin{align*}
[\overline{X(w)}]=\sum_{\substack{u,v\in W\\ uwF(v)^{-1}=w_0\\ l(w_0)-l(w)=l(u)+l(v)}}([\Omega_{w_0u}]\cdot [\Omega_{w_0v}])q^{l(v)}.
\end{align*}
\end{theorem}

We start by stating a useful lemma.

\begin{lemma}\label{central_longest_2}
Let $W$ be a finite irreducible Weyl group.\begin{rqenumerate}\label{rqcentral_longest_2}
    \item \label{rqcentral_longest_211}There is an equality of sets \[\{sw_0:s\in S\}=\{w\in W:l(w)=l(w_0)-1\},\]and its cardinality is $|S|$. 
    Moreover, if $W$ is not of type $A_1$ then each element of the set above has full support. 
    \item \label{rqcentral_longest_21}If $W$ is of type $A_{n-1}$, then for $1\leq i\leq n-1$ we have $s_i=w_0s_{n-i}w_0$.
    \item \label{rqcentral_longest_22}The element $w_0$ is central in $W$ if and only if $W$ is of type $C_n$, $D_n (n\text{ even})$, $F_4$, $E_7$, $E_8$ or $G_2$.
    \item \label{rqcentral_longest_23}If $W$ is of type $D_n (n\text{ odd})$, then \begin{align*}
    w_0s_{n-1}w_0&=s_n,\\
    w_0s_iw_0&=s_i
    \end{align*}
    for $i<n-1$. 
    \item \label{rqcentral_longest_24}If $W$ is of type $E_6$ then \begin{align*}
        w_0s_1w_0&=s_6,\\
        w_0s_3w_0&=s_5,\\
        w_0s_2w_0&=s_2,\\
        w_0s_4w_0&=s_4.
\end{align*}
\end{rqenumerate}
\end{lemma}
\begin{proof}
    We start with \cref{rqcentral_longest_211}.
    The equality of sets and its cardinality follows from standard properties of $w_0$ (see for example \cite[Chapter 1, Section 1.4; Propositions 2.3.1, 2.3.2]{bjorner2006combinatorics}).
    The final claim follows from \cite[Chapter 5, Exercices, paragraph 6, Exercice 2]{bourbaki1968}, which says that for each $s\in S$ there exists a reduced expression of $w_0$ that contains $s$ at least twice. 
    \Cref{rqcentral_longest_21} is \cite[Planche I, (XI)]{bourbaki1968}.
    \Cref{rqcentral_longest_22} is precisely \cite[Exercise 1.13]{Geck2000CharactersOF}.
    \Cref{rqcentral_longest_23} and \cref{rqcentral_longest_24} follow from \cite[Planche I V, (XI)]{bourbaki1968} and \cite[Planche V, (XI)]{bourbaki1968} respectively.
\end{proof}

We also recall that the Chow groups of a flag variety (full or partial) are determined by Schubert varieties.

\begin{proposition}[Chevalley]\label{Chow_Flag}
	Let $J\subset S$, $P\coloneq P_J$ and $0\leq k\leq \dim G/P$.

    The Chow group $\CH_k(G/P)$ is freely generated by $[\Omega^P_w]$ where $w\in W^J$ such that $l(w)=k$. 
    Moreover, $G$ acts trivially on $\CH_k(G/P)$. 
\end{proposition}
\begin{proof}
    The full flag variety case is an unpublished result by Chevalley (see \cite[Section 3]{ASENS_1974_4_7_1_53_0}).
    For the more general case, this is \cite[Example 1.9.1]{fulton2012intersection} as Schubert varieties form a cellular decomposition of $G/P$.
    The last claim follows because $G$ is connected, see \cite[Expos\'e 5, Lemme 1]{Grothendieck1958}.
\end{proof}

\begin{proposition}\label{DL_in_terms_of_Sch_2}
Let $s\in S$, then \Cref{table:1} expresses the cycle $[\overline{X(sw_0)}]$ in $\Pic(X)$ in terms of Schubert divisor cycles (see \cref{Chow_Flag}).
\end{proposition}

\begin{table}[h!]
\centering
\begin{tabular}{ m{10em} | m{15em}} 
 \makecell{Type of $W$} & \makecell{$[\overline{X(sw_0)}]$} \\ [0.5ex] 
 \hline\hline
 \\[-1em]
 \makecell{$A_{n}$} & \makecell{$q[\Omega_{sw_0}]+[\Omega_{w_0s}]$} \\
 \\[-1em]
 \hline 
 \\[-1em]
 \makecell{$\prescript{2}{}{A_{n}}$, $C_n$, $D_n (n\text{ even})$,\\
$\prescript{2}{}{D_n} (n\text{ odd})$, $F_4$, $E_7$,\\
$E_8$, $\prescript{2}{}{E_6}$, $G_2$} & \makecell{$(q+1)[\Omega_{w_0s}]$}  \\
 \\[-1em]
 \hline 
 \\[-1em]
 \makecell{$D_n (n\text{ odd})$,\\
 $\prescript{2}{}{D_n} (n\text{ even})$} & \makecell{$[\Omega_{w_0s_n}]+[\Omega_{w_0s_{n-1}}]q$ if $s=s_n$,\\
 $[\Omega_{w_0s_n}]q+[\Omega_{w_0s_{n-1}}]$ if $s=s_{n-1}$,\\
 $(q+1)[\Omega_{w_0s}]$ else}\\
 \\[-1em]
 \hline 
 \\[-1em]
 \makecell{$\prescript{3}{}{D_4}$} & \makecell{$[\Omega_{w_0s_1}]+[\Omega_{w_0s_3}]q$ if $s=s_1$,\\
 $[\Omega_{w_0s_3}]+[\Omega_{w_0s_4}]q$ if $s=s_3$,\\
 $[\Omega_{w_0s_4}]+[\Omega_{w_0s_1}]q$ if $s=s_4$,\\
 $(q+1)[\Omega_{w_0s_2}]$ if $s=s_2$} \\
 \\[-1em]
 \hline 
 \\[-1em]
 \makecell{$E_6$} & \makecell{
$[\Omega_{w_0s_1}]+[\Omega_{w_0s_6}]q$ if $s=s_1$,\\
$[\Omega_{w_0s_6}]+[\Omega_{w_0s_1}]q$ if $s=s_6$,\\
$[\Omega_{w_0s_3}]+[\Omega_{w_0s_5}]q$ if $s=s_3$,\\
$[\Omega_{w_0s_5}]+[\Omega_{w_0s_3}]q$ if $s=s_5$,\\
$(q+1)[\Omega_{w_0s_2}]$ else}
\end{tabular}
\caption{$[\overline{X(sw_0)}]$ in terms of Schubert divisor cycles.}
\label{table:1}
\end{table}

\begin{proof}
Let $s\in S$.
Note that for $u,v\in W$ with $l(u)+l(v)=1$ (which is equal to $l(w_0)-l(sw_0)$), it follows that $l(u)\neq l(v)$ and $0\leq l(u),l(v)\leq 1$.
In particular, $u$ and $v$ are involutions.
We have the following equivalence\begin{align*}
usw_0F(v)=w_0\Leftrightarrow
sw_0=uw_0F(v),
\end{align*}
hence the terms in the sum in \Cref{kim} correspond to two pairs of options $(u,v)=(s,e)$ and $(u,v)=(e,F^{-1}(w_0sw_0))$. 
We claim that the equalities now follow from \cref{central_longest_2}.

Indeed, if $F$ acts trivially on $W$ then this follows directly.
If $W$ has type $\prescript{2}{}{A_{n}}$, then the Frobenius switches $s_i$ and $s_{n+1-i}$.
If $W$ has type $\prescript{2}{}{D_n}$, the Frobenius fixes $s_i$ for $i\geq 2$ and permutes $s_{n-1}$ and $s_n$.
Hence, if $W$ is of type $\prescript{2}{}{A_{n}}$ or $\prescript{2}{}{D_n} (n\text{ odd})$, then $v=s$. 
Similarly, if $W$ is of type $D_n (n\text{ odd})$ or $\prescript{2}{}{D_n} (n\text{ even})$ then $v=s$ if $s=s_i$ with $i\leq n-2$ and $v=s_{n-1}$ if $s=s_n$ and $v=s_n$ if $s=s_{n-1}$. 
The remaining $\prescript{3}{}{D_4}$ and $\prescript{2}{}{E_6}$ cases can be checked by hand.
\end{proof}

Now we are ready to compute $\Pic(X(w_0))$.
Note that if $G^F$ is of type $A_1$, then $w_0$ is Coxeter and in that case $X(w_0)\cong \Pro^1\setminus \Pro^1(\F_q)$ is the \textit{Drinfeld upper half space} of dimension $1$ (see \cite[2.2]{originalpaper}).
In particular, we already know $\Pic(X(w_0))=0$ (see for example \cite[Theorem 5]{hansen2002picard}). 
This is why we exclude it from \cref{Picard_longest_elem}.

\begin{cor}\label{Picard_longest_elem}
Suppose $G^F$ is not of type $A_1$.
An explicit description of $\Pic(X(w_0))$ is given in \Cref{table:2}.
Moreover, the action of $G^F$ on $\Pic(X(w_0))$ is trivial. 
\end{cor}

\begin{table}[H]
\centering
\begin{tabular}{ m{10em} | m{15em}} 
 \makecell{Type of $W$} & \makecell{$\Pic(X(w_0))$} \\ [0.5ex] 
 \hline\hline
 \\[-1em]
 \makecell{$A_{2m}$} & \makecell{$(\Z/(q^2-1)\Z)^{m}$} \\
 \\[-1em]
 \hline 
 \\[-1em]
 \makecell{$A_{2m+1}$, $m>0$} & \makecell{$(\Z/(q^2-1)\Z)^{m}\oplus \Z/(q+1)\Z$}  \\
 \\[-1em]
 \hline 
 \\[-1em]
 \makecell{$\prescript{2}{}{A_{n}}$, $C_n$, $D_n (n\text{ even})$, \\
$\prescript{2}{}{D_n} (n\text{ odd})$, $F_4$, $E_7$,\\ $E_8$, $G_2$, $\prescript{2}{}{E_6}$} & \makecell{$(\Z/(q+1)\Z)^m$ where $m\coloneq |S|$}\\
 \\[-1em]
 \hline 
 \\[-1em]
 \makecell{$D_n (n\text{ odd})$,\\
$\prescript{2}{}{D_n} (n\text{ even})$} & \makecell{$(\Z/(q+1)\Z)^{n-2}\oplus \Z/(q^2-1)\Z$} \\
 \\[-1em]
 \hline 
 \\[-1em]
 \makecell{$\prescript{3}{}{D_4}$} & \makecell{$\Z/(q+1)\Z\oplus \Z/(q^3+1)\Z$} \\
\\[-1em]
 \hline 
 \\[-1em]
 \makecell{$E_6$} & \makecell{$(\Z/(q+1)\Z)^{2}\oplus (\Z/(q^2-1)\Z)^2$} 
\end{tabular}
\caption{Explicit description of $\Pic(X(w_0))$.}
\label{table:2}
\end{table} 

\begin{proof}
    Note $\overline{X(w_0)}=X$ and \begin{align*}
        X\setminus X(w_0)=\bigcup_{s\in S}\overline{X(sw_0)}
    \end{align*}is of codimension $1$ in $X$.
    Moreover, assuming $W$ is not of type $A_1$, it consists of $|S|$ irreducible components by \cref{rqcentral_longest_211} of \cref{central_longest_2} and \cref{rqgen_iso_normal4} of \cref{gen_iso_normal}. 
    Hence, we have a localization (exact) sequence (\cite[Proposition 1.8]{fulton2012intersection})
    \[
        \begin{tikzcd}
        \displaystyle \bigoplus_{s\in S}\CH^0(\overline{X(sw_0)}) \arrow[r] & \Pic(X) \arrow[r] & \Pic(X(w_0)) \arrow[r] & 0.
        \end{tikzcd}
    \]
    Recall that $\Pic(X)$ is freely generated by $[\Omega_{w_0s}]$ for $s\in S$ (see \Cref{Chow_Flag}).
    The presentation of $\Pic(X(w_0))$ follows from \Cref{DL_in_terms_of_Sch_2}, the exact sequence above and the following claim. 
    There are isomorphisms of $\Z$-modules
    \begin{align*}
        \dfrac{\Z x\oplus \Z y}{(x+qy,qx+y)}&\cong\Z/(q^2-1)\Z,\\        
        \dfrac{\Z x\oplus \Z y\oplus \Z z}{(x+qy,y+qz,z+qx)}&\cong\Z/(q^3+1)\Z.
    \end{align*}
    For the first isomorphism, one notices the relation $y=-qx$ in the quotient, and hence 
    \begin{align*}
        \dfrac{\Z x\oplus \Z y}{(x+qy,qx+y)}\cong\dfrac{\Z x}{(x-q^2x)}\cong \Z/(q^2-1)\Z.
    \end{align*}
    The second isomorphism is obtained similarly. 
    
    To show that $G^F$ acts trivially, let $g\in G^F$. 
    The open embedding $X(w_0)\subset X$ induces a $G^F$-equivariant map $\Pic(X)\to \Pic(X(w_0))$ which sits in a commutative diagram 
    \[
        \begin{tikzcd}
        \Pic(X) \arrow[r] \arrow[d, "g"] & \Pic(X(w_0)) \arrow[d, "g"] \arrow[r] & 0 \\
        \Pic(X) \arrow[r]                & \Pic(X(w_0)) \arrow[r]                & 0
        \end{tikzcd}
    \]
    where the rows are exact.
    By \Cref{Chow_Flag}, we know that the map $g\colon \Pic(X)\to \Pic(X)$ is the identity. 
    Hence, the same is true for $g\colon \Pic(X(w_0))\to \Pic(X(w_0))$ by the surjectivity of $\Pic(X)\to \Pic(X(w_0))$. 
\end{proof}

\begin{rem}
    Notice that $G^F$ acts trivially on $\Pic(X(w_0))$ even if $G^F$ is of type $\prescript{2}{}{B_2}$, $\prescript{2}{}{G_2}$ or $\prescript{2}{}{F_4}$; the same argument as in the proof of \cref{Picard_longest_elem} goes through.
\end{rem}

\subsection{Other Weyl group elements}\label{sec:other_weyl_gp_elems}

For this section, suppose $G^F$ is of untwisted type $A_{n-1}$ for $n\geq 3$. 
For some suitable $s\in S$, combining \Cref{kim} with the theory of fine Deligne--Lusztig varieties allows us to show the vanishing of the rational Picard group of $X(sw_0)$ (see \Cref{rat_Pic_last_to_longest}). 
However, some Schubert calculus is required. 
The techniques even allow us to compute this Picard group integrally, which we illustrate with an example in the case of $n=4$ in \Cref{GL_4_last_to_longest}. 
The statement and proof of \cref{fine_DL_as_parab_Sch} are technical, and therefore the reader may prefer to read \cref{GL_4_last_to_longest} first to gain intuition.

We first recall some facts that we will need from Schubert calculus.
For $1\leq k\leq n$, define the following polynomials in $\Z[x_1,\ldots, x_n]$ called the \textit{elementary symmetric functions}
\begin{align*}
    e_k\coloneq \displaystyle \sum_{1\leq i_1<\cdots<i_k\leq n}x_{i_1}\cdots x_{i_k}.
\end{align*} 
The Chow ring of the full flag variety $X$ has the following description. 

\begin{proposition}[Borel]\label{Sch_var_poly_iso}
There is an isomorphism of rings \begin{align*}
\CH_*(X)&\to \dfrac{\Z[x_1,\ldots, x_n]}{(e_1,\ldots,e_n)}.
\end{align*}
\end{proposition}
\begin{proof}
See for example \cite[Proposition 3.6.15]{manivel2001symmetric}.
\end{proof}

\begin{definition}
For $w\in W$, the \textit{Schubert polynomial associated to} $w$ is the image of $[\Omega_{w_0w}]$ under the isomorphism in \Cref{Sch_var_poly_iso} (see \cite[Theorem 3.6.18]{manivel2001symmetric}).
It is denoted by $\mathfrak{S}_w$.
\end{definition}

The following is often referred to as Monk's formula. 

\begin{theorem}\cite[Theorem 2.7.1]{manivel2001symmetric}
For all $w\in W$ and $1\leq m\leq n-1$ we have \begin{align*}
\mathfrak{S}_w\cdot \mathfrak{S}_{s_m}=\sum_{\substack{j\leq m< k\\ l(w(j\ k))=l(w)+1}}\mathfrak{S}_{w(j\ k)}. 
\end{align*}
\end{theorem}

As an immediate consequence, this gives us the following. 

\begin{lemma}\label{monk_cons}
Let $1\leq i,j\leq n-1$, then 
\begin{align*}
\mathfrak{S}_{s_i}\cdot \mathfrak{S}_{s_j}&=\begin{cases}
	\mathfrak{S}_{s_{i+1}s_{i}}&\text{if }i=j\neq n-1\\
\mathfrak{S}_{s_{i-1}s_{i}}&\text{if }i=j=n-1\\
\mathfrak{S}_{s_{i}s_{j}}+\mathfrak{S}_{s_{j}s_{i}}&\text{if }s_is_j\neq s_js_i\\\mathfrak{S}_{s_is_j}&\text{if }s_is_j=s_js_i\text{ and }i\neq j.
\end{cases}
\end{align*}
\end{lemma}

We can now determine fine Deligne--Lusztig cycles in the Picard group of a certain partial flag variety in terms of Schubert cycles.

\begin{proposition}\label{fine_DL_as_parab_Sch}
Let $1\leq i\neq j\leq n-1$ and define for $I\coloneq \{s_i\}$ the following
\begin{alignat*}{2}
w_j&\coloneq w^i_j&&\coloneq s_is_jw_0\in \prescript{I}{}W\\
x_j&\coloneq x^i_j&&\coloneq [\overline{X_I\{w_j\}}]\\
e_j&\coloneq e^i_j&&\coloneq [\Omega_{w_j^{-1}}^{P_I}]
\end{alignat*}where the cycles are seen as elements of $\Pic(G/P_I)$. 
Then, up to potential $\Q$-scaling in the case where $2i=n$, we have 
\begin{align*}
    x_j=\begin{cases}
    (q^2+1)e_j+qe_{n-j-1}&\text{if }i=n-1\text{ and } j=1\\
    (q^2+1)e_j+qe_{n-j+1}&\text{if }i\neq n-1\text{ and }i=n-j\\
    (q+1)e_j+(q+q^2)e_{n-j}&\text{if }s_is_j=s_js_i\text{ and }i=n-i\\
    e_j+qe_{n-j}&\text{else}.
    \end{cases}
\end{align*}
\end{proposition}
\begin{proof}
Using \Cref{kim}, we get \begin{equation}\label{eq:DL_in_terms_Schproof}
[\overline{X(w_j)}]=\begin{cases}
[\Omega_{w_j^{-1}}]+[\Omega_{w_0s_i}]\cdot [\Omega_{w_0s_{n-j}}]q+[\Omega_{w_0s_{n-i}s_{n-j}}]q^2&\text{if }s_is_j\neq s_js_i\\
[\Omega_{w_j^{-1}}]+[\Omega_{w_0s_i}]\cdot [\Omega_{w_0s_{n-j}}]q+[\Omega_{w_0s_j}]\cdot [\Omega_{w_0s_{n-i}}]q+[\Omega_{w_0s_{n-i}s_{n-j}}]q^2&\text{else}.
\end{cases}
\end{equation}
The tactic of the proof is as follows. 
We want to push $[\overline{X(w_j)}]$ and the Schubert cycles in \cref{eq:DL_in_terms_Schproof} forward under $\CH_{l(w_0)-2}(X)\to \CH_{l(w_0)-2}(G/P_I)=\Pic(G/P_I)$.
In order to do this for the Schubert cycles, we need to compute the intersection products in \cref{eq:DL_in_terms_Schproof} first.
Computing the intersection products in the formula comes down to computing $\mathfrak{S}_{s}\cdot \mathfrak{S}_{t}$ for particular $s,t\in S$ as specified above (see \Cref{Sch_var_poly_iso}).
To compute this product of Schubert polynomials, we use \cref{monk_cons}.
Depending on the Weyl group element, the Schubert cycles are either mapped to $0$ or to a unique Schubert cycle in $\Pic(G/P_I)$ by \cref{Sch_birat}.

We start with computing the pushforward of $[\overline{X(w_j)}]$.
If $2i\neq n$, then $[\overline{X(w_j)}]$ maps to $x_j$ as $X(w_j)\cong X_I\{w_j\}$ by \cref{Hoeven_3} of \cref{Hoeven}.
Indeed,  
\begin{align*}
    w_js_iw_j^{-1}&=s_is_js_{n-i}s_js_i
\end{align*}
and this is equal to $s_i$ if and only if $s_js_{n-i}s_j=s_i$.
Regardless of whether $s_j$ and $s_{n-i}$ commute, the reader may check that the latter is impossible (we obtain contradictions with $2i\neq n$ if they commute and with $l(s_js_{n-i}s_j)=3$ if they do not).
If $2i=n$ then $X(w_j)$ and $X_I\{w_j\}$ may not be isomorphic and hence $[\overline{X(w_j)}]\mapsto dx_j$ for some $d\in \N\setminus \{0\}$ by \cref{Hoeven_2} of \cref{Hoeven}.
Hence, up to potential $\Q$-scaling, we may assume $[\overline{X(w_j)}]\mapsto x_j$.

We now move on to computing the intersection products in \cref{eq:DL_in_terms_Schproof}, and then we apply the pushforward.
As stated in \Cref{monk_cons}, we need to distinguish four cases depending on $k\coloneq |i-(n-j)|=|(n-i)-j|$ as follows \begin{enumerate}
\item \label{case1_proof}$k=0$ and $i\neq n-1$,
\item $k=0$ and $i= n-1$,
\item $k=1$,
\item $k>1$.
\end{enumerate}
We illustrate the case in \cref{case1_proof}, as the other cases are similar.
Suppose $k=0$ and $i\neq n-1$, which in particular means $i=n-j$. 
If $s_is_j\neq s_js_i$ then by \cref{eq:DL_in_terms_Schproof}, \cref{monk_cons} and \cref{Sch_var_poly_iso}
\begin{align*}
    [\overline{X(w_j)}]&=[\Omega_{w_j^{-1}}]+[\Omega_{w_0s_i}]\cdot [\Omega_{w_0s_{n-j}}]q+[\Omega_{w_0s_{n-i}s_{n-j}}]q^2\\
    &=[\Omega_{w_j^{-1}}]+[\Omega_{w_0s_{i+1}s_i}]q+[\Omega_{w_0s_{n-i}s_{n-j}}]q^2.
\end{align*}
Note 
\begin{align*}
    w_j^{-1}&=w_0s_js_i=w_0s_{n-i}s_{n-j}\in W^I,\\
    w_{n-j+1}^{-1}&=(s_is_{n-j+1}w_0)^{-1}=w_0s_{n-j+1}s_i=w_0s_{i+1}s_i\in W^I
\end{align*}
and hence by \cref{Sch_birat} we obtain
\begin{align*}
    x_j=(q^2+1)e_j+qe_{n-j+1}
\end{align*}
after applying the pushforward.
If $s_is_j= s_js_i$ then we have to compute one more intersection product, namely $[\Omega_{w_0s_j}]\cdot [\Omega_{w_0s_{n-i}}]=[\Omega_{w_0s_j}]\cdot [\Omega_{w_0s_{j}}]$.
By \cref{monk_cons}, this product is either $[\Omega_{w_0s_{j+1}s_j}]$ or $[\Omega_{w_0s_{j-1}s_j}]$. 
As $i\neq j$, we have $w_0s_{j+1}s_j, w_0s_{j-1}s_j\notin W^I$ because multiplying by $w_0$ on the left is an antiautomorphism on $(W,\leq)$ (\cite[Proposition 2.3.4]{bjorner2006combinatorics}).
Hence, the intersection product $[\Omega_{w_0s_j}]\cdot [\Omega_{w_0s_{j}}]$ is mapped to $0$.
Therefore, we again get 
\begin{align*}
    x_j=(q^2+1)e_j+qe_{n-j+1}.
\end{align*}

Applying the same tactic for the other three cases allows us to obtain the following. 
Suppose first $s_is_j\neq s_js_i$.
Then \begin{align*}
x_j=\begin{cases}
(q^2+1)e_j+qe_{n-j+1}&\text{if }k=0\text{ and }i\neq n-1\\
(q^2+1)e_j+qe_{n-j-1}&\text{if }k=0\text{ and }i=n-1\\
e_j+qe_{n-j}&\text{else}.
\end{cases}
\end{align*}
If $s_is_j=s_js_i$, then 
\begin{align*}
x_j=\begin{cases}
(q^2+1)e_j+qe_{n-j+1}&\text{if }k=0\text{ and }i\neq n-1\\
(q^2+1)e_j+qe_{n-j-1}&\text{if }k=0\text{ and }i=n-1\\
(q+1)e_j+(q+q^2)e_{n-j}&\text{if }k\geq 1\text{ and }i=n-i\\
e_j+qe_{n-j}&\text{else}.
\end{cases}
\end{align*}
Summarizing gives us the claimed expressions. 
\end{proof}

In \cite[Proposition 3.19]{hansen}, Hansen showed that (classical) Deligne--Lusztig cycles form a basis for $\CH_*(X)_{\Q}$ (see also \cite[Theorem 3.11]{kim2016} for a stronger result).
We are now in a position to show an analogous statement in codimension $1$ for the partial flag varieties $G/P(s)$ for $s\in S$ and in the case where $G^F$ is of type $A_{n-1}$.

\begin{cor}\label{basis_partial}
Keep using the same notation as \Cref{fine_DL_as_parab_Sch}.
For $1\leq i\leq n-1$, the set $\{x_j:1\leq j\neq i\leq n-1\}$ forms a $\Q$-basis for $\Pic(G/P(s_i))_{\Q}$. 
\end{cor}
\begin{proof}
As $\Pic(G/P(s_i))$ is freely generated by the $e_j$ (see \Cref{Chow_Flag}), it is enough to write each $e_j$ as a $\Q$-linear combination of the $x_j$ and to show that the $x_j$ form a $\Q$-linearly independent set.
To do this, we make a case distinction. 

Suppose $2i\neq n$.
Let $j$ such that $n-j\neq i$, then \begin{align*}
x_j&=e_j+qe_{n-j}\\
x_{n-j}&=e_{n-j}+qe_{j}
\end{align*}which means that we can write $e_j$ and $e_{n-j}$ as a $\Q$-linear combination of $x_j$ and $x_{n-j}$, which are $\Q$-linearly independent. 
We will be done with this case when we show that we can also obtain $e_{n-i}$ in this way. 
We have \begin{align*}
x_{n-i}=\begin{cases}
(q^2+1)e_{1}+qe_{n-2}&\text{if }i=n-1\\
(q^2+1)e_{n-i}+qe_{i+1}&\text{else}.
\end{cases}
\end{align*}
If $i=n-1$, then we have shown above that $e_{n-2}$ can be written as a $\Q$-linear combination of $x_2$ and $x_{n-2}$. 
If $i\neq n-1$, then we need to distinguish if $n-(i+1)$ equals $i$ or not. 
If $n-(i+1)=i$ then $x_{n-i}=(q^2+q+1)e_{n-i}$.
If not, then we already know that $e_{i+1}$ can be written as a $\Q$-linear combination of $x_{i+1}$ and $x_{n-(i+1)}$, which implies that we obtain $e_{n-i}$ in the desired way. 

Now suppose $2i=n$, then we may work up to potential $\Q$-scaling. 
Let $j$ such that $s_is_j=s_js_i$.
In particular, $s_is_{n-j}=s_{n-j}s_i$ as $1<|i-j|=|n-i-(n-j)|=|i-(n-j)|$ and so \begin{align*}
x_j&=(q+1)e_j+(q^2+q)e_{n-j}\\
x_{n-j}&=(q+1)e_{n-j}+(q^2+q)e_{j}
\end{align*}
which means that we can write $e_j$ and $e_{n-j}$ as a $\Q$-linear combination of $x_j$ and $x_{n-j}$, which are $\Q$-linearly independent.
Now let $j$ be such that $s_is_j\neq s_js_i$ (there are precisely two such options; $j=i+1,i-1$), then \begin{align*}
x_j&=e_j+qe_{n-j}\\
x_{n-j}&=e_{n-j}+qe_{j}
\end{align*} which implies that we obtain the two remaining vectors $e_{i+1}$ and $e_{i-1}$ in the desired way.
\end{proof}

\begin{cor}\label{rat_Pic_last_to_longest}
Let $1\leq i\leq n$ such that $2i\neq n$. 
Then $\Pic(X(s_iw_0))_{\Q}=0$. 
Moreover, the $G^F$-action on $\Pic(X(s_iw_0))$ is trivial. 
\end{cor}
\begin{proof}
Set $w\coloneq s_iw_0\in \prescript{\{s_i\}}{}W$.
By \cref{Hoeven_3} of \Cref{Hoeven}, $X(w)\cong X_{s_i}\{w\}$ and so it is equivalent to compute $\Pic(X_{s_i}\{w\})$.
As $X_{s_i}\{w\}$ is open in $G/P(s_i)$, and \[G/P(s_i)\setminus X_{s_i}\{w\}=\bigcup^{n-1}_{\substack{j=1\\j\neq i}}\overline{X_{s_i}\{w_j\}}\](where $w_j$ is defined as in \cref{fine_DL_as_parab_Sch}) we have an exact sequence 
\[
    \begin{tikzcd}
    \displaystyle\bigoplus^{n-1}_{\substack{j=1\\j\neq i}}\CH^0(\overline{X_{s_i}\{w_j\}}) \arrow[r,"f"] & \Pic(G/P(s_i)) \arrow[r] & \Pic(X_{s_i}\{w\}) \arrow[r] & 0.
    \end{tikzcd}
\]
Tensoring this sequence with $\Q$ produces an exact sequence of $\Q$-vector spaces, of which the pushforward $f$ is surjective by \Cref{basis_partial}.
The $G^F$-action being trivial on $\Pic(X_{s_i}\{w\})$ follows from \Cref{Chow_Flag}, in the same way as in the proof of \Cref{Picard_longest_elem}.
\end{proof}

\begin{rem}\label{rem:DL_iso_fineDL}
Combining the proof of \Cref{rat_Pic_last_to_longest} and \Cref{fine_DL_as_parab_Sch} allows us to even compute the Picard group of  $X(s_iw_0)\cong X_{s_i}\{s_iw_0\}$ integrally for $1\leq i\leq n-1$ with $2i\neq n$.
We treat the $A_3$ case in \Cref{GL_4_last_to_longest}.
Moreover, the condition $2i\neq n$ in \Cref{rat_Pic_last_to_longest} is equivalent to $X(s_iw_0)\cong X_{s_i}\{s_iw_0\}$.
Indeed, by \cite[Example 3.10(b)]{alma994473953502466} we know that the quotient map $X\to G/P(s_i)$ restricts to a $\Pro^1(\F_q)$-bundle $X(s_iw_0)\to X_{s_i}(w_0)=X_{s_i}(s_iw_0)=X_{s_i}\{s_iw_0\}$. 
This means that a similar tactic cannot be used to compute the integral Picard group of $X(s_iw_0)$ if $2i=n$.  
\end{rem}

\begin{exmp}\label{GL_4_last_to_longest}
Suppose $G^F$ is of type $A_3$.
We will show that the integral Picard groups of $X(s_1w_0)$ and $X(s_3w_0)$ are isomorphic to \begin{align*}
\Z/(q^2+1)(q+1)\Z.
\end{align*}
As the varieties $X(s_1w_0)$ and $X(s_3w_0)$ are isomorphic to $X_{s_1}\{s_1w_0\}$ and $X_{s_3}\{s_3w_0\}$ respectively (see proof of \cref{rat_Pic_last_to_longest}), it is equivalent to compute the integral Picard groups of the latter.
We are done by \Cref{Chow_Flag} and the proof of \Cref{rat_Pic_last_to_longest} as \begin{align*}
x_2^1&=(q+1)e_2,\\
x_3^1&=(q^2+1)e_3+qe_2,\\
x_1^3&=(q^2+1)e_1+qe_2,\\
x_2^3&=(q+1)e_2,
\end{align*}by \Cref{fine_DL_as_parab_Sch}, and we have an isomorphism of $\Z$-modules \begin{align*}
    \dfrac{\Z x\oplus \Z y}{((q^2+1)x+qy,(q+1)y)}\cong \Z/(q^2+1)(q+1)\Z.
\end{align*}
Indeed, as $-y=qy=-(q^2+1)x$ we have 
\begin{align*}
    \dfrac{\Z x\oplus \Z y}{((q^2+1)x+qy,(q+1)y)}&\cong \dfrac{\Z x}{((q^2+1)(q+1)x)}\\
    &\cong \Z/(q^2+1)(q+1)\Z.
\end{align*}
We can also compute $\Pic(X_{s_2}\{s_2w_0\})\cong\Z/(q^2-1)\Z$, as 
\begin{align*}
    x_1^2&=e_1+qe_3,\\
    x_3^2&=e_3+qe_1
\end{align*}
by \cref{fine_DL_as_parab_Sch} (note that these equalities are true without needing to scale).
\end{exmp}

\section{Rational Picard group}\label{sec_rat_Picard_group}

Let $G^F$ be a finite group of Lie type as in \cref{sec:fin_gps_weyl_dynkin}.
As mentioned in the \nameref{Introduction}, Hansen conjectured\footnote{Recall the footnote for \Cref{conjintroduction}.} the following (\cite[Conjecture 1]{hansen}).

\begin{conj}[Hansen]\label{conj}
    Let $w\in W$ be $F$-Coxeter and $i\in \N$ such that $i<l(w)$.
    Then $\CH_i(X(w))$ has rank zero, i.e., $\CH_i(X(w))_{\Q}=0$. 
\end{conj} 

As field extensions are faithfully flat, we obtain an equivalent conjecture if we tensor with $\Ql$ instead of $\Q$.
In this section, we will see that there are $w\in W$ such that $\Pic(X(w))$ has a free part as an abelian group.
More precisely, $\Pic(X(w))_{\Ql}\coloneq \Pic(X(w))\otimes_{\Z}\Ql$ is a (nonzero) $G^F$-representation over $\Ql$ which we will henceforth refer to as the rational Picard group. 
We study these representations for all finite groups of Lie type, and compare them to the corresponding $\ell$-adic cohomology groups $H^2(X(w))$ via cycle maps.
In order to do this, we first establish a connection (\cref{cycle_map_lemma_app}) between $\ell$-adic cohomology of a variety and its first higher Chow group in a particular degree in \cref{sec:motivic_cohom}.
We exploit this connection by invoking known vanishing results for $\ell$-adic cohomology groups of Deligne--Lusztig varieties in the untwisted type $A$ case (\cref{orlik_codim_11}), and proving related vanishing results for any type (\cref{H1_vanishing_other_gps}). 
Along the way, we will see how a result by Lusztig (\cite[Corollary 2.7]{Lusztig1976}) implies that \cref{conj} holds for the principal series part of $\CH^k(X(w))_{\Ql}$ for any type and $w\in W$ an $F$-Coxeter element.

\subsection{Motivic cohomology}\label{sec:motivic_cohom}

Let $Y$ be a smooth separated scheme of finite type over a perfect field $k$.
In \cite{Mazza2006LectureNO}, the motivic cohomology groups $H^{n,i}_\mot(Y,\Z)$ are defined as the hypercohomology of motivic complexes $\Z(i)$ of sheaves with respect to the Zariski topology (see \cite[Definition 3.1, Lemma 3.2]{Mazza2006LectureNO}).
It turns out that these groups are merely a renumbering of Bloch's higher Chow groups $\text{CH}^m(Y,l)$ (known as the Fundamental Comparison Theorem, see \cite[Theorem 19.1]{Mazza2006LectureNO}), which constitute a way of extending the basic exact sequence of Chow groups to a long exact sequence (see \cite{BLOCH1986267}).
Indeed, for all $n$ and $i\geq 0$ there is a natural isomorphism $H^{n,i}_\mot(Y,\Z)\cong \CH^i(Y,2i-n)$.
Moreover, we also have a natural isomorphism $\CH^1(Y,1)\cong \mO_Y(Y)^\times$ (\cite[Theorem 6.1]{BLOCH1986267}).
We will often identify these abelian groups.
We also use the Beilinson--Lichtenbaum conjecture, which is now a theorem, stated as follows. 

\begin{theorem}[Suslin--Voevodsky, Geisser--Levine]\label{beilinson_lichtenbaum}
Let $m$ be a positive integer that is invertible in $k$. 
Then the cycle map \begin{align*}
H^{i,j}_\mot(Y,\Z/m)\to H^i_{\et}(Y,\mu_m^{\otimes j})
\end{align*}is an isomorphism if $i\leq j$, and injective if $i=j+1$. 
\end{theorem} 
\begin{proof}
The cycle map $H^{i,j}_\mot(Y,\Z/m)\to H^i_\Lich(Y,\Z/m(j))$ to Lichtenbaum motivic cohomology (also called \'etale motivic cohomology, \cite[Definition 10.1]{Mazza2006LectureNO}) is an isomorphism if $i\leq j$ and injective if $i=j+1$ (see for example \cite[Theorem 6.17]{voevodsky_23} and \cite[Corollary 1.2]{Levine}).
Now apply \cite[Theorem 10.2]{Mazza2006LectureNO} to obtain a composition of isomorphisms.
\end{proof}

For an abelian group $N$, denote by $T(N)$ its torsion subgroup (which is in particular normal) and define $N/\text{tors}\coloneq N/T(N)$.
Recall from \cref{sec:fin_gps_weyl_dynkin} that $\F$ denotes an algebraic closure of a finite field.

\begin{proposition}\label{fin_gen}
Let $Y$ be a smooth projective variety over $\F$, and $U\subset Y$ an open subvariety such that the complement is of pure codimension $1$. 
Let $N\coloneq \mO_{U}(U)^\times$, then $N/\text{tors}$ is a finitely generated free abelian group. 
\end{proposition}
\begin{proof}
As $U$ is open in $Y$, we obtain a long exact sequence of higher Chow groups (see \cite[Theorem 3.1]{BLOCH1986267})
\[
    \begin{tikzcd}
    \cdots \arrow[r] & \mO_{Y}(Y)^\times \arrow[r,"g"] & N \arrow[r,"f"]          & \CH^0(Y\setminus U) \arrow[lld] \\
                    & \Pic(Y) \arrow[r]                         & \Pic(U) \arrow[r] & 0.                                             
    \end{tikzcd}
\]
We have that $\mO_{Y}(Y)^\times\cong \F^\times$ by \cite[Theorem 3.4, page 18]{hartshorne}, which means that every element of $\mO_{Y}(Y)^\times$ is torsion. 
On the other hand, $Y$ has finitely many irreducible components by \cite[Lemma 0BA8]{stacksproject} as it is a variety over $\F$.
This means that also $Y\setminus U$ consists of finitely many irreducible components.
Hence, $\CH^0(Y\setminus U)$ is a finitely generated free abelian group (see \cite[Example 1.3.2]{fulton2012intersection}). 
We now show $T(N)=\ker(f)$.
Let $x\in \ker f=\im g$, this means that there exists a $a\in \mO_Y(Y)^\times$ such that $g(a)=x$, and a $m\in\N$ such that $a^m=1$. 
In particular, $x^m=g(a)^m=1$ and so $\ker(f)\subset T(N)$.
The other inclusion is immediate as $\CH^0(Y\setminus U)$ is torsion-free.
Hence,\begin{align*}
N/\text{tors}=N/\ker(f)\cong \im(f)\subset \CH^0(Y\setminus U).
\end{align*}
As $\Z$ is a noetherian ring, we have that $N/\text{tors}$ is finitely generated as a $\Z$-module.
From the structure theorem for finitely generated abelian groups (see \cite[5.2.6]{stillwell2012classical}), it follows that $N/\text{tors}$ is also free.
\end{proof}

\begin{rem}\label{app_rem}
In particular, \cref{fin_gen} holds in the following two situations. 
\begin{enumerate}
\item \label{app_rem_1} For $U\coloneq X(w)$ and $Y\coloneq \overbar{X}(w)$ for $w\in F^+$,
\item \label{app_rem_2} For $U\coloneq X(\underline tv)$ and $Y\coloneq \overbar{X}(t,s_{i_1},\ldots,s_{i_r})$ where $t\in S$ and $v=s_{i_1}\cdots s_{i_r}\in F^+$ (see \cref{partial_compactification}). 
\end{enumerate}
We illustrate \cref{app_rem_1}.
If $w=s_{i_1}\cdots s_{i_r}\in F^+$ then 
\begin{align*}
    \partial \overbar{X}(w)=\bigcup_{j=1}^r\overbar{X}(w_j)
\end{align*}
where 
\begin{align*}
    w_j=s_{i_1}\cdots \widehat{s_{i_j}}\cdots s_{i_r}
\end{align*}
is a word of length $r-1=l(w)-1$ in the monoid $F^+$.
By \cite[Example 1.3.1(c), Proposition 1.8]{fulton2012intersection}, we get \begin{align*}\CH_{l(w)-1}(\partial \overbar{X}(w))&\cong \bigoplus_{j=1}^r\CH_{l(w)-1}(\overbar{X}(w_j))\\
&\cong \bigoplus_{j=1}^r\CH_{l(w)-1}(X^{F^+}(w_j)).
\end{align*}
As each $X^{F^+}(w_j)$ consists of finitely many irreducible components, $\CH_{l(w)-1}(\partial \overbar{X}(w))$ is a finitely generated free abelian group (see \cite[Example 1.3.2]{fulton2012intersection}).
The case in \cref{app_rem_2} is analogous. 
\end{rem}

\begin{proposition}\label{cycle_map_lemma}
Let $Y$ be a smooth scheme of finite type over an algebraically closed field $k$. 
Suppose that for $N\coloneq H^{1,1}_\mot(Y,\Z)$, the quotient $N/\text{tors}$ is a flat and finitely generated $\Z$-module.
If $H^1(Y)$ vanishes then so does $\CH^1(Y,1)_{\Ql}$. 
\end{proposition}
\begin{proof}
As field extensions are faithfully flat, it suffices to show that $H^1_{\et}(Y,\Q_{\ell})=0$ implies $\CH^1(Y,1)_{\Q_{\ell}}=0$.
Consider the short exact sequence 
\[
    \begin{tikzcd}
    0 \arrow[r] & T(N) \arrow[r] & N \arrow[r] & N/\text{tors} \arrow[r] & 0
    \end{tikzcd}
\] 
of $\Z$-modules.
Tensoring this with $\Z/\ell^m$ over $\Z$ gives us a short exact sequence \[\begin{tikzcd}
0 \arrow[r] & T(N)/\ell^m \arrow[r] & N/\ell^m \arrow[r] & (N/\text{tors})/\ell^m \arrow[r] & 0
\end{tikzcd}\]as $\text{Tor}^{\Z}_1(N/\text{tors},\Z/\ell^m)=0$ because $N/\text{tors}$ is a flat $\Z$-module by assumption (see \cite[Appendix B, The Tor functors]{Ma}). 
The inverse system $(T(N)/\ell^m)_m$ is Mittag-Leffler as the transition maps are surjective, hence taking the inverse limit and then tensoring with $\Q_{\ell}$ gives us a short exact sequence 
\begin{equation}\label{eq:tors_ex_seq}
    \begin{tikzcd}
        0 \arrow[r] & {\displaystyle\lim_{\leftarrow m}T(N)/\ell^m\otimes_{\Z_\ell}\Q_{\ell}} \arrow[r] & {\displaystyle\lim_{\leftarrow m}N/\ell^m\otimes_{\Z_\ell}\Q_{\ell}} \arrow[r] & {\displaystyle\lim_{\leftarrow m}(N/\text{tors})/\ell^m\otimes_{\Z_\ell}\Q_{\ell}} \arrow[r] & 0
    \end{tikzcd}
\end{equation}
as $\Q_{\ell}=\text{Frac}(\Z_\ell)$ is a flat $\Z_\ell$-module (\cite[Lemma 0598]{stacksproject}). 
We now show that the last term is isomorphic to $\CH^1(Y,1)_{\Q_{\ell}}$.
By assumption, $N/\text{tors}$ is finitely generated and so there is an isomorphism \begin{align*}
N/\text{tors}\otimes_\Z\Z_\ell\cong \lim_{\leftarrow m} (N/\text{tors})/\ell^m
\end{align*}by \cite[Theorem 8.7]{Ma}.
Hence, \begin{align*}
\CH^1(Y,1)_{\Q_\ell}&\cong H^{1,1}_\mot(Y,\Z)_{\Q_\ell}\\
&\cong N\otimes_{\Z}\Z_\ell\otimes_{\Z_\ell}\Q_\ell\\
&\cong N/\text{tors}\otimes_{\Z}\Z_\ell\otimes_{\Z_\ell}\Q_\ell\\
&\cong \lim_{\leftarrow m}((N/\text{tors})/\ell^m)\otimes_{\Z_\ell}\Q_\ell.
\end{align*}
This means that $\CH^1(Y,1)_{\Q_\ell}$ vanishes as soon as $\lim_{\leftarrow m}N/\ell^m\otimes_{\Z_\ell}\Q_{\ell}$ does.
Therefore, we are done if we show that there is an injection \begin{align*}
\displaystyle\lim_{\leftarrow m}N/\ell^m\otimes_{\Z_\ell}\Q_{\ell}\hookrightarrow H^1_{\et}(Y,\Q_{\ell}).
\end{align*}
By \cref{beilinson_lichtenbaum}, the cycle map \begin{align*}
\text{cl}^1\colon H^{1,1}_\mot(Y,\Z/\ell^m)\to H^1_\et(Y,\mu_{\ell^m})
\end{align*}
is a natural isomorphism for all positive integers $m$ (\cite[Lecture 4]{Mazza2006LectureNO}).
As $k$ contains all $\ell^m$-th roots of unity, we get natural isomorphisms of sheaves $\mu_{\ell^m}\cong \Z/\ell^m$ (see \cite[Section 16]{milneLEC}) which in particular induce natural isomorphisms $H^1_{\et}(Y,\mu_{\ell^m})\cong H^1_{\et}(Y,\Z/\ell^m)$.
In particular, this means that the inverse systems $(H^{1,1}_\mot(Y,\Z/\ell^m))_m$ and $(H^1_\et(Y,\Z/\ell^m))_m$ are isomorphic.
Hence,\begin{align*}
\lim_{\leftarrow m} H^{1,1}_\mot(Y,\Z/\ell^m)\cong \lim_{\leftarrow m} H^1_\et(Y,\Z/\ell^m)=H^1_\et(Y,\Z_\ell)
\end{align*}
as $\Z_\ell$-modules, where the equality holds by definition (see \cite[I,19]{milneLEC}).
Consider the short exact sequence of motivic complexes
\[
    \begin{tikzcd}
    0 \arrow[r] & \Z(1) \arrow[r, "\ell^m"] & \Z(1) \arrow[r] & \Z(1)/\ell^m \arrow[r] & 0.
    \end{tikzcd}
\] 
This induces a long exact sequence of motivic cohomology 
\[
    \begin{tikzcd}
    \cdots \arrow[r] & {H^{1,1}_\mot(Y,\Z)} \arrow[r, "\ell^m"]                                   & {H^{1,1}_\mot(Y,\Z)} \arrow[ld, "g"] \\
                    & {H^{1,1}_\mot(Y,\Z/\ell^m)\cong H^1_\et(Y,\Z/\ell^m)} \arrow[r] & \Pic(Y)=H^{2,1}_\mot(Y,\Z) \arrow[ld, "\ell^m"]      \\
                    & \Pic(Y) \arrow[r]                                                       & \cdots .                         
    \end{tikzcd}
\]
The morphism $g$ induces a well-defined injection \begin{align*}
N/\ell^m&\hookrightarrow H^{1,1}_\mot(Y,\Z/\ell^m)\\
[x]&\mapsto g(x).
\end{align*}
Taking the inverse limit (which is left exact) and then tensoring with $\Q_\ell$ provides us with the desired injection
\begin{equation}\label{eq:injection}
    \lim_{\leftarrow m}(N/\ell^m)\otimes_{\Z_\ell}\Q_\ell\hookrightarrow H^1_{\text{\' et}}(Y,\Q_\ell).\qedhere    
\end{equation}
\end{proof}

\begin{cor}\label{cycle_map_lemma_app_general}
    Let $U$ and $Y$ be as in \cref{fin_gen}. 
    If $H^1(U)=0$ then $\CH^1(U,1)_{\Ql}=0$.
\end{cor}
\begin{proof}
    This is the combination of \cref{fin_gen} and \cref{cycle_map_lemma}.
\end{proof}

We now want to explain how we can modify the statement of \cref{cycle_map_lemma_app_general} in a situation where $U$ and $Y$ are Deligne--Lusztig varieties (more specifically, as in \cref{app_rem}).
To do this, we need to recall some preliminaries on $G^F$-representations over $\Ql$.
For $P\subset G$ an $F$-stable parabolic subgroup and $L\subset P$ an $F$-stable Levi, we denote by 
\begin{align*}
    \Rr^{G}_{L}&\colon \Mod(\Ql[L^F])\to \Mod(\Ql[G^F])\\
    \prescript{*}{}\Rr^{G}_{L}&\colon \Mod(\Ql[G^F])\to \Mod(\Ql[L^F])
\end{align*}the Harish-Chandra induction and restriction functors respectively.
For a reference, see for example \cite[Chapter 5]{Digne_Michel_2020}.
An irreducible $G^F$-representation $\theta$ is called cuspidal if for any $F$-stable Levi $L$ of an $F$-stable parabolic subgroup of $G$ we have that $\prescript{*}{}\Rr^{G}_{L}\theta$ vanishes. 
It is called \textit{unipotent} if it appears as a constituent of some $H^i_c(X(w))$ for $i\in\N$ and $w\in W$.

\begin{definition}\label{def:PS}
    We call the set of irreducible constituents (up to isomorphism) of $\Rr^{G}_{T}\1$ the \textit{unipotent principal series} of $G^F$ (\cite[Example 3.2.6]{Geck2020TheCT}). 
    The union of the sets of irreducible constituents of $\Rr^{G}_{T}\theta$ for $\theta$ an irreducible $T^F$-representation is called the \textit{principal series} of $G^F$.
    If $V$ is a $G^F$-representation, then we denote by $V_\upr$ and $V_\pr$ the direct sum of irreducible subrepresentations with multiplicity of $V$ that are in the unipotent principal series and principal series respectively.
\end{definition}

Note that any irreducible representation in the unipotent principal series is unipotent.
Indeed, as $X(e)\cong G^F/B^F$ we get $H^0_c(X(e))\cong \Ql[G^F/B^F]\cong \Ql[G^F/U^F]\otimes_{\Ql[T^F]}\1=\Rr^{G}_{T}\1$.
Moreover, a unipotent $G^F$-representation is in the principal series if and only if it is in the unipotent principal series by \cite[Corollary 6.3]{originalpaper}.
We will use these facts implicitly throughout the preprint.

\begin{lemma}\label{B^FvsPS_and_BF_fix_exact}
    \begin{enumerate}
        \item \label{B^FvsPS_and_BF_fix_exact_1}Let $\rho$ be an irreducible $G^F$-representation.
        It is in the unipotent principal series if and only if $\rho^{B^F}\neq 0$. 
        Similarly, it is in the principal series if and only if $\rho^{U^F}\neq 0$.
        \item \label{B^FvsPS_and_BF_fix_exact_2}Taking $B^F$-fixed points is an exact functor $\Mod(\Ql[G^F])\to \Mod(\Ql)$.
        Likewise, taking $U^F$-fixed points is an exact functor $\Mod(\Ql[G^F])\to \Mod(\Ql[T^F])$.
    \end{enumerate}
\end{lemma}
\begin{proof}
    \Cref{B^FvsPS_and_BF_fix_exact_1} follows from Frobenius reciprocity. 
    For \cref{B^FvsPS_and_BF_fix_exact_2} we use that the group algebra $\Ql[B^F]$ is semisimple by Maschke's theorem, and hence any $\Ql[B^F]$-module is injective (and projective).
    The described functor can be identified with $\Hom_{\Ql[B^F]}(\1,-)$, which is exact as $\1$ is injective over $\Ql[B^F]$.
	As $T^F$ normalises $U^F$, taking $U^F$-fixed points is a functor $\Mod(\Ql[G^F])\to \Mod(\Ql[T^F])$.
	Its exactness is analogous to the $B^F$-fixed points case. 
\end{proof}

In the Deligne--Lusztig case, we obtain the following modification of \cref{cycle_map_lemma_app_general} where it is sufficient for the unipotent principal series part of $H^1(U)$ to vanish.

\begin{cor}\label{cycle_map_lemma_app}
Let $U$ and $Y$ be as in either case of \cref{app_rem}. 
If $H^1(U)^{B^F}=0$ then $\CH^1(U,1)_{\Ql}=0$. 
\end{cor}
\begin{proof}
We show first that $\CH^1(U,1)_{\Ql}$ is contained in the unipotent principal series. 
Following the proof of \cref{fin_gen}, we see that there is an exact sequence 
\[
    \begin{tikzcd}
    \cdots \arrow[r] & \F^\times \arrow[r] & {\mO_U(U)^\times\cong\CH^1(U,1)} \arrow[r] & \CH^0(Y\setminus U) \arrow[r] & \cdots
    \end{tikzcd}
\]
which gives us an injection $\CH^1(U,1)_{\Ql}\hookrightarrow \CH^0(Y\setminus U)_{\Ql}$ where the codomain is contained in the unipotent principal series (\cref{rqgen_iso_normal4} of \cref{gen_iso_normal}).
By \cref{B^FvsPS_and_BF_fix_exact_1} of \cref{B^FvsPS_and_BF_fix_exact}, the domain vanishes if and only if $\CH^1(U,1)_{\Ql}^{B^F}=0$. 
We can now apply the functor of taking $B^F$-fixed points to the short exact sequence \cref{eq:tors_ex_seq} and the injection \cref{eq:injection} of $\Ql[G^F]$-modules constructed in the proof of \cref{cycle_map_lemma}.
As this functor is exact (\cref{B^FvsPS_and_BF_fix_exact_2} of \cref{B^FvsPS_and_BF_fix_exact}), this yields the result.
\end{proof}

We want to apply \cref{cycle_map_lemma_app} to Deligne--Lusztig varieties as in \cref{app_rem}.
In \cref{sec:the_general_case}, we determine a class of $w\in F^+$ such that the vanishing condition in \cref{cycle_map_lemma_app} is satisfied.
In the untwisted type $A$ case, all $w\in F^+$ of full support and non-Coxeter satisfy these conditions (see \cref{section_type_A}). 
We will repeatedly use the fact that, in the situation of \cref{genratchow}, universal homeomorphisms induce isomorphisms on rational (higher) Chow groups.
The version that we state here is weaker that the one in \cite[Proposition 11.2.7]{Cisinski_2019} (see also \cite[11.2.2]{Cisinski_2019}). 

\begin{lemma}[Cisinski--Déglise]\label{genratchow}
Let $X\to Y$ be a separated universal homeomorphism of finite type between smooth $\F$-schemes. 
Then $H^{i,j}_\mot(X,\Q)\cong H^{i,j}_\mot(Y,\Q)$ for all $i,j$. 
In particular, $\CH^i(X)_{\Q}\cong \CH^i(Y)_{\Q}$ for every $i\in \N$.  
\end{lemma}

\subsection{Cyclic shifts}\label{sec:cyclic_shifts}

In this section, we let $W$ denote an irreducible finite Coxeter group with generating set $S$, and $F\colon W\to W$ a group automorphism such that $F(S)=S$.
As Weyl groups are in particular Coxeter groups, we can use the definitions and results of this section in the rest of this text. 
The following definition is given in \cite[Sections 2 and 6.1]{GECK2000570}. 

\begin{definition}[Geck--Pfeiffer]\label{Fcyc}
    \begin{rqenumerate}
        \item \label{Fcyc_Fconj}We say $w,w'\in W$ are \textit{$F$-conjugate} if there exists $x\in W$ such that $w'=x^{-1}wF(x)$. 
        This defines an equivalence relation on $W$, and the corresponding equivalence classes are called \textit{$F$-conjugacy classes} of $W$. 
        \item \label{Fcyc_nontriv}For $w,w'\in W$ and $s\in S$, we write $w\xrightarrow{s}_F w'$ if $swF(s)=w'$ and $l(w)\geq l(w')$. 
        We write $w\xrightarrow{s,=}_Fw'$ and $w\xrightarrow{s,>}_Fw'$ if $w\xrightarrow{s}_Fw'$ and additionally $l(w)=l(w')$ or $l(w)>l(w')$ respectively. 

        For $w,w'\in W$, we write $w\xrightarrow{}_Fw'$ if there are sequences $s_1,s_2,\ldots,s_r\in S$ and $w=w_0,w_1,\ldots,w_r=w'$ such that $w_i\in W$ and $w_{i-1}\xrightarrow{s_i}_Fw_i$ for all $i=1,\ldots,r$. 
        By convention, we always have $w\xrightarrow{}_Fw$.
        We also define the \textit{cyclic shift class} of $w$ as \[\Cyc(w)=\{v\in W:w\xrightarrow{}_Fv\quad\text{and}\quad l(w)=l(v)\}.\] 
        \item \label{Fcyc_triv}If $F$ acts trivially on $W$, we simply write $\xrightarrow{}$ for the relations defined in \cref{Fcyc_nontriv}. 
    \end{rqenumerate}
\end{definition} 

The following theorem is fundamental in our approach.
This was first shown in \cite[Theorem 2.6]{GECK2000570} and then in \cite[Theorem 3.1]{He_2012} with an alternative proof. 

\begin{theorem}[Geck--Pfeiffer, Kim, He--Nie]\label{GKP_thm}
Let $C$ be an $F$-conjugacy class of $W$.
For each $w\in C$, there exists an element $w_{\minn}\in C_{\minn}$ such that $w\xrightarrow{}_Fw_{\minn}$, where \[C_{\minn}=\{w\in C:l(w)\leq l(v)\text{ for all }v\in C\}.\] 
\end{theorem}

Naturally, elements of $C_{\minn}$ are called \textit{elements of minimal length} of $C$ (where $C$ and $C_\minn$ are as in \cref{GKP_thm}).
One reason why cyclic shift classes are important is due to the upcoming \cref{cyclic_shift_iso}.
To prove it, we first need \cref{equiv_cyc_conj}.
This lemma is related to a different way of thinking about cyclic shift classes due to Brou\'e and Michel (see \cite[Definitions 3.16, 6.7]{Broué1997} and \cite[Definition 3.1.3, Exercise 3.7]{Geck2000CharactersOF}).

\begin{lemma}\label{equiv_cyc_conj}
    Let $w,v\in W$ such that $l(w)=l(v)$.
    The following are equivalent. 
    \begin{enumerate}
        \item \label{equiv_cyc_conj_1} There exist $x,y\in W$ such that $l(x)\leq 1$, $\{xy,yF(x)\}=\{w,v\}$ and $l(w)=l(x)+l(y)$,
        \item \label{equiv_cyc_conj_2} We have $w=v$ or $w=xvF(x)$ for some $x\in S$. 
    \end{enumerate} 
\end{lemma}
\begin{proof}
    This follows from \cite[Lemma 1.2.6]{Geck2000CharactersOF}. 
\end{proof}

\begin{proposition}\label{cyclic_shift_iso}
Suppose $G^F$ is a finite group of Lie type with Weyl group $W$, as in \cref{sec:fin_gps_weyl_dynkin}.
For $w\in W$ and $v\in \Cyc(w)$, there are $G^F$-equivariant isomorphisms \begin{align*}
    \CH^i(X(w))_{\Q}&\cong \CH^i(X(v))_{\Q},\\
    H^i_c(X(w))&\cong H^i_c(X(v))
\end{align*} for all $i\in \N$.
\end{proposition}
\begin{proof}
By \cref{Fcyc_nontriv} of \cref{Fcyc}, there are sequences $s_1,s_2,\ldots,s_r\in S$ and $w=w_0,w_1,\ldots,w_r=v$ such that $w_i\in W$ and $w_{i-1}\xrightarrow{s_i,=}_Fw_i$ for all $i=1,\ldots,r$.
Let $1\leq i\leq r$.
We first want to show that there is a universal homeomorphism between $X(w_{i-1})$ and $X(w_i)$ (the direction of the arrow will not matter). 
By \cref{equiv_cyc_conj}, there exist in particular $x,y\in W$ such that $\{xy,yF(x)\}=\{w_{i-1},w_i\}$ and $l(w_{i-1})=l(x)+l(y)$. 
By \cref{univ_homeo_DL_general}, we have a universal homeomorphism $f\colon X(xy)\to X(yF(x))$ which is separated and of finite type as the source and target are.
By \cref{genratchow} and \cite[Proposition 3.1.6]{digne2006cohomologie}, we get the claim for $w_{i-1}$ and $w_i$. 
As this holds for every $i=1,\ldots,r$, we are now done. 
\end{proof}

\begin{proposition}[Steinberg]\cite[Lemma 1.8]{Lusztig1976}\label{cyc_cox}
Any two $F$-Coxeter elements have the same cyclic shift class. 
\end{proposition}

This gives us the following corollary, which verifies \cref{conj} for codimension $1$ and classical types. 

\begin{cor}[Hansen]\label{Cox_pic}
Suppose $G^F$ is a finite group of Lie type with Weyl group $W$, as in \cref{sec:fin_gps_weyl_dynkin}, and let $w\in W$ be $F$-Coxeter.
\begin{enumerate}
    \item If $G^F$ is of type $A$ then $\CH_i(X(w))_{\Q}=0$ for all $i<l(w)$. 
    \item If $G^F$ is of type $\prescript{2}{}A$, $C$, $D$ or $\prescript{2}{}D$ then $\Pic(X(w))_{\Q}=0$. 
\end{enumerate}
\end{cor}
\begin{proof}
    By \cref{cyclic_shift_iso,cyc_cox}, it is enough to show the vanishing for $w$ a particular $F$-Coxeter element $w$.
    This is implied by \cite[Corollary 4, Theorem 5]{hansen2002picard}. 
\end{proof}

We will see in \cref{sec:the_general_case} that the notion defined in \cref{def_terminal} is important to determine $\Pic(X(w))_{\Ql}$ for $w\in W$ of full support.

\begin{definition}\label{def_terminal}
    We say $w\in W$ is \textit{terminal} if it is of full $F$-support and of minimal length in its $F$-conjugacy class.
\end{definition}

$F$-Coxeter elements are terminal, but the converse is not always true.
For example, in type $B_3$ the (distinct) elements $w_0$ and $s_3s_2s_3s_2s_1$ are terminal but not Coxeter. 
However, it is true in untwisted type $A$ (\cite[Example 3.1.16]{Geck2000CharactersOF}).

\subsection{Cycle map}\label{sec:cycle_class_map}

Let $Y$ be a scheme and $0\leq i\leq \dim Y$.
The cycle map $\cl\colon \CH^i(Y)_{\Ql}\to H^{2i}(Y)$ is a comparison morphism from Chow groups to $\ell$-adic cohomology of a scheme.
For schemes with a cellular decomposition (e.g. Schubert varieties), the cycle map is an isomorphism (see \cite[Examples 1.9.1, 19.1.11.]{fulton2012intersection}). 
Discarding very trivial cases (e.g. the smooth compactification of a Coxeter variety in type $A_1$ is $\Pro^1$), the smooth compactifications $\overbar{X}(w)$ for $w\in F^+$ cannot be expected to have a cellular decomposition. 
However, if we let $G^F$ act diagonally on a product $\overbar{X}(w)\times\overbar{X}(w')$ of two such varieties then Lusztig was able to show in \cite[Sections 2 and 3]{Lusztig1998} that the amalgamated product 
\begin{align*}
    \overbar{X}(w)\times_{G^F}\overbar{X}(w')
\end{align*}
practically admits a cellular decomposition.
More precisely, the product $\overbar{X}(w)\times \overbar{X}(w')$ is stratified by locally closed $G^F$-stable subvarieties $Z_a$ such that each $G^F\backslash Z_a$ has the same cohomology as affine space of dimension $\dim Z_a$. 
Digne, Michel and Rouquier were then able to push this result further in see \cite[Section 3.4]{digne2006cohomologie}.
We state their result in a less general manner, as it is sufficient for our needs.

\begin{proposition}[Digne--Michel--Rouquier]\label{DMR_cycle_map_ver2}
    Let $v,v'\in F^+$ and consider the diagonal action of $G^F$ on $\overbar{X}(v)\times \overbar{X}(v')$. 
    Then for every $i\in \N$, there is an isomorphism\[ \CH^i(\overbar{X}(v)\times \overbar{X}(v'))^{G^F}\otimes_{\Z}\Ql(-i)\to H^{2i}(\overbar{X}(v)\times \overbar{X}(v'))^{G^F}.\]
\end{proposition}
\begin{proof}
     This is \cite[Proposition 3.4.2, Remarque 3.4.3]{digne2006cohomologie}.
\end{proof}

\begin{rem}
    The result \cite[Proposition 3.4.2]{digne2006cohomologie} also holds for products $X(v)\times X(v')$ of generalized Deligne--Lusztig varieties, not only for (smooth) projective ones as stated in \cref{DMR_cycle_map_ver2}.
    Then one gets a direct sum of motivic cohomology groups in the domain instead of a single term. 
    A result like \cref{cycle_map_lemma_app} for generalized Deligne--Lusztig varieties can also be derived from that version.
\end{rem}

To derive from \cref{DMR_cycle_map_ver2} a result for single Deligne--Lusztig varieties (as opposed to a product), we need the following.

\begin{lemma}\label{basic_repthy}
    Let $H$ be a finite group and $k$ a field of characteristic not dividing $|H|$.
	Let $V$ be an $H$-representation over $k$ and $B$ a subgroup of $H$. 
    Denote by $V^H$ and $V_H$ the $H$-invariants and coinvariants respectively.
	\begin{enumerate}
		\item \label{basic_repthy_1}There is an isomorphism $V^H\cong V_H$ of $k$-vector spaces.
		\item \label{basic_repthy_2}There is an isomorphism \begin{align*}
			(\Ind^H_B(\1_B)\otimes_k V)_H\cong V_B 
		\end{align*}of $k$-vector spaces.
	\end{enumerate}
\end{lemma}
\begin{proof}
    \Cref{basic_repthy_1} follows after checking that the morphism 
	\begin{align*}
		f&\colon V\to V^H\\
		v&\mapsto \frac{1}{|H|}\sum_{h\in H}hv
	\end{align*}
    in $\Mod(k)$ is surjective and its kernel is generated by $\{v-gv:v\in V, g\in H\}$.
	For \cref{basic_repthy_2}, one can use the isomorphism $\Ind^H_B(\1_B)\cong k[H]\otimes_{k[B]}\1_B$ and the fact that taking coinvariants is the same as tensoring with the trivial representation over the group algebra.
\end{proof}

\begin{cor}\label{DMR_cycle_iso_on_BF_ver2}
    Let $v\in F^+$, then for every $i\in \N$ there is an isomorphism 
    \begin{align*}
        \CH^i(\overbar{X}(v))^{B^F}_{\Ql}\to H^{2i}(\overbar{X}(v))^{B^F}(i).
    \end{align*}
    In particular, $\cl\colon \CH^i(\overbar{X}(v))_{\Ql,\upr}\to H^{2i}(\overbar{X}(v))_{\pr}(i)$ is an isomorphism of $G^F$-representations.
\end{cor}
\begin{proof}
    This follows from \cref{DMR_cycle_map_ver2} when we take $v=e$.
    We have $\overbar{X}(e)=G^F/B^F$ and note that $\overbar{X}(e)\times \overbar{X}(v)$ may be viewed as $|G^F/B^F|$ disjoint copies of $\overbar{X}(v)$ on which $G^F$ acts transitively.
    Hence there are natural isomorphisms of $\Ql$-vector spaces
    \begin{align*}
        \CH^i(G^F/B^F\times\overbar{X}(v))^{G^F}_{\Ql}&\cong (\Ql[G^F/B^F]\otimes_{\Ql}\CH^i(\overbar{X}(v))_{\Ql})^{G^F},\\
        H^{2i}(G^F/B^F\times\overbar{X}(v))^{G^F}&\cong (\Ql[G^F/B^F]\otimes_{\Ql}H^{2i}(\overbar{X}(v)))^{G^F}.
    \end{align*}
    We are done after applying \cref{basic_repthy}, as $\Ql[G^F/B^F]\cong \Ql[G^F]\otimes_{\Ql[B^F]}1_{B^F}$. 
    It follows from \cref{B^FvsPS_and_BF_fix_exact_1} of \cref{B^FvsPS_and_BF_fix_exact} that $\cl\colon \CH^i(\overbar{X}(v))_{\Ql}\to H^{2i}(\overbar{X}(v))(i)$ is an isomorphism on the unipotent principal series parts.
    Finally, $H^{2i}(\overbar{X}(v))_{\pr}(i)=H^{2i}(\overbar{X}(v))_{\upr}(i)$ by \cite[Proposition 4.2]{Lusztig1998}.
\end{proof}

\begin{proposition}\label{rq_DMR_cycle_iso_on_BF_3}
    Let $w\in F^+$, then $\cl\colon \Pic(X(w))_{\Ql,\upr}\hookrightarrow H^2(X(w))_\pr(1)$ is an injection of $G^F$-representations.
\end{proposition}
\begin{proof}
    Note that the statement is vacuously true if $w=e$, so we will assume $l(w)\geq 1$.
    We first claim that the open embedding $X(w)\subset \overbar{X}(w)$ induces a commutative diagram 
    \[
        \begin{tikzcd}
        \CH^0(\partial \overbar{X}(w))_{\Ql} \arrow[r] \arrow[d,"\cong"] & \Pic(\overbar{X}(w))_{\Ql} \arrow[r] \arrow[d] & \Pic(X(w))_{\Ql} \arrow[d] \arrow[r] & 0 \\
        H^2_{\partial \overbar{X}(w)}(\overbar{X}(w))(1) \arrow[r]                   & H^2(\overbar{X}(w))(1) \arrow[r]                  & H^2(X(w))(1)                            &  
        \end{tikzcd}
    \]
    in $\Mod(\Ql[G^F])$ where the rows are exact, the vertical maps are cycle maps and the left vertical map is an isomorphism. 
    The rows are exact by \cite[Proposition 1.8]{fulton2012intersection} and \cite[Theorem 9.4]{milneLEC}.
    The left square commutes by definition of the cycle map (see \cite[Section 23]{milneLEC}) and the right square commutes by \cite[Expos\'e VI, 3.3(ii)]{Illusie1977SminaireDG}.
    What is left to do for the claim is to show the left vertical map is an isomorphism. 
    Set \begin{align*}
    Z\coloneq \partial\overbar{X}(w)=\bigcup_{{\substack{v\preccurlyeq w\\l(v)=l(w)-1}}}\overbar{X}(v)
    \end{align*}(we write $v\preccurlyeq w$ if $v$ is a subword of $w$) and consider its open and smooth subscheme 
    \begin{align*}
        U\coloneq \bigcup_{{\substack{v\preccurlyeq w\\l(v)=l(w)-1}}}X(v).
    \end{align*}
    We note $\CH^0(Z)\cong \CH^0(U)$ as $U$ and $Z$ have the same number of irreducible components (one can also use \cite[Proposition 1.8]{fulton2012intersection} to see this).
    As $U$ is smooth, the cycle map $\CH^0(U)_{\Ql}\to H^0(U)$ is an isomorphism. 
    We now show $H^0(U)\cong H^2_{Z}(\overbar{X}(w))(1)$. 
    If $l(w)=1$ then $Z=U=X(e)$ is smooth, and the desired isomorphism follows from purity (see \cite[Theorem 16.1, Theorem 16.7]{milneLEC}).
    If $l(w)>1$, we proceed as follows.
    Consider the smooth varieties \begin{align*}
        U'&\coloneq \bigcup_{{\substack{v\preccurlyeq w\\l(v)\geq l(w)-1}}}X(v),\\
        V'&\coloneq X(w).
    \end{align*}
    As $U$ is smooth, closed in $U'$ and of codimension $1$, we again get by purity $H^0(U)\cong H^2_U(U')(1)$.
    As we have open embeddings $V'\subset U'\subset \overbar{X}(w)$, we have a long exact sequence (\cite[Page 138]{milneLEC})
    \[
        \begin{tikzcd}
        \cdots \arrow[r] & H^2_{\overbar{X}(w)\setminus U'}(\overbar{X}(w))(1) \arrow[r] & H^2_{\overbar{X}(w)\setminus V'}(\overbar{X}(w))(1) \arrow[r] & H^2_{U'\setminus V'}(U')(1) \arrow[ld] \\
                        &                                 & H^3_{\overbar{X}(w)\setminus U'}(\overbar{X}(w))(1) \arrow[r] & \cdots.                          
        \end{tikzcd}
    \]
    Now, $\dim \overbar{X}(w)\setminus U'=l(w)-2$ and this means $H^i_{\overbar{X}(w)\setminus U'}(\overbar{X}(w))(1)=0$ for $2\leq i\leq 3$ by semi-purity (\cite[Lemma 23.1]{milneLEC}).
    Putting everything together gives us 
    \begin{align*}
        H^0(U)&\cong H^2_U(U')(1)\\
        &=H^2_{U'\setminus V'}(U')(1)\\
        &\cong H^2_{\overbar{X}(w)\setminus V'}(\overbar{X}(w))(1)\\
        &=H^2_{Z}(\overbar{X}(w))(1)
    \end{align*}
    which finishes the proof of the claim.
    Taking $B^F$-fixed points is exact (\cref{B^FvsPS_and_BF_fix_exact_2} of \cref{B^FvsPS_and_BF_fix_exact}), which means that we also have a commutative diagram where the rows are exact
    \[
        \begin{tikzcd}
        \CH^0(\partial \overbar{X}(w))_{\Ql} \arrow[r] \arrow[d,"\cong"] & {\Pic(\overbar{X}(w))_{\Ql,\upr}} \arrow[r] \arrow[d] & {\Pic(X(w))_{\Ql,\upr}} \arrow[r] \arrow[d] & 0 \\
        H^2_{\partial \overbar{X}(w)}(\overbar{X}(w))(1) \arrow[r]                   & H^2(\overbar{X}(w))_{\upr}(1) \arrow[r]                   & H^2(X(w))_{\upr}(1)                            &  
        \end{tikzcd}
    \]
    in $\Mod(\Ql[G^F])$ after restricting to the unipotent principal series (see \cref{B^FvsPS_and_BF_fix_exact_1} of \cref{B^FvsPS_and_BF_fix_exact} and note that $\CH^0(\partial \overbar{X}(w))_{\Ql}$ and $H^2_{\partial \overbar{X}(w)}(\overbar{X}(w))(1)$ are in the unipotent principal series). 
    The left vertical arrow is an isomorphism, and so is the middle vertical arrow by \cref{DMR_cycle_iso_on_BF_ver2}.
    A diagram chase gives us that $\cl\colon \Pic(X(w))_{\Ql,\upr}\to H^2(X(w))_\upr(1)$ is an injection.
    We are done as $H^2(X(w))_\upr=H^2(X(w))_\pr$ as $H^2(X(w))$ is a direct sum of unipotent $G^F$-representation.
\end{proof}

As we are only interested in $G^F$-representations, we will often disregard the Tate twist in \cref{DMR_cycle_iso_on_BF_ver2,rq_DMR_cycle_iso_on_BF_3}.
\Cref{rq_DMR_cycle_iso_on_BF_3} motivates us to compute rational Picard groups of Deligne--Lusztig varieties, as we will be computing at least part of the principal series part of the $\ell$-adic cohomology in degree $2$ of the Deligne--Lusztig variety. 
If $w\in W$ is $F$-Coxeter, we will see in \cref{sec:the_cox_case} that \cref{rq_DMR_cycle_iso_on_BF_3} does not provide us with interesting information as the domain vanishes.
However, outside of the $F$-Coxeter case we will see (\cref{sec:the_general_case,section_type_A,small_rk_comps}) that $\cl\colon \Pic(X(w))_{\Ql,\upr}\hookrightarrow H^2(X(w))_\pr$ is actually very often an isomorphism if $w\in W$ is of full $F$-support (especially in type $A$). 

\begin{rem}
    It is important to point out the following.
    For $w\in F^+$, the cohomology groups $H^i_c(\overbar{X}(w))$ are known for all $i$ in the sense of \cite[Corollaire 3.3.8 (ii)]{digne2006cohomologie}, and can be explicitly computed using the command $\mathsf{DeligneLusztigLefschetz}$ in GAP. 
    We will not use this in this preprint.
\end{rem}

\begin{cor}\label{Pic_PS}
    Let $w\in F^+$.
    If $\Pic(X(w))_{\Ql}=0$ then $\Pic(\overbar{X}(w))_{\Ql}=\Pic(\overbar{X}(w))_{\Ql,\upr}=\Pic(\overbar{X}(w))_{\Ql,\pr}$ and we have an isomorphism of $G^F$-representations \begin{align*}
        \Pic(\overbar{X}(w))_{\Ql}\cong \CH^0(\partial \overbar{X}(w))_{\Ql}/\CH^1(X(w),1)_{\Ql}.
    \end{align*}  
\end{cor}
\begin{proof}
    We have a localization sequence 
    \[
        \begin{tikzcd}
        0 \arrow[r] & {\CH^1(X(w),1)_{\Ql}} \arrow[r] & \CH^0(\partial \overbar{X}(w))_{\Ql} \arrow[r] & \Pic(\overbar{X}(w))_{\Ql} \arrow[r] & 0.
        \end{tikzcd}
    \]
    As $\Ql[G^F]$ is semisimple, $\Pic(\overbar{X}(w))_{\Ql}$ is a quotient of $\CH^0(\partial \overbar{X}(w))_{\Ql}$.
    By \cref{rqgen_iso_normal4} of \cref{gen_iso_normal}, the latter is a direct sum of representations of the form $\Ind^{G^F}_{P^F}\1$ where $P\subset G$ is a standard parabolic subgroup. 
\end{proof}

Recall that the vanishing condition $\Pic(X(w))_{\Ql}=0$ in \cref{Pic_PS} holds conjecturally (\cref{conj}) for any $F$-Coxeter element $w\in W$, and is proved if $G^F$ is of classical type (\cref{Cox_pic}).

\subsection{The Coxeter case}\label{sec:the_cox_case}

Recall the unipotent radical $U\subset B$ with $B=U\rtimes T$. 
In \cite[Corollary 2.7]{Lusztig1976}, Lusztig showed that the quotient of a Coxeter variety by $U^F$ is a torus.
Hence, \cref{conj} is true on the principal series part as stated in \cref{conj_Hansen_BF_allcodim} below (see also \cite[3.1.3]{hansen}).
Recall the function $\varphi\colon S\to S_F$ from \cref{sec:fin_gps_weyl_dynkin}.
To state our results, we need the following definition.
For an $F$-stable $J\subset S$, denote by $i_J$ the Harish-Chandra induction $\Rr^{G}_{L_J}\1$.
As we are inducing the trivial representation, this coincides with ordinary induction $\Ind^{G^F}_{P_J^F}\1$.
This is what Orlik denotes by $i^{G}_{P_J}$ in \cite{orlik2014cohomology}\footnote{Orlik denotes by $\mathbf{G}$ the reductive group $\Gl_n$ over $\F$ and by $G$ the finite group of Lie type $\Gl_n(\F_q)$ obtained via the standard Frobenius $(a_{ij})\mapsto (a_{ij}^q)$.}.
The representations $i_J$ have been well-studied for any type, and their decomposition in terms of irreducible representations comes down to decomposing the ordinary induction $\Ind^{W^F}_{W_J^F}\1$ (see \cite[Example 3.2.6, Theorem 3.2.7]{Geck2020TheCT}).
This is easily done using GAP (see \cite{Sch97}), and we will use it for some examples in \cref{small_rk_comps}.

\begin{proposition}\label{conj_Hansen_BF_allcodim}
    Let $w\in W$ be $F$-Coxeter and $k>0$.
    \begin{enumerate}
        \item \label{conj_Hansen_BF_allcodim_1} The fixed points $\CH^k(X(w))_{\Q}^{U^F}$ vanish. 
        In particular, $\CH^k(X(w))_{\Ql,\pr}=0$.
        \item \label{conj_Hansen_BF_allcodim_2} There is an isomorphism of $G^F$-representations $\CH^1(X(w),1)_{\Ql}\cong H^1(X(w))_\pr$.
        \item \label{conj_Hansen_BF_allcodim_3}There is an isomorphism of $G^F$-representations
        \begin{align*}
            \Pic(\overbar{X}(w))_{\Ql,\pr}&\cong (\bigoplus_{J\in S_F}i_{\varphi^{-1}(S_F\setminus \{J\})})/H^1(X(w))_\pr.
        \end{align*}
        In particular, $\Pic(\overbar{X}(w))_{\Ql,\pr}=\Pic(\overbar{X}(w))_{\Ql,\upr}$.
        %If $\Pic(X(w))_{\Q}=0$ then $\Pic(\overbar{X}(w))_{\Ql,\pr}\cong \Pic(\overbar{X}(w))_{\Ql}$.
    \end{enumerate}
\end{proposition}
\begin{proof}
    The result \cite[Corollary 2.7]{Lusztig1976} says $U^F\backslash X(w)\cong \G_m^{l(w)}$, and hence $\CH^k(X(w))^{U^F}_{\Q}\cong \CH^k(\G_m^{l(w)})_{\Q}$ by \cite[Example 1.7.6]{fulton2012intersection}.
    As $\G_m^{l(w)}$ is open in $\A^{l(w)}$, its Chow groups vanish (even integrally).
    We are done with \cref{conj_Hansen_BF_allcodim_1} by \cref{B^FvsPS_and_BF_fix_exact_1} of \cref{B^FvsPS_and_BF_fix_exact}.
    For \cref{conj_Hansen_BF_allcodim_2}, we revisit the proof of \cref{cycle_map_lemma} with $Y\coloneq U^F\backslash X(w)\cong \G_m^{l(w)}$.
    We first claim that $N\coloneq H^{1,1}_\mot(Y)\cong \F^{\times}\times\Z^{l(w)}$.
    Indeed, the closed complement $Z$ of $\G_m^{l(w)}$ in $\A^{l(w)}$ is a union of hyperplanes
    \begin{align*}
        Z&=\bigcup^{l(w)}_{i=1}\{x\in \A^{l(w)}:x_i=0\}\\
        &\cong \bigcup^{l(w)}_{i=1}\A^{l(w)-1}
    \end{align*}
    and hence has $l(w)$ irreducible components.
    This means that the localization sequence 
    \[
        \begin{tikzcd}
        \cdots \arrow[r] & {\CH^0(Z,1)} \arrow[r] & {\CH^1(\A^{l(w)},1)} \arrow[r] & {\CH^1(\G_m^{l(w)},1)} \arrow[r] & \CH^0(Z) \arrow[r] & \Pic(\A^{l(w)}) \arrow[r] & \cdots
        \end{tikzcd}
    \]
    simplifies to a short exact sequence 
    \[
        \begin{tikzcd}
        0 \arrow[r] & \F^\times \arrow[r] & {\CH^1(\G_m^{l(w)},1)} \arrow[r] & \Z^{l(w)} \arrow[r] & 0
        \end{tikzcd}
    \]
    of abelian groups by the fact that $\CH^0(Z,1)\cong H^{-1,0}_\mot(Z,\Z)=0$ by \cite[Corollary 4.2]{Mazza2006LectureNO}.
    As $\Z^{l(w)}$ is a free (and hence projective) $\Z$-module, the short exact sequence splits and we get our claim. 
    In particular, $T(N)=\F^\times$ and we have $\F^\times/\ell^m=0$ for all $m\in \N$ as $\F$ is algebraically closed.
    This means 
    \begin{align*}
        \displaystyle\lim_{\leftarrow m}T(N)/\ell^m=0
    \end{align*}
    and so the proof of \cref{cycle_map_lemma} gives us an injection 
    \begin{align*}
        \CH^1(Y,1)_{\Ql}\hookrightarrow H^1(Y)\cong H^1(X(w))^{U^F}.
    \end{align*}
    As $X(w)$ is affine\footnote{It has now been shown by Brosnan--Hong--Lee that all Deligne--Lusztig varieties are affine \cite[Corollary 3.12]{brosnan2026geometryregularsemisimplelusztig}. We will not need this result.} (\cite[Corollary 2.8]{Lusztig1976}), $\mO_Y(Y)\cong \mO_{X(w)}(X(w))^{U^F}$ (\cite[Proposition 12.27]{görtz2010algebraic}) and hence we may simplify the domain as $\CH^1(X(w),1)_{\Ql}^{U^F}$. 
    As $\Ql$-vector spaces, we have an isomorphism $\CH^1(X(w),1)_{\Ql}^{U^F}\cong \Ql^{l(w)}$.
    This induces an injection of $G^F$-representations $\CH^1(X(w),1)_{\Ql,\pr}\hookrightarrow H^1(X(w))_\pr$ where the domain is nonzero (\cref{B^FvsPS_and_BF_fix_exact_1} of \cref{B^FvsPS_and_BF_fix_exact}).
    As $H^1(X(w))_\pr$ is an irreducible $G^F$-representation by \cite[Section 7.3, Theorem 6.1]{Lusztig1976}, this is actually an isomorphism of $G^F$-representations.
    We are done with \cref{conj_Hansen_BF_allcodim_2} once we show $\CH^1(X(w),1)_{\Ql}=\CH^1(X(w),1)_{\Ql,\pr}$.
    This follows from the claim that $\CH^1(X(w),1)_{\Ql}$ lies in the unipotent principal series, which we now show.
    The open embedding $X(w)\subset \overbar{X}(w)$ gives us a localization sequence 
    \[
        \begin{tikzcd}
        0 \arrow[r] & {\CH^1(X(w),1)_{\Ql}} \arrow[r] & \CH^0(\partial \overbar{X}(w))_{\Ql} \arrow[r] & \Pic(\overbar{X}(w))_{\Ql} \arrow[r] & \Pic(X(w))_{\Ql} \arrow[r] & 0
        \end{tikzcd}
    \]
    which tells us in particular that $\CH^1(X(w),1)_{\Ql}$ is contained in $\CH^0(\partial \overbar{X}(w))_{\Ql}$.
    By \cref{rqgen_iso_normal4} of \cref{gen_iso_normal}, this latter representation is in the unipotent principal series (we will actually compute it explicitly in the proof of \cref{conj_Hansen_BF_allcodim_3}).
    Let us prove \cref{conj_Hansen_BF_allcodim_3}.
    By \cref{conj_Hansen_BF_allcodim_1}, $\Pic(X(w))_{\Ql}^{U^F}=0$.
    Taking $U^F$-fixed points of the localization sequence above gives us a short exact sequence
    \[
        \begin{tikzcd}
        0 \arrow[r] & {\CH^1(X(w),1)_{\Ql}^{U^F}} \arrow[r] & \CH^0(\partial \overbar{X}(w))_{\Ql}^{U^F} \arrow[r] & \Pic(\overbar{X}(w))_{\Ql}^{U^F} \arrow[r] & 0
        \end{tikzcd}
    \]
    by \cref{B^FvsPS_and_BF_fix_exact_2} of \cref{B^FvsPS_and_BF_fix_exact}.
    Due to \cref{B^FvsPS_and_BF_fix_exact_1} of \cref{B^FvsPS_and_BF_fix_exact}, we have 
    \begin{align*}
        \Pic(\overbar{X}(w))_{\Ql,\pr}&\cong \CH^0(\partial \overbar{X}(w))_{\Ql,\pr}/\CH^1(X(w),1)_{\Ql,\pr}
    \end{align*}
    (one can compare with \cref{Pic_PS}, where something similar happens).
    We have $\CH^1(X(w),1)_{\Ql}\cong H^1(X(w))_\pr\cong \CH^1(X(w),1)_{\Ql,\pr}$ by \cref{conj_Hansen_BF_allcodim_2}, and so we are done once we compute $\CH^0(\partial \overbar{X}(w))_{\Ql}=\CH^0(\partial \overbar{X}(w))_{\Ql,\pr}$.
    As $w$ is an $F$-Coxeter element we know that for any $v\leq w$ of length $l(w)-1$ the variety $\overbar{X}(v)$ has $|\varphi^{-1}(S_F\setminus \{J\})|$ irreducible components and $G^F$ acts transitively on them, where $J\in S_F$ is the unique $F$-orbit such that $\supp_F(v)=S_F\setminus \{J\}$ (see \cref{rqgen_iso_normal4} of \cref{gen_iso_normal}).
    This means \begin{align*}
        \CH^0(\partial \overbar{X}(w))_{\Ql}\cong \bigoplus_{J\in S_F}i_{\varphi^{-1}(S_F\setminus \{J\})}
    \end{align*}
    by \cref{subsec:elsofpropsupp}.
\end{proof}

\subsection{The general case}\label{sec:the_general_case}

In this section, we consider $\Pic(X(w))_{\Ql}$ for $w\in W$ of full support and not of minimal length in its $F$-conjugacy class (in particular, non-terminal elements of $W$). 
\cref{inductive_step_general} allows us to compute it recursively in some suitable cases. 
For instance, if $G^F$ is of untwisted type $A$, it allows us to determine all rational Picard groups for $w\in W$ of full support (see \cref{section_type_A}). 
In practice, we use the vanishing of a particular $\ell$-adic cohomology group.
This reduces to the vanishing of a cohomology group of a variety of codimension $1$, as specified in the upcoming \cref{H1_vanishing_lower_impl_higher}.
To prove it, we first prove the following helpful \cref{helpfullemma} which we will use multiple times.

\begin{lemma}\label{helpfullemma}
    Let $s\in S$ and $v\in F^+$ such that $w\coloneq svF(s)$ has full $F$-support and set $t\coloneq F(s)$. 
    Then $H^{2l(w)-1}_c(X(w))\cong H^{2l(vt)-1}_c(X(v\underline t))(-1)$.
\end{lemma}
\begin{proof}
    There is a universal homeomorphism $X(w)\to X(vtt)$ (\cref{univ_homeo_DL_general}), it is equivalent to prove the statement for $X(vtt)$ (\cref{cyclic_shift_iso}).  
    The open embedding of $X(vtt)$ in $X(vt\underline t)$ (see \cref{partial_compactification}) induces a long exact sequence 
    \[
        \begin{tikzcd}
        \cdots \arrow[r] & H^{2l(w)-2}_c(X(vt)) \arrow[r,"f"] & H^{2l(w)-1}_c(X(vtt)) \arrow[r,"g"] & H^{2l(w)-1}_c(X(vt\underline t)) \arrow[r] & 0.
        \end{tikzcd}
    \]
    As $vt$ has full $F$-support and $2l(w)-2=2l(vt)$, the domain of the morphism $f$ is $\1(-l(vt))$ (\cite[Proposition 3.3.14.]{digne2006cohomologie}). 
    In fact, the same result tells us that $\1$ cannot appear in $H^{2l(w)-1}_c(X(vtt))$, which means $f=0$ by Schur's lemma.
    Hence $g$ is an isomorphism.
    Moreover, we can apply \cref{affine_bundle} to obtain \begin{align*}
        H^{2l(w)-1}_c(X(vtt))&\cong H^{2l(w)-1}_c(X(vt\underline t))\\
        &\cong H^{2l(w)-3}_c(X(v\underline t))(-1)\\
        &=H^{2l(vt)-1}_c(X(v\underline t))(-1).
    \end{align*}
\end{proof}

\begin{proposition}\label{H1_vanishing_lower_impl_higher}
    Let $s,v,w$ and $t$ be as in \cref{helpfullemma}.
    If $H^1(X(vt))=0$ then $H^1(X(w))=0$.
    Analogously, $H^1(X(vt))^{B^F}=0$ implies $H^1(X(w))^{B^F}=0$.
\end{proposition}
\begin{proof}
    By \cref{helpfullemma} and Poincar\'e duality, it is equivalent to show $H^{2l(vt)-1}_c(X(v\underline t))=0$.
    The open embedding $X(vt)\subset X(v\underline t)$ gives a long exact sequence \[\begin{tikzcd}
    \cdots \arrow[r] & H^{2l(vt)-1}_c(X(vt)) \arrow[r] & H^{2l(vt)-1}_c(X(v\underline t)) \arrow[r] & H^{2l(vt)-1}_c(X(v)) \arrow[r] &\cdots.
    \end{tikzcd}\]
    As $2l(vt)-1=2l(v)+1$, we know that $H^{2l(vt)-1}_c(X(v))=0$.
    We are done with the first implication because $H^{2l(vt)-1}_c(X(vt))$ also vanishes by Poincar\' e duality and our assumption. 
    The second also follows from the proof above, as taking $B^F$-fixed points is exact (\cref{B^FvsPS_and_BF_fix_exact_2} of \cref{B^FvsPS_and_BF_fix_exact}).
\end{proof}

We now want to show that there are Weyl group elements $w$ such that $H^1(X(w))_\pr=0$ to apply \cref{cycle_map_lemma_app} (recall \cref{B^FvsPS_and_BF_fix_exact_1} of \cref{B^FvsPS_and_BF_fix_exact}).
More precisely, we want to prove \cref{H1_vanishing_other_gps} and to do this we need \cref{cohom_quot,gen_by_I,funny_root_lemma}.
For \cref{cohom_quot}, note that a stronger result is treated in \cite[Proposition 8.1.10]{Digne_Michel_2020} for example.

\begin{proposition}\label{cohom_quot}
    Let $H$ be a finite group acting from the left on a quasi-projective variety $Y$ over $\F$. 
    The $\Ql$-vector spaces $H^i_c(H\backslash Y)$ and $H^i_c(Y)^H$ are isomorphic for all $i$.
\end{proposition}

\begin{lemma}\label{gen_by_I}
    Let $w\in W$, $I\subset S$ and recall $\Phi^+$ and $\Phi_I$ from \cref{sec:fin_gps_weyl_dynkin}.
    The following are equivalent
    \begin{rqenumerate}
    \item \label{gen_by_I_1}For all $\alpha\in \Phi^+\setminus \Phi_I$ we have $w(\alpha)\geq 0$,
    \item \label{gen_by_I_2}The element $w$ is in $W_I$.
    \end{rqenumerate}
\end{lemma}
\begin{proof}
    \Cref{gen_by_I_2} implies \cref{gen_by_I_1} is clear.
    For the other implication, one can use \cite[Last Corollary before III.10.3]{Humphreys1973IntroductionTL} and induction.
\end{proof}

\begin{lemma}\label{funny_root_lemma}
    Let $I\subset S$.
    For $g\in G$, we have $^gB\cap U_I=1$ if and only if $g\in P_Iw_0B$. 
\end{lemma}
\begin{proof}
    By the Bruhat decomposition, we know $g=bww_0b'$ for some $w\in W$ and $b,b'\in B$ (see \cite[Proposition 2.3.4(i)]{bjorner2006combinatorics}). 
    As $B\subset P_I=U_I\rtimes L_I$, we get $\prescript{b}{}{U_I}=U_I$ and so $\prescript{g}{}B\cap U_I= \prescript{b}{}{(\prescript{w}{}{B^-}\cap U_I)}$.
    This is equal to $1$ if and only if $\prescript{w}{}{B^-}\cap U_I=1$. 
    Moreover, as $U_I=\prod_{\alpha\in \Phi^+\setminus \Phi_I}U_\alpha$ (see for example \cite[Proposition 3.4.1]{Digne_Michel_2020}) we get \begin{align*}
    \prescript{w}{}{B^-}\cap U_I&=B^-\cap \prescript{w^{-1}}{}{U_I}\\
    &=\prod_{\substack{\alpha\in\Phi^+\setminus \Phi_I\\w^{-1}(\alpha)\leq 0}}U_{w^{-1}(\alpha)}
    \end{align*}
    which means that this equals $1$ if and only if $\{\alpha\in \Phi^+\setminus \Phi_I:w^{-1}(\alpha)\leq 0\}=\varnothing$.
    By \cref{gen_by_I}, this is equivalent to $w\in W_I$, which is in turn equivalent to $g\in P_Iw_0B$.
\end{proof}

We obtain the following vanishing result for $G^F$ of any type.

\begin{proposition}\label{H1_vanishing_other_gps}
    Let $s\in S$ and $v\in W_{S\setminus F\cdot s}$ such that $vF(s)$ is $F$-Coxeter.
    \begin{enumerate}
        \item \label{H1_vanishing_other_gps1}
        If $G^F$ is not of type $\prescript{2}{}{A_{2n}}, \prescript{2}{}{B_2}, \prescript{2}{}{F_4}$ or $\prescript{2}{}{G_2}$ then $H^{1}(X(svF(s)))=0$.
        \item \label{H1_vanishing_other_gps2} 
        If $G^F$ is of type $\prescript{2}{}{A_{2n}}, \prescript{2}{}{B_2}, \prescript{2}{}{F_4}$ or $\prescript{2}{}{G_2}$, then $H^{1}(X(svF(s)))_\pr=0$.
    \end{enumerate}
\end{proposition}
\begin{proof}
Set $t\coloneq F(s)$ and $w\coloneq svt$.
Note that $w\in W$ is of full $F$-support and $l(w)=l(v)+2$. 
By \cref{helpfullemma}, we can reduce to show the vanishing of $H^{2l(sv)-1}_c(X(v\underline t))$.
Like in the proof of \cref{H1_vanishing_lower_impl_higher}, the long exact sequence induced by the open embedding $X(vt)\subset X(v\underline t)$ in particular gives us a surjection \begin{align}\label{surj_1}
H^{2l(sv)-1}_c(X(vt))\twoheadrightarrow H^{2l(sv)-1}_c(X(v\underline t)).
\end{align}
To deduce the desired vanishing, we will use that $vt$ is $F$-Coxeter and the results in \cite[Proof of Corollary 2.10]{Lusztig1976}.
For $I\coloneq \varphi^{-1}(S_F\setminus \{F\cdot (\prescript{w_0}{}t)\})$, Lusztig considers the variety \[Y\coloneq \{B_0\in X(v\underline t):B_0\cap U_I=1\}\]and proves that the quotient $\widetilde Y\coloneq U_I^F\backslash Y$ is isomorphic to $X_{L_I}(c)\times \A^1$ where $c\in W_I$ is $F$-Coxeter, and that $Y$ contains $X(vt)$ as an open subscheme.  
This implies $l(c)=l(sv)-1$ and $\dim X(v\underline t)\setminus Y\leq \dim X(v\underline t)\setminus X(vt)=\dim X(v)=l(v)$.
By \cref{funny_root_lemma}, $Y=X(v\underline t)\cap P_Iw_0B/B$ which shows that $Y\subset X(v\underline t)$ is open as $P_Iw_0B/B\subset X$ is open. 
The open embedding $Y\subset X(v\underline t)$ induces a long exact sequence 
\[
    \begin{tikzcd}
    \cdots \arrow[r] & H^{2l(sv)-1}_c(Y) \arrow[r] & H^{2l(sv)-1}_c(X(v\underline t)) \arrow[r] & H^{2l(sv)-1}_c(X(v\underline t)\setminus Y) \arrow[r] & \cdots
    \end{tikzcd}
\]
where $H^{2l(sv)-1}_c(X(v\underline t)\setminus Y)$ vanishes as $2l(sv)-1=2l(v)+1>2\dim X(v\underline t)\setminus Y$. 
Taking $U_I^F$-fixed points (this is a right-exact\footnote{It is even exact, by the same proof as for \cref{B^FvsPS_and_BF_fix_exact_2} of \cref{B^FvsPS_and_BF_fix_exact}, but we do not need this.} functor $\Mod(\Ql[G^F])\to \Mod(\Ql[L_I^F])$) of the surjection $H^{2l(sv)-1}_c(Y)\twoheadrightarrow H^{2l(sv)-1}_c(X(v\underline t))$ gives us a surjection 
\begin{align*}
    H^{2l(sv)-1}_c(\widetilde Y) \twoheadrightarrow H^{2l(sv)-1}_c(X(v\underline t))^{U^F_I}
\end{align*}
by \cref{cohom_quot}.
By the Künneth formula, we can simplify the domain as \begin{align*}
H^{2l(sv)-1}_c(\widetilde Y)\cong H^{2l(sv)-3}_c(X_{L_I}(c))(-1)=H^{2l(v)-1}_c(X_{L_I}(c))(-1)
\end{align*}
giving us a surjection \begin{align}\label{surj_2}
H^{2l(v)-1}_c(X_{L_I}(c))(-1)\twoheadrightarrow H^{2l(sv)-1}_c(X(v\underline t))^{U^F_I}.
\end{align}
We now prove \cref{H1_vanishing_other_gps1}.
Suppose $H^{2l(sv)-1}_c(X(v\underline t))\neq 0$.
By our assumption on the type of $G^F$, the $G^F$-representation  $H^{2l(sv)-1}_c(X(vt))$ is irreducible (and in the principal series) by \cite[Section 7.3, Theorem 6.1, Proposition 6.3]{Lusztig1976}\footnote{Lusztig denotes type $\prescript{2}{}{A_3}$ as $\prescript{2}{}{D_3}$, see \cite[(1.5)]{Lusztig1976}.}.
By \cref{surj_1}, $H^{2l(sv)-1}_c(X(v\underline t))\cong H^{2l(sv)-1}_c(X(vt))$.
Hence we also have an isomorphism after taking $U_I^F$-fixed points, and the result \cite[Corollary 2.10]{Lusztig1976} gives us\begin{align*}
H^{2l(sv)-1}_c(X(v\underline t))^{U_I^F}&\cong H^{2l(sv)-1}_c(X(vt))^{U_I^F}\\
&\cong H^{2l(sv)-2}_c(X_{L_I}(c))\oplus H^{2l(sv)-3}_c(X_{L_I}(c))(-1)\\
&= H^{2l(v)}_c(X_{L_I}(c))\oplus H^{2l(v)-1}_c(X_{L_I}(c))(-1).
\end{align*}This contradicts the surjectivity of \cref{surj_2}, as $H^{2l(v)}_c(X_{L_I}(c))$ contains the trivial representation (see \cite[Proposition 1.22]{Lusztig1976}).
Therefore, $H^{2l(sv)-1}_c(X(v\underline t))=0$.
To prove \cref{H1_vanishing_other_gps2}, suppose $H^{2l(sv)-1}_c(X(v\underline t))_\pr\neq 0$.
In this case, $H^{2l(sv)-1}_c(X(vt))$ is a sum of an irreducible principal series representation and a finite number of irreducibles not in the principal series (this finite number is one less than the number of rows in the tables for the considered types in \cite[Section 7.3]{Lusztig1976}).
Hence, $H^{2l(sv)-1}_c(X(vt))_\pr$ is irreducible and \cref{surj_1} gives an isomorphism $H^{2l(sv)-1}_c(X(vt))_\pr\cong H^{2l(sv)-1}_c(X(v\underline t))_\pr$.
Using \cite[Proposition 5.3.9]{Digne_Michel_2020}, we see that the same proof as for \cref{H1_vanishing_other_gps1} gives us $(H^{2l(sv)-1}_c(X(v\underline t))_\pr)^{U_I^F}\cong H^{2l(v)}_c(X_{L_I}(c))_\pr\oplus H^{2l(v)-1}_c(X_{L_I}(c))_\pr(-1)$ which contradicts the surjectivity of \cref{surj_2}.
\end{proof}

Before treating our main theorem, we first observe that the vanishing of $H^1(X(w))$ holds for $w=w_0$ in almost all cases, as a consequence of the following. 

\begin{lemma}\label{cohom_codim_1_general}
    Let $w\in W$ and suppose \begin{enumerate}
        \item \label{cohom_codim_1_general1}
        The cohomology group $H^{2l(w)-1}_c(\overline{X(w)})$ vanishes, 
        \item \label{cohom_codim_1_general2}
        The representation $H^{2l(w)-2}_c(\partial \overline{X(w)})$ is a (finite) direct sum of trivial representations.
    \end{enumerate} 
    Then $H^1(X(w))=0$.
    The same implication is true if $w\in F^+$ and we replace $\overline{X(w)}$ by $\overbar{X}(w)$ in \cref{cohom_codim_1_general1} and $\partial\overline{X(w)}$ by $\partial\overbar{X}(w)$ in \cref{cohom_codim_1_general2}. 
\end{lemma}
\begin{proof}
    We first consider the case where $w\in W$.
    By Poincar\'e duality, it is equivalent to show $H^{2l(w)-1}_c(X(w))=0$.
    By the first assumption, the long exact sequence in cohomology gives us a surjection $f\colon H^{2l(w)-2}_c(\partial \overline{X(w)})\to H^{2l(w)-1}_c(X(w))$. 
    Suppose that $f\neq 0$, then one of the irreducible constituents $\rho$ of the domain gives us a nonzero subrepresentation $\im (f|_\rho)\subset H^{2l(w)-1}_c(X(w))$.
    As $\rho\cong \1$ by the second assumption and $f$ is a morphism of $G^F$-representations, $\im (f|_{\1})\cong \1$. 
    This contradicts \cite[Proposition 3.3.14.]{digne2006cohomologie}, and so $f=0$. 
    We are done as $f$ is surjective. 
    The case where $w\in F^+$ is analogous.
\end{proof}

\begin{cor}
    If $G^F$ is not of type $A_1$, then $H^1(X(w_0))=0$.
\end{cor}
\begin{proof}
    Note $\overline{X(w_0)}=X$ only has cohomology in even degrees (it has a cellular decomposition given by Schubert cells), and that $X\setminus X(w_0)$ is of codimension $1$ in $X$.
    As $G^F$ is not of type $A_1$, we have seen in the proof of \cref{Picard_longest_elem} that $X\setminus X(w_0)$ consists of $|S|$ irreducible components.  
    This means $H^{2l(w_0)-2}_c(X\setminus X(w_0))\cong \1^{|S|}$.
    Now we can apply \cref{cohom_codim_1_general}.
\end{proof}

We can now state and prove our main theorem.

\begin{theorem}\label{inductive_step_general}
    Let $s\in S$ and $v\in F^+$ such that $w\coloneq svF(s)$ has full $F$-support.
    \begin{rqenumerate}
        \item \label{inductive_step_general_1}Suppose $\CH^1(X(w),1)_{\Ql}$ and $\CH^1(X(vF(s)),1)_{\Ql}$ vanish (by \cref{H1_vanishing_lower_impl_higher,cycle_map_lemma_app}, this is the case whenever $H^1(X(vF(s)))^{B^F}$ does). 
        Then we have an isomorphism of $G^F$-representations \[\Pic(X(w))_{\Ql}\cong \Pic(X(vF(s)))_{\Ql}\oplus (i_{\varphi^{-1}(\supp_F(v))}/\1_{G^F}).\]
        \item \label{inductive_step_general_2}Suppose $vF(s)$ is $F$-Coxeter, then we have an isomorphism of $G^F$-representations \[\Pic(X(w))_{\Ql}\cong \Pic(\overline{X(vF(s))})_{\Ql}/(\1_{G^F}\oplus \bigoplus_{J\in S_F\setminus \{F\cdot s\}}i_{\varphi^{-1}(S_F\setminus \{J\})}).\]
    \end{rqenumerate}
\end{theorem}
\begin{proof}
    Set $t\coloneq F(s)$.
    We start with \cref{inductive_step_general_1}.
    As there is a universal homeomorphism $X(w)\to X(vtt)$, it is enough to show the statement for $X(vtt)$ (see \cref{univ_homeo_DL_general,genratchow}). 
    By assumption and the latter result, $\CH^1(X(vtt),1)_{\Ql}=0$ and this gives us a localization sequence 
    \[
        \begin{tikzcd}
        0 \arrow[r] & \CH^0(X(vt))_{\Ql} \arrow[r] & \Pic(X(vt\underline t))_{\Ql} \arrow[r] & \Pic(X(vtt))_{\Ql} \arrow[r] & 0.
        \end{tikzcd}
    \]
    As $vt$ is of full $F$-support, $X(vt)$ is irreducible (\cref{rqgen_iso_normal4} of \cref{gen_iso_normal}) and hence $\CH^0(X(vt))_{\Ql}\cong \1$.
    By \cref{affine_bundle} and because $\Ql[G^F]$ is semisimple, we get \begin{align*}
        \Pic(X(vtt))_{\Ql}&\cong \Pic(X( v t\underline t))_{\Ql}/\1\\
        &\cong \Pic(X(v \underline t))_{\Ql}/\1.
    \end{align*}
    The assumption on the vanishing of $\CH^1(X(vt),1)_{\Ql}$ gives us a localization sequence \[\begin{tikzcd}
    0 \arrow[r] & \CH^0(X(v))_{\Ql}  \arrow[r] & \Pic(X(v\underline{t}))_{\Ql} \arrow[r] & \Pic(X(vt))_{\Ql} \arrow[r] & 0.
    \end{tikzcd}\]
    As before, $X(v)$ has $|G^F/P^F_{\varphi^{-1}(\supp_F(v))}|$ irreducible components and so $\CH^0(X(v))_{\Ql}\cong i_{\varphi^{-1}(\supp_F(v))}$.
    We now have 
    \begin{align*}
        \Pic(X(vtt))_{\Ql}&\cong \Pic(X(v \underline t))_{\Ql}/\1\\
        &\cong (\Pic(X(vt))_{\Ql}\oplus i_{\varphi^{-1}(\supp_F(v))})/\1\\
        &=\Pic(X(vt))_{\Ql}\oplus (i_{\varphi^{-1}(\supp_F(v))}/\1).
    \end{align*}
    Now we continue with \cref{inductive_step_general_2}.
    By \cref{cycle_map_lemma_app,H1_vanishing_other_gps} we get as before
    \begin{align*}
        \Pic(X(w))_{\Ql}\cong \Pic(X(v\underline{t}))_{\Ql}/\1.
    \end{align*}
    Now, $X(v\underline{t})\subset \overline{X(vt)}$ is an open embedding where both schemes are smooth (see \cref{partial_compactification,smooth_cox_comp}). 
    Using \cref{helpfullemma} we see that $H^1(X(v\underline{t}))_\pr=0$ by Poincar\'e duality. 
    Hence we get a localization sequence 
    \[\begin{tikzcd}
    0 \arrow[r] & \CH^0(\overline{X(vt)}\setminus X(v\underline{t}))_{\Ql}  \arrow[r] & \Pic(\overline{X(vt)})_{\Ql} \arrow[r] & \Pic(X(v\underline{t}))_{\Ql} \arrow[r] & 0.
    \end{tikzcd}\]
    We will be done once we compute the first term.
    For $r\coloneq l(vt)$ we can write $vt=s_{i_1}\cdots s_{i_{r-1}}t$ and see \begin{align*}
        \overline{X(vt)}\setminus X(v\underline{t})=\bigcup_{j=1}^{r-1}\overline{X(s_{i_1}\cdots \widehat{s_{i_j}}\cdots s_{i_{_{r-1}}}t)}.
    \end{align*}
    To shorten notation, set $vt^{\wedge_j}\coloneq s_{i_1}\cdots \widehat{s_{i_j}}\cdots s_{i_{r-1}}t$ for $1\leq j\leq r-1$. 
    With this decomposition, we see \begin{align*}
        \CH^0(\overline{X(vt)}\setminus X(v\underline{t}))_{\Ql}&\cong \bigoplus_{j=1}^{r-1}\CH^0(\overline{X(vt^{\wedge_j})})_{\Ql}\\
        &\cong \bigoplus_{j=1}^{r-1}i_{\varphi^{-1}(\supp_F(vt^{\wedge_j}))}\\
        &\cong \bigoplus_{j=1}^{r-1}i_{\varphi^{-1}(S_F\setminus \{F\cdot s_{i_j}\})}\\
        &\cong \bigoplus_{J\in S_F\setminus \{F\cdot s\}}i_{\varphi^{-1}(S_F\setminus \{J\})}
    \end{align*}where we have used that $vt$ is $F$-Coxeter in the last two steps. 
\end{proof}

\begin{rem}\label{Ftriv_simpl}
    We may simplify \cref{inductive_step_general_2} of \cref{inductive_step_general} in the following two situations. 
    First, if $F$ acts trivially on $W$ then $\varphi\colon S\to S_F$ becomes a bijection and hence \begin{align*}
        \bigoplus_{J\in S_F\setminus \{F\cdot s\}}i_{\varphi^{-1}(S_F\setminus \{J\})}=\bigoplus_{u\in S\setminus \{s\}}i_{S\setminus \{u\}}.
    \end{align*}
    Second, \cref{conj_Hansen_BF_allcodim_3} of \cref{conj_Hansen_BF_allcodim} tells us that we can easily compute the principal series part of $\Pic(\overbar{X}(w))_{\Ql}$ for $w\in W$ an $F$-Coxeter element.
    In particular, for any $s\in S$ we have
    \begin{align*}
        \Pic(\overline{X(w)})_{\Ql,\pr}/(\1_{G^F}\oplus \bigoplus_{J\in S_F\setminus \{F\cdot s\}}i_{\varphi^{-1}(S_F\setminus \{J\})})&\cong i_{\varphi^{-1}(S_F\setminus \{F\cdot s\})}/(H^1(X(w))_\pr\oplus \1_{G^F}).
    \end{align*}
    Moreover, $\Pic(\overbar{X}(w))_{\Ql,\pr}=\Pic(\overbar{X}(w))_{\Ql}$ if $\Pic(X(w))_{\Ql}=0$ by \cref{Pic_PS}.
    Note that the $\ell$-adic cohomology of $X(w)$ has been computed in \cite[Section 7.3]{Lusztig1976} (the actual irreducible representation $H^1(X(w))_\pr$ can be derived from \cite[Corollary 2.10]{Lusztig1976} in particular cases, as we will see in some examples in \cref{small_rk_comps}).
\end{rem}

\begin{rem}
    Let $s\in S$ and $v\in W$ be as in \cref{inductive_step_general}.
    \Cref{inductive_step_general_2} deals with the case where $vF(s)$ is a particular type of terminal element; namely an $F$-Coxeter element.  
    If $vF(s)$ is some other terminal element $\varpi$, \cref{inductive_step_general} can only be applied assuming some vanishing results.  
\end{rem}

In untwisted type $A$, the results in this section are enough to compute $\Pic(X(w))_{\Ql}$ for all $w\in W$ of full support, due to the fact that the terminal elements are precisely the Coxeter elements (see \cite[Example 3.1.16]{Geck2000CharactersOF}). 
We elaborate on this in \cref{section_type_A}. 
In \cref{small_rk_comps}, we show how these results can be used to compute $\Pic(X(w))_{\Ql}$ for $w\in W$ of full support in the cases of $\prescript{2}{}{A_3}, B_2, B_3$ and $G_2$.
However, we only get partial results for the last two, due to the existence of terminal elements that are non-Coxeter and unequal to $w_0$.

\begin{rem}    
    As mentioned in \cref{sec:cycle_class_map}, we will see many examples where the cycle map $\Pic(X(w))_{\Ql,\upr}\hookrightarrow H^2(X(w))_\pr$ is an isomorphism.
    However, this is not always the case as the following example illustrates.
    Let $G^F$ be of type $\prescript{2}{}{A_2}$, $\prescript{2}{}{B_2}$ or $\prescript{2}{}{G_2}$ and let $w\in W$ be of length $2$.
    In particular, $w$ is of full $F$-support and not $F$-Coxeter. 
    \Cref{inductive_step_general_2} of \cref{inductive_step_general} gives us that $\Pic(X(w))_{\Ql,\upr}=0$ as for any $s\in S$ there are isomorphisms of $G^F$-representations $\Pic(\overline{X(s)})_{\Ql,\upr}\cong H^2(\overline{X(s)})_\pr\cong \1$ by \cref{DMR_cycle_iso_on_BF_ver2}.
    On the other hand, there are isomorphisms of $G^F$-representations $H^2(X(w))_\pr\cong H^2_c(X(w))_\pr\cong \St$ by \cite[Proposition 3.3.15]{digne2006cohomologie}.
    The result \cite[Proposition 3.4.2]{digne2006cohomologie} guarantees that there is a natural number $n\neq 1$ such that $H^{2,n}_\mot(X(w))_{\Ql,\pr}\cong \St$.
    As $\dim X(w)=2$ and $H^{2,n}_\mot(X(w))\cong \CH^n(X(w),2n-2)$, we know that $n=2$.
\end{rem}

\subsection{The untwisted type \texorpdfstring{$A$}{A} case}\label{section_type_A}

Throughout this entire section, we restrict ourselves to $G^F$ of untwisted type $A_{n-1}$ where $n\geq 2$. 
We start by stating \cref{cycle_map_iso}, which is the assertion of \cite[Proposition 7.15]{orlik2014cohomology} in codimension $1$. 
We will give another proof of this using \cref{Pic_PS} and \cref{orlik_codim_11}.

\begin{proposition}\label{cycle_map_iso}
Let $w\in F^+$.
The cycle map $\Pic(\overbar{X}(w))_{\Ql}\to H^2(\overbar{X}(w))$ is an isomorphism of $G^F$-representations. 
\end{proposition}

The following \cref{dim_cox_closure} can also be proved using the blow up formula for Chow groups (\cite[Proposition 6.7]{fulton2012intersection}); the Zariski closure of the standard Coxeter variety can be obtained as a sequence of blow ups starting from $\Pro^{n-1}$ (\cite[Proposition 1.6]{Genestier_espaces_sym_Drin}, see also \cite[Proof of Proposition 3.4]{orlik2014cohomology}).

\begin{proposition}\label{dim_cox_closure}
Let $w\in W$ be a Coxeter element, then $\cl\colon \Pic(\overline{X(w)})_{\Ql}\to H^2(\overline{X(w)})$ is an isomorphism of $G^F$-representations and 
\begin{align*}
\Pic(\overline{X(w)})_{\Ql}\cong \1\oplus \bigoplus_{i=1}^{n-2}i_{S\setminus \{s_i\}}.
\end{align*}
\end{proposition}
\begin{proof}
First, we have $\overline{X(w)}\cong \overbar{X}(w)$ by \cref{smooth_cox_comp}.  
To show that the cycle map is an isomorphism, it is enough to show $\Pic(\overline{X(w)})_{\Ql,\upr}=\Pic(\overline{X(w)})_{\Ql}$ and $H^2(\overline{X(w)})_\upr=H^2(\overline{X(w)})$ by \cref{DMR_cycle_iso_on_BF_ver2}. 
The $\ell$-adic cohomology group has been computed in \cite[Proposition 3.4]{orlik2014cohomology} as follows.
Set $P_1=G$ and $P_i=P_{S\setminus \{s_{i-1}\}}$ for $1\leq i\leq n-1$, then 
\begin{align*}
H^2(\overline{X(w)})&\cong \bigoplus_{i=1}^{n-1}\Ind^{G^F}_{P_i^F}(\1)
\end{align*}and we see that this representation is contained in the unipotent principal series.
To see $\Pic(\overline{X(w)})_{\Ql,\upr}=\Pic(\overline{X(w)})_{\Ql}$ we apply \cref{Cox_pic,Pic_PS}.
The claim follows.
\end{proof}

\subsubsection{\texorpdfstring{$\ell$}{l}-adic cohomology}\label{subsec_cohomology}

To invoke \cref{inductive_step_general}, we need to determine when the vanishing requirements in \cref{inductive_step_general_1} are satisfied. 
\Cref{orlik_codim_11} by Orlik tells us when $H^1(X(w))$ vanishes for nonterminal $w\in W$. 
To add details, we allow ourselves to correct the statement and update its proof. 
Orlik uses the notion of \textit{height} of elements in $F^+$ and $W$ \cite[Section 5]{orlik2014cohomology}.
Instead of asking the reader to look up the definition, we give a shorter and equivalent definition below.

\begin{definition}
    Let $w\in F^+$.
    We define the \textit{height} $\htt(w)$ of $w$ as 
    \begin{align*}
        \htt(w)\coloneq\dfrac{l(w)-l(\gamma(w)_{\minn})}{2}
    \end{align*}
    where $\gamma\colon F^+\to W$ is the natural projection and $\gamma(w)_{\minn}$ is defined as in \cref{GKP_thm}.
    After replacing $\gamma$ with the identity, we can make the same definition for elements in $W$.
\end{definition}

\begin{rem}
    The height takes values in $\N$.
    Indeed, this is a priori not clear from our definition, but can be shown using \cite[Section 5, Lemma 5.4]{orlik2014cohomology}.
    In type $A$, we can make the following reinterpretations of the height using \cite[Example 3.1.16]{Geck2000CharactersOF}.
    For $w=s_{i_1}\cdots s_{i_r}\in F^+$, the statement $\htt(w)\geq 1$ is equivalent to the existence of a $s\in S$ and $1\leq j\neq k\leq r$ such that $s=s_{i_j}=s_{i_k}$.
    Similarly, $\htt(w)=0$ is equivalent to the $s_{i_j}$ being pairwise distinct.
    We have the same reinterpretations in $W$ if we require $w=s_{i_1}\cdots s_{i_r}\in W$ to be a reduced decomposition.
\end{rem}

For a subset $I\subset S$, we say $W_I$ is \textit{irreducible} if the underlying Dynkin diagram is connected.

\begin{proposition}[Orlik]\cite[Proposition 7.6]{orlik2014cohomology}\label{orlik_codim_11}
Let $w\in F^+$ with $\htt(w)\geq 1$.
If $W_{\supp(w)}$ is irreducible then $H^{2l(w)-1}_c(X(w))=0$. 
\end{proposition}
\begin{proof}
    Set $J\coloneq \supp(w)$.
    Note that $L_J$ is isomorphic to a product of $\Gl_d$ for some $2\leq d\leq n$ with $n-d$ copies of $\G_m$.
    This means that the full flag variety of $L_J$ is isomorphic to the full flag variety of $\Gl_d$, and hence we know that $X_{L_J}(w)$ is a classical Deligne--Lusztig variety in the full flag variety of $\Gl_d$ (as opposed to a nontrivial product of such in different full flag varieties, see \cite[Proposition 3.6]{orlik2014cohomology}). 
    As $H^{2l(w)-1}_c(X(w))\cong \Rr^{G}_{L_J}H^{2l(w)-1}_c(X_{\Gl_d}(w))$, we may assume without loss of generality that $J=S$.
    We prove the claim by induction on the length of $w$.
    Consider the set \begin{align*}
    A\coloneq \{w\in F^+:\htt(w)\geq 1, \supp(w)=S\}.
    \end{align*}
    Set $m\coloneq \min\{l(w):w\in A\}=|S|+1$. 
    The start of induction concerns all elements $w\in A$ such that $l(w)=m$. 
    This is covered by \cref{H1_vanishing_other_gps} (see also \cite[Corollary 5.9]{orlik2014cohomology}). 
    For the induction step, suppose that for all $w\in A$ with $l(w)=k\geq m$ we have $H^{2l(w)-1}_c(X(w))=0$. 
    Let $w\in A$ with $l(w)=k+1$.
    Without loss of generality, we may assume $w=sw's$ with $w'\in F^+$ and $s\in S$ (see for example \cite[Proposition 3.1.6]{digne2006cohomologie}, \cite[Proposition 7.9]{orlik2014cohomology}).
    By the induction hypothesis and Poincar\'e duality, $H^1(X(w's))=0$. 
    \Cref{H1_vanishing_lower_impl_higher} now gives us $H^1(X(w))=0$ and hence we are done by Poincar\'e duality. 
\end{proof}

\begin{rem}\label{Saschas_fouten}
\cref{orlik_codim_11} is \cite[Proposition 7.6]{orlik2014cohomology} without the assumption on the support.  
Without this assumption, the statement is not true, as illustrated with a counterexample in \cref{not_full_supp}.
This causes the crucial step of the proof of \cite[Corollary 7.8]{orlik2014cohomology} to be false unless an extra assumption is added.
Such a sufficient condition is added in \cref{Sascha_fixed}.  
\end{rem}

\begin{exmp}\label{not_full_supp}
Let $G^F$ be of type $A_4$ and $w=s_1s_2s_1s_4$.
Note that $\supp(w)=\{s_1,s_2,s_4\}$ corresponds to a disconnected Dynkin diagram (i.e., $W_{\supp(w)}$ is reducible). 
By \cite[Proposition 3.6]{orlik2014cohomology}, \begin{align*}
H^7_c(X(w))=\Rr^{G}_{L_{\supp(w)}}H^7_c(X_{\Gl_3}(s_1s_2s_1)\times X_{\Gl_2}(s_4)). 
\end{align*} 
By the K\"unneth formula, \begin{align*}
0\neq H^6_c(X_{\Gl_3}(s_1s_2s_1))\otimes_{\Ql} H^1_c(X_{\Gl_2}(s_4))\subset H^7_c(X_{\Gl_3}(s_1s_2s_1)\times X_{\Gl_2}(s_4)).
\end{align*}
Taking the Harish-Chandra induction $\Rr^{G}_{L_{\supp(w)}}$ on both sides gives us a nonzero subrepresentation of $H^7_c(X(w))$. 
\end{exmp}

\begin{cor}[Orlik]\label{Sascha_fixed}
Let $w=sw's\in F^+$ with $\htt(w')\geq 1$ and $\supp(w)=S$.
Suppose $H^1(X(w'))=0$ (by \cref{orlik_codim_11}, a sufficient condition for this is that $W_{\supp(w)}$ is irreducible). 
Then \begin{align*}
H^{2l(w)-2}_c(X(w))=H^{2l(w's)-2}_c(X(w's))(-1)\oplus (i_{\supp(w')}/\1)(-l(w)+1).
\end{align*}
\end{cor}
\begin{proof}
See \cite[Proof of Corollary 7.8]{orlik2014cohomology}. 
The vanishing assumption now guarantees $H^{2l(w')-1}_c(X(w'))=0$ by Poincar\' e duality. 
\end{proof}

\subsubsection{Picard group}

We obtain the following reformulation of \cref{inductive_step_general} in the untwisted type $A$ case.

\begin{theorem}\label{rat_thm}
Let $w'\in F^+$ and $s\in S$ such that $w\coloneq sw's$ has full support.  
There is an isomorphism of $G^F$-representations
\begin{align*}
    \Pic(X(w))_{\Ql}\cong \begin{cases}
    \Pic(X(w's))_{\Ql}\oplus (i_{\text{supp}(w')}/\1)&\quad\text{ if }\text{ht}(w's)\geq 1\\
    i_{S\setminus \{s\}}/i_{S\setminus \{s_{n-1}\}}&\quad\text{ else.}
    \end{cases}
\end{align*}
\end{theorem}
\begin{proof}
    The case $\htt(w's)\geq 1$ follows from \cref{inductive_step_general_1} of \cref{inductive_step_general} and \cref{orlik_codim_11}.
    The other case follows from \cref{inductive_step_general_2} and \cref{dim_cox_closure} (see also \cref{Ftriv_simpl}), as \begin{align*}
        \Pic(\overline{X(w's)})_{\Ql}/(\bigoplus_{u\in S\setminus \{s\}}i_{S\setminus \{u\}}\oplus\1)&\cong \bigoplus_{i=1}^{n-2}i_{S\setminus \{s_i\}}/(\bigoplus_{u\in S\setminus \{s\}}i_{S\setminus \{u\}})\\
        &\cong \bigoplus_{u\in S\setminus \{s_{n-1}\}}i_{S\setminus \{u\}}/(\bigoplus_{u\in S\setminus \{s\}}i_{S\setminus \{u\}})\\
        &\cong i_{S\setminus \{s\}}/i_{S\setminus \{s_{n-1}\}}.
    \end{align*}
\end{proof}

\cref{rat_thm} allows us to determine any $\Pic(X(w))_{\Ql}$ where $w\in W$ is of full support and $\htt(w)\geq 1$ recursively, using the following algorithm. 
Recall that we are not interested in Coxeter elements due to \cref{Cox_pic}. 

\begin{algorithm}\label{algorithm}
Let $w\in W$ of full support and $\htt(w)\geq 1$. 

\begin{rqenumerate}[label=\textbf{A.\arabic*}]
  \item \label{l} \underline{If $w\xrightarrow{s,>}sws$ for some $s\in S$.} 
	Then \cref{rat_thm} allows us to 			compute $\Pic(X(w))_{\Ql}$ via a case 			distinction. 
	Set $w'\coloneq sws$. 
  \begin{rqenumerate}[label=\textbf{A.1.\arabic*}]
    \item \underline{If $\htt(w's)\geq 1$.}Then \[\Pic(X(w))_{\Ql}\cong \Pic(X(w's))_{\Ql}\oplus (i_{\text{supp}(w')}/\1).\]
    Note $|\supp(w')|\geq |S|-1$, and so $i_{\text{supp}(w')}$ is easily determined by \cite[Theorem 2.1]{orlik2014cohomology}.
    To compute $\Pic(X(w's))_{\Q}$, return to the beginning of the algorithm by taking $w\coloneq w's$. 
    \item \underline{If $\htt(w's)=0$.} Then \[\Pic(X(w))_{\Ql}\cong i_{S\setminus \{s\}}/i_{S\setminus \{s_{n-1}\}}.\]
    The algorithm stops here. 
  \end{rqenumerate}
  \item \label{ll}\underline{If $w\xrightarrow{s,=}sws$ for all $s\in S$.} 
  Then we claim that there is an element $v\in \Cyc(w)$ such that $v\xrightarrow{s,>}svs$ for some $s\in S$.
  If this claim is true then $\Pic(X(w))_{\Ql}\cong \Pic(X(v))_{\Ql}$ by \cref{cyclic_shift_iso} and we may restart the algorithm by taking $w\coloneq v$ as $v$ has full support and $\htt(v)=\htt(w)\geq 1$.
  Note that in particular, one lands in \cref{l}.

  To show the claim, let $C(w)\subset W$ be the conjugacy class of $w$.
  Then by \cref{GKP_thm} there exists a $w_{\minn}\in C(w)_{\minn}$ such that $w\xrightarrow{}w_{\minn}$. 
Let $s_1,s_2,\ldots,s_r\in S$ and $w=w_0,w_1,\ldots,w_r=w_{\minn}$ be such a pair of sequences. 
As $\htt(w_{\minn})=0$, we have $l(w)>l(w_{\minn})$.
This means that the set \[\{1\leq i\leq r:w_i\xrightarrow{s_{i+1},>}w_{i+1}\}\]is nonempty and contains a minimum $i_{\minn}$.
The element $v\coloneq w_{i_{\minn}}$ satisfies the claim. 
\end{rqenumerate}
This algorithm ends because each time we land in \cref{l}, we are either done or calculating the rational Picard group of a Deligne--Lusztig variety whose corresponding Weyl group element has colength $1$ with respect to the element we started with. 
\end{algorithm}

We have applied \cref{algorithm} to compute all rational Picard groups for $w\in W$ of full support for types $A_3$ and $A_4$ (note that elements of full support in $W(A_1)=S_2$ or $W(A_2)=S_3$ are Coxeter or $w_0$, and hence already covered by Hansen results) in \cref{subsec_gln}.
We now give another proof of \cref{cycle_map_iso}.
\begin{proof}[Proof of \cref{cycle_map_iso}]\prooflabell{cycle_map_iso}{Pf_cycle_map_iso}
    We may assume that $w$ is of full support by \cref{subsec:elsofpropsupp}.
    As we have shown the Coxeter case in \cref{dim_cox_closure}, we can assume $w$ is of height at least $1$. 
    By \cref{DMR_cycle_iso_on_BF_ver2}, it is enough to show that $\Pic(\overbar{X}(w))_{\Ql}$ and $H^2(\overbar{X}(w))$ are entirely contained in the unipotent principal series.
    As all unipotent representations lie in the (unipotent) principal series (we are in the untwisted type $A$ case, see \cite[Section 13.8]{Carter1972}), $H^2(\overbar{X}(w))$ lies in the unipotent principal series by \cite[Proposition 4.2]{Lusztig1998}.
    We make a case distinction to show that $\Pic(\overbar{X}(w))_{\Ql}$ is in the unipotent principal series.
    If $w\in F^+$ is of full support and of height at least $1$ then by \cref{cycle_map_lemma_app,orlik_codim_11} we get an exact sequence 
    \[
        \begin{tikzcd}
        0 \arrow[r] & \CH^0(\partial \overbar{X}(w))_{\Ql} \arrow[r] & \Pic(\overbar{X}(w))_{\Ql} \arrow[r] & \Pic(X(w))_{\Ql} \arrow[r] & 0.
        \end{tikzcd}
    \]
    The representation $\CH^0(\partial \overbar{X}(w))_{\Ql}$ is in the unipotent principal series, as it is a direct sum of Harish-Chandra inductions of the trivial representation. 
    By \cref{algorithm}, the representation $\Pic(X(w))_{\Ql}$ is also in the unipotent principal series. 
    Hence, $\Pic(\overbar{X}(w))_{\Ql}$ is in the unipotent principal series.
\end{proof}

\Cref{algorithm} allows us to state the following more precise version of \cref{rq_DMR_cycle_iso_on_BF_3} in this untwisted type $A$ case. 

\begin{proposition}\label{cycle_inj_GLn}
    For $w\in W$, the cycle map $\cl\colon \Pic(X(w))_{\Ql}\to H^2(X(w))$ is an injection of $G^F$-representations.
\end{proposition}
\begin{proof}
    By \cref{rq_DMR_cycle_iso_on_BF_3}, we need to show that both $\Pic(X(w))_{\Ql}$ and $H^2(X(w))$ lie in the unipotent principal series. 
    This is the case for the former by \cref{algorithm}. 
    For the latter, this follows from the fact that any unipotent representation is in the principal series (\cite[Section 13.8]{Carter1972}).
\end{proof}

There are many cases where the cycle map in \cref{cycle_inj_GLn} is actually an isomorphism.
More specifically, if $G^F$ is of type $A_n$ with $1\leq n\leq 4$ then we can show this for Weyl group elements of height at least $1$ and of full support.
Indeed, the cases $n=1$ or $2$ are vacuous/obvious and for $n=3$ or $4$ see \cref{GL_4_rat,GL_5_rat_prop}.
The author does not know of any examples in untwisted type $A$ where this cycle map is not surjective.
To prove these aforementioned results, we use the following.

\begin{theorem}\label{iso_crit}
Let $w\in W$ be of full support. 
Then the cycle map $\cl\colon \Pic(X(w))_{\Ql}\to H^{2}(X(w))$ is an isomorphism if either of the following conditions are satisfied\begin{rqenumerate}\label{rqiso_crit}
\item \label{rqiso_crit1}The element $w$ has height $1$,
\item \label{rqiso_crit2}The element $w$ is of the form $w=sw's$ with $w'\in W$, $s\in S$ and $l(w)=l(w')+2$ such that $\htt(w')\geq 1$, $\Pic(X(w's))_{\Ql} \cong H^2(X(w's))$ and $H^1(X(w'))=0$.
\end{rqenumerate}
\end{theorem}
\begin{proof}
\Cref{cycle_inj_GLn} tells us that $\cl\colon \Pic(X(w))_{\Ql}\to H^2(X(w))$ is injective, so we only need to show surjectivity. 
This comes down to showing that $\dim_{\Ql}\Pic(X(w))_{\Ql}=\dim_{\Ql}H^2(X(w))$.
Suppose $w\in W$ with $\htt(w)=1$.
Without loss of generality, we may assume $w=sw's$ for some $w'\in W$, $s\in S$ such that $l(w)=l(w')+2$ (see \cref{ll} of \cref{algorithm}, as $\ell$-adic cohomology is also invariant under cyclic shift by \cref{cyclic_shift_iso}).
The element $w'\in W$ is a Coxeter element in $W$ or $W_{S\setminus \{s\}}$. 
In the first case, we know that $H^2(X(w))=0$ by \cite[Corollary 5.13]{orlik2014cohomology} and Poincar\'e duality, which means that the claim is trivially true. 
In the second case, note $l(w)=n$ and $l(w')=n-2$.
Then by \cref{rat_thm}
\begin{align*}
    \Pic(X(w))_{\Ql}=i_{S\setminus \{s\}}/i_{S\setminus \{s_{n-1}\}}.
\end{align*}
On the other hand, by \cite[Corollary 5.9]{orlik2014cohomology} we get the following isomorphisms of $G^F$-representations
\begin{align*}
    H^{2l(w)-2}_c(X(w))&\cong H^{2(n-1)}_c(X(w))\\
    &\cong H^{2n-4}_c(X(w'))/(j_{(n-1,1)}\oplus\1)\\
    &\cong H^{2l(w')}_c(X(w'))/i_{S\setminus \{s_{n-1}\}}\\
    &\cong i_{\text{supp}(w')}/i_{S\setminus \{s_{n-1}\}}.
\end{align*}
As $\text{supp}(w')=S\setminus \{s\}$ and $\dim_{\Ql}H^{2l(w)-2}_c(X(w))=\dim_{\Ql}H^2(X(w))$ by Poincar\'e duality, this proves that if \cref{rqiso_crit1} holds then $\text{cl}$ is an isomorphism.
The fact that \cref{rqiso_crit2} implies that $\cl$ is an isomorphism follows from \cref{rat_thm,Sascha_fixed}. 
\end{proof}

\section{Appendix}\label{Appendix}

\subsection{Dynkin diagrams}

The connected Dynkin diagrams with Frobenius action are displayed below.
The left column displays the trivial action ($\delta=1$, untwisted) and the right column the action of order $\delta=2$ (twisted).
In the right column, the arrows display how the vertices are permuted by $F$, and vertices
without an incoming or outgoing arrow are understood to be fixed.
In type $\prescript{2}{}{A_n}$ with $n$ odd,
$s_{\frac{n+1}{2}}$ is the unique vertex that is fixed by $F$.
If $n$ is even then no vertex is fixed.

\begin{minipage}[t]{0.48\textwidth}
\[
\begin{aligned}
	A_n (n\geq 1)&:\dynkin[
		make indefinite edge/.list={2-3},
		labels={1,2,n},
		label macro/.code={s_{#1}},
		edge length=.75cm
	]A{3}\\
	B_n (n\geq 2)&:\dynkin[
		make indefinite edge/.list={2-3},
		labels={1,2,n-1,n},
		arrows=true,
		label macro/.code={s_{#1}},
		edge length=.75cm
	]B{4}\\
	C_n (n\geq 3)&:\dynkin[
		make indefinite edge/.list={2-3},
		labels={1,2,n-1,n},
		arrows=true,
		label macro/.code={s_{#1}},
		edge length=.75cm
	]C{4}\\
	D_n (n\geq 4)&:\dynkin[
		make indefinite edge/.list={2-3},
		labels={1,2,n-2,n-1,n},
		label directions={,,right,,},
		label macro/.code={s_{#1}},
		edge length=.75cm
	]D{5}\\
	E_6&:\dynkin[
		labels={1,2,3,4,5,6},
		label macro/.code={s_{#1}},
		edge length=.75cm
	]E{6}\\
	E_7&:\dynkin[
		labels={1,2,3,4,5,6,7},
		label macro/.code={s_{#1}},
		edge length=.75cm
	]E{7}\\
	E_8&:\dynkin[
		labels={1,2,3,4,5,6,7,8},
		label macro/.code={s_{#1}},
		edge length=.75cm
	]E{8}\\
	F_4&:\dynkin[
		labels={1,2,3,4},
		arrows=true,
		label macro/.code={s_{#1}},
		edge length=.75cm
	]F{4}\\
	G_2&:\dynkin[
		labels={1,2},
		arrows=true,
		label macro/.code={s_{#1}},
		edge length=.75cm
	]G{2}
\end{aligned}
\]

\end{minipage}
\hfill
\begin{minipage}[t]{0.48\textwidth}

\[
\begin{aligned}
\prescript{2}{}{A_n}&:\dynkin[
	make indefinite edge/.list={2-3},
	label macro/.code={s_{#1}},
	labels*={s_1,s_2,s_{n-1},s_n},
	involutions={1<below>4;23},
	edge length=.75cm
]A{4}\\
\prescript{2}{}{D_n}&:\dynkin[
	make indefinite edge/.list={2-3},
	labels={1,2,n-2,n,n-1},
	label directions={,,right,left,left},
	involutions={[out=80,in=100,relative]45},
	label macro/.code={s_{#1}},
	edge length=.75cm
]D{5}\\
\prescript{2}{}{E_6}&:\dynkin[
	labels*={s_1,s_2,s_3,s_4,s_5,s_6},
	involutions={16;35},
	label macro/.code={s_{#1}},
	edge length=.75cm
]E{6}\\
\prescript{2}{}{B_2}&:\dynkin[
	labels={1,2},
	arrows=true,
	involutions={12},
	label macro/.code={s_{#1}},
	edge length=.75cm
]B{2}\\
\prescript{2}{}{G_2}&:\dynkin[
	labels={1,2},
	arrows=true,
	involutions={12},
	label macro/.code={s_{#1}},
	edge length=.75cm
]G{2}\\
\prescript{2}{}{F_4}&:\dynkin[
	labels*={s_1,s_2,s_3,s_4},
	arrows=true,
	involutions={14;23},
	label macro/.code={s_{#1}},
	edge length=.75cm
]F{4}
\end{aligned}
\]

\end{minipage}

The following is the unique Dynkin diagram where $F$ acts with order $3$. 

\[\prescript{3}{}{D_4}:\begin{tikzcd}
	{\bullet_{s_1}} \\
	& {\bullet_{s_2}} & {\bullet_{s_4}} \\
	{\bullet_{s_3}}
	\arrow[no head, from=1-1, to=2-2]
	\arrow[curve={height=-18pt}, from=1-1, to=2-3]
	\arrow[no head, from=2-2, to=2-3]
	\arrow[curve={height=-18pt}, from=2-3, to=3-1]
	\arrow[curve={height=-12pt}, from=3-1, to=1-1]
	\arrow[no head, from=3-1, to=2-2]
\end{tikzcd}\]

\subsection{Computations in groups of small rank}\label{small_rk_comps}

In this section, we compute $\Pic(X(w))_{\Ql}$ for $w\in W$ of full $F$-support in examples of $G^F$ of small rank. 
One indeed finds many examples where the rational Picard group does not vanish. 
The groups considered in \cref{subsec_gln,small_rk_comps:Sp} are all of classical type, hence we have the vanishing in \cref{Cox_pic} at our disposal.
In \cref{subsec:typeG2} we are dealing with an exceptional group, but at least know that the principal series part of the rational Picard group of a Coxeter variety vanishes (\cref{conj_Hansen_BF_allcodim_1} of \cref{conj_Hansen_BF_allcodim}).
The computations all use the same tactic by combining the results \cref{H1_vanishing_lower_impl_higher,H1_vanishing_other_gps,inductive_step_general} (see also \cref{Ftriv_simpl}) and the following computations. 
\begin{itemize}
    \item For $J\subset S$, the decomposition $i_J$ into irreducible $G^F$-representations can be computed using GAP with the functions $\mathsf{CoxeterGroup}$, $\mathsf{Twistings}$, $\mathsf{ReflectionSubgroup}$ and $\mathsf{InductionTable}$.
    The reader is warned that GAP uses Dynkin diagram conventions different from ours: the double edge in types $B$ and $C$ is at the beginning and the long and short roots are swapped (also for type $G_2$).
    \item For all the types that we consider except $G_2$, the $G^F$-representation $H^1(X(c))_\pr$ (for $c\in W$ an $F$-Coxeter element) can be computed using \cite[Corollary 2.10]{Lusztig1976}, \cite[Proposition 3.3.14]{digne2006cohomologie} and the previous computation. 
    For type $G_2$, we additionally use \cite[Théorème 4.4.3]{digne2006cohomologie}.
\end{itemize}

For this reason, we only give a full proof of \cref{GL_4_rat}.

\subsubsection{Types \texorpdfstring{$A_3$}{A_3}, \texorpdfstring{$A_4$}{A_4} and \texorpdfstring{$\prescript{2}{}{A_3}$}{^2A_3}}\label{subsec_gln}

Let $n\geq 2$ and suppose $G^F$ is of type $A_{n-1}$ or $\prescript{2}{}{A_{n-1}}$. 
By \cite[Section 13.8]{Carter1972}, the unipotent $G^F$-representations are parametrized by partitions of $n$. 
For $\lambda$ such a partition, we denote by $j_\lambda$ the corresponding irreducible representation with the convention that $j_{(n)}=\1$ and $j_{(1^n)}=\St$. 
If $G^F$ is untwisted, we can decompose $i_J$ for $J\subset S$ in terms of $j_\lambda$ as explained in \cite[Theorem 3.2.1]{dat_orlik_rapoport_2010}. 

\begin{proposition}\label{GL_4_rat}
Consider $G^F$ of type $A_3$ and let $W'\subset W$ be the subset of elements of full support and height at least $1$.
For every $w\in W'$, the cycle map $\cl\colon \Pic(X(w))_{\Ql}\to H^2(X(w))$ is an isomorphism. 
In particular, for $w\in W'$,\begin{align*}
\Pic(X(w))_{\Ql}\cong \begin{cases}
j_{(2,2)}&\quad\text{if }w=(1\ 3)(2\ 4)\\
j_{(3,1)}&\quad\text{if }w=(1\ 4)\\
0&\quad\text{else}.
\end{cases}
\end{align*} 
\end{proposition}
\begin{proof}
    We start by computing the rational Picard groups. 
    Using GAP, we determine $|W'|=9$.
    By \cref{cyclic_shift_iso}, it is enough to compute the rational Picard group for the $5$ cyclic shift class representatives given in \cref{table:gl4}.
    \begin{table}[h!]
    \centering
    \begin{tabular}{ m{10em} | m{10em} | m{10em}} 
    \makecell{$w$} & \makecell{$l(w)$} & \makecell{$|\text{Cyc}(w)|$}\\ [0.5ex] 
    \hline\hline
    \makecell{$(1\ 2\ 4)$} & \makecell{$4$} &\makecell{$4$}\\
    \hline 
    \makecell{$(1\ 3)(2\ 4)$} & \makecell{$4$} &\makecell{$1$} \\
    \hline 
    \makecell{$(1\ 4\ 2\ 3)$} & \makecell{$5$} & \makecell{$2$}\\
    \hline 
    \makecell{$(1\ 4)$} & \makecell{$5$} & \makecell{$1$} \\
    \hline 
    \makecell{$(1\ 4)(2\ 3)$} & \makecell{$6$} & \makecell{$1$}
    \end{tabular}
    \caption{Five cyclic shift representatives in $W'$.}
    \label{table:gl4}
    \end{table}
    If $w=(1\ 2\ 4)=s_3s_1s_2s_3=sw's$ for $s\coloneq s_3$ and $w'\coloneq s_1s_2$, \cref{algorithm} gives $\Pic(X(w))_{\Ql}=0$ after one step. 
    If $w=(1\ 4\ 2\ 3)=s_2s_3s_2s_1s_2$, we get $\Pic(X(w))_{\Ql}\cong \Pic(X(s_3s_2s_1s_2))_{\Ql}$ and this vanishes as $s_3s_2s_1s_2\in \Cyc((1\ 2\ 4))$. 
    If $w=(1\ 4)(2\ 3)=w_0$ then \cref{rat_w0_vanishes} gives $\Pic(X(w_0))_{\Ql}=0$ (using $w=s_1s_2s_3s_1s_2s_1$ one can indeed check that \cref{algorithm} gives the same answer). 
    What is left to do is to consider $w_1\coloneq (1\ 3)(2\ 4)$ and $w_2\coloneq (1\ 4)$. 

    First consider $w_1=s_2s_1s_3s_2$, which is of height $1$.
    By \cref{iso_crit}, the cycle map is an isomorphism for $w_1$.
    By \cref{algorithm} and \cite[Theorem 3.2.1]{dat_orlik_rapoport_2010}, we obtain the following isomorphisms of $G^F$-representations:
    \begin{align*}
        \Pic(X(w_1))_{\Ql}&\cong i_{S\setminus \{s_2\}}/i_{S\setminus \{s_3\}}\\
        &\cong (j_{(2,2)}\oplus j_{(3,1)}\oplus \1)/(j_{(3,1)}\oplus \1)\\
        &\cong j_{(2,2)}.
    \end{align*}
    For $w_2$, we may choose two different reduced expressions; $w_2=s_1s_2s_3s_2s_1=s_3s_1s_2s_1s_3$ which give rise to different choices of $s$ and $w'$.
    Indeed, we can take $w'\coloneq s_2s_3s_2$ and $s\coloneq s_1$ or $w'\coloneq s_1s_2s_1$ and $s\coloneq s_3$. 
    We show that both choices give the same answer. 
    For either choice, $H^1(X(w'))=0$ by \cref{orlik_codim_11} and $\Pic(X(w's))_{\Ql}=0=H^2(X(w))$ by the steps before, Poincar\'e duality and \cite[Appendix B]{orlik2014cohomology}.
    By \cref{iso_crit}, the cycle map is also an isomorphism for $w_2$. 
    Moreover, we have that $\text{ht}(w's)=1$ in either case and so \cref{rat_thm} and \cite[Theorem 3.2.1]{dat_orlik_rapoport_2010} gives us
    \begin{align*}
        \Pic(X(w_2))_{\Ql}&\cong i_{\text{supp}(w')}/\1\\
    &\cong j_{(3,1)}
    \end{align*} 
    independent of our choice of $w'$ and $s$.
    By \cite[Appendix B]{orlik2014cohomology} and Poincar\'e duality, we know $H^2(X(w))=0$ for all $w\in W'\setminus \{w_1,w_2\}$ which means that the corresponding injective cycle map from \cref{cycle_inj_GLn} is trivially surjective.
\end{proof}

The following \cref{GL_5_rat_prop} is obtained with the same tactic as \cref{GL_4_rat}, after decomposing the following inductions using GAP.
If $G^F$ is of type $A_4$ then 
\begin{align*}
    i_{S\setminus \{s_1\}}&\cong j_{(4,1)}\oplus \1\\
    &\cong i_{S\setminus \{s_4\}},\\
    i_{S\setminus \{s_2\}}&\cong j_{(3,2)}\oplus j_{(4,1)}\oplus \1\\
    &\cong i_{S\setminus \{s_3\}}.
\end{align*}

\begin{proposition}\label{GL_5_rat_prop}
Consider $G^F$ of type $A_4$, and let $w\in W$ of full support and non-Coxeter. 
Up to choice of cyclic shift representative, the representation $\Pic(X(w))_{\Ql}$ is given in \cref{table:Gl5}.
Moreover, the cycle map $\cl\colon \Pic(X(w))_{\Ql}\to H^2(X(w))$ is an isomorphism.
\end{proposition}

\begin{table}[h!]
    \centering
    \begin{tabular}{ m{8em} | m{5em} | m{5em} | m{15em}} 
    \makecell{$w$} & \makecell{$l(w)$} & \makecell{$|\Cyc(w)|$} & \makecell{$\Pic(X(w))_{\Ql}$} \\ [0.5ex] 
    \hline\hline
    \makecell{$(1\ 3)(2\ 4\ 5)$} & \makecell{$5$} & \makecell{$2$} & \makecell{$j_{(3,2)}$} \\
    \hline 
    \makecell{$(3\ 5)(1\ 4\ 2)$} & \makecell{$5$} & \makecell{$2$} & \makecell{$j_{(3,2)}$} \\
    \hline 
    \makecell{$(5\ 3\ 2\ 1)$} & \makecell{$5$} & \makecell{$12$} & \makecell{$0$}  \\
    \hline 
    \makecell{$(1\ 5\ 4\ 2\ 3)$} & \makecell{$6$} & \makecell{$10$} & \makecell{$0$} \\
    \hline 
    \makecell{$(1\ 3)(2\ 5)$} & \makecell{$6$} & \makecell{$1$} & \makecell{$j_{(3,2)}\oplus j_{(4,1)}$}   \\
    \hline 
    \makecell{$(1\ 4)(3\ 5)$} & \makecell{$6$} & \makecell{$1$} & \makecell{$j_{(3,2)}\oplus j_{(4,1)}$}   \\
    \hline 
    \makecell{$(1\ 5\ 2)$} & \makecell{$6$} & \makecell{$6$}& \makecell{$j_{(4,1)}$} \\
    \hline 
    \makecell{$(1\ 3\ 4)(2\ 5)$} & \makecell{$7$} & \makecell{$10$} & \makecell{$0$} \\
    \hline 
    \makecell{$(1\ 5)$} & \makecell{$7$} & \makecell{$1$} & \makecell{$j_{(4,1)}^2$}  \\
    \hline 
    \makecell{$(5\ 1\ 4\ 3)$} & \makecell{$7$} & \makecell{$2$} & \makecell{$j_{(4,1)}$} \\
    \hline 
    \makecell{$(1\ 5\ 2\ 3)$} & \makecell{$7$} & \makecell{$2$} & \makecell{$j_{(4,1)}$} \\
    \hline 
    \makecell{$(1\ 5\ 3\ 2\ 4)$} & \makecell{$8$} &\makecell{$6$} & \makecell{$0$} \\
    \hline 
    \makecell{$(1\ 5)(2\ 3)$} & \makecell{$8$} & \makecell{$1$} &\makecell{$j_{(4,1)}$} \\
    \hline
    \makecell{$(1\ 5)(3\ 4)$} & \makecell{$8$} & \makecell{$1$} &\makecell{$j_{(4,1)}$} \\
    \hline 
    \makecell{$(1\ 4)(2\ 5)$} & \makecell{$8$} & \makecell{$1$} &\makecell{$0$} \\
    \hline 
    \makecell{$l(w)>8$} & \makecell{$-$} & \makecell{$-$} &\makecell{$0$}
    \end{tabular}
    \caption{Rational Picard group of $X(w)$ for $G^F$ of untwisted type $A_4$, where $w$ is of full support and non-Coxeter.}
    \label{table:Gl5}
\end{table}

\begin{rem}
    The results above show that if $G^F$ is of type $A_n$ for $1\leq n\leq 3$ and $w\in W$ of full support, then $\Pic(X(w))_{\Ql}=0$ is equivalent to the smoothness of $\overline{X(w)}$.
    Indeed, this follows by combining \cref{gen_iso_normal,smooth_cox_comp,Cox_pic}, \cref{GL_4_rat} and the fact that in a full flag variety of type $A_3$ the only singular Schubert varieties are the ones associated to $(1\ 3)(2\ 4)$ and $(1\ 4)$ (recall the pattern avoidance criterion mentioned in \cref{rem_pattern_avoidance}).
    This equivalence is false in general.
    In particular, both implications already fail for $G^F$ of type $A_4$: by \cref{GL_5_rat_prop,gen_iso_normal} the Weyl group elements $(1\ 5\ 3)$ and $(1\ 3\ 4)(2\ 5)$ in $S_5$ establish this fact as the former corresponds to a smooth Schubert variety and the latter does not.
\end{rem}

We now turn our attention to the twisted case $\prescript{2}{}A_3$, where we use the tactic explained in \cref{small_rk_comps} after computing the following $G^F$-isomorphisms 
\begin{alignat*}{2}
    i_{\varphi^{-1}(F\cdot s_1)}&=i_{\{s_1,s_3\}}&&\cong j_{(3,1)}\oplus j_{(2,2)}\oplus \1,\\
    i_{\varphi^{-1}(F\cdot s_2)}&=i_{\{s_2\}}&&\cong j_{(2,1,1)}\oplus j_{(2,2)}\oplus \1,\\
    H^1(X(c))&\cong j_{(2,2)}.
\end{alignat*}

\begin{proposition}\label{type2A3}
    Let $G^F$ be of type $\prescript{2}{}{A_3}$ and $w\in W\setminus \{w_0\}$ of full support and not $F$-Coxeter. 
    Up to choice of cyclic shift representative, the representation $\Pic(X(w))_{\Ql}$ is given in \cref{table:2A3}.  
\end{proposition}

\begin{table}[h!]
    \centering
    \begin{tabular}{ m{8em} | m{5em} | m{5em} | m{15em}} 
    \makecell{$w$} & \makecell{$l(w)$} & \makecell{$|\Cyc(w)|$} & \makecell{$\Pic(X(w))_{\Ql}$} \\ [0.5ex] 
    \hline\hline
    \makecell{$s_2s_3s_2$} & \makecell{$3$} & \makecell{$4$} & \makecell{$j_{(3,1)}$} \\
    \hline 
    \makecell{$s_1s_2s_3$} & \makecell{$3$} & \makecell{$1$} & \makecell{$j_{(2,1,1)}$}  \\
    \hline 
    \makecell{$s_3s_2s_1$} & \makecell{$3$} & \makecell{$1$} & \makecell{$j_{(2,1,1)}$}  \\
    \hline 
    \makecell{$s_1s_3s_2s_3$} & \makecell{$4$} & \makecell{$2$} & \makecell{$j_{(3,1)}$} \\
    \hline 
    \makecell{$s_3s_2s_3s_1$} & \makecell{$4$} & \makecell{$2$} & \makecell{$j_{(3,1)}$} \\
    \hline
    \makecell{$s_2s_1s_3s_2$} & \makecell{$4$} & \makecell{$1$} & \makecell{$j_{(3,1)}^2\oplus j_{(2,2)}$}   \\
    \hline 
    \makecell{$s_2s_3s_2s_1s_2$} & \makecell{$5$} & \makecell{$1$}& \makecell{$j_{(3,1)}$}\\
    \hline 
    \makecell{$s_2s_1s_2s_3s_2$} & \makecell{$5$} & \makecell{$1$}& \makecell{$j_{(3,1)}$}\\
    \hline 
    \makecell{$s_1s_3s_2s_1s_3$} & \makecell{$5$} & \makecell{$1$}& \makecell{$j_{(3,1)}$}  
    \end{tabular}
    \caption{Rational Picard group of $X(w)$ for $G^F$ of type $\prescript{2}{}{A_3}$, where $w\neq w_0$ is of full support and not $F$-Coxeter.}
    \label{table:2A3}
\end{table}

\subsubsection{Types \texorpdfstring{$B_2$}{B_2} and \texorpdfstring{$B_3$}{B_3}}\label{small_rk_comps:Sp}

Let $n\geq 2$ and suppose $G^F$ is of type $B_n$.
By \cite[Proposition 11.4.2, Section 13.8]{Carter1972}, the irreducible $G^F$-representations that are in the unipotent principal series are parametrized by bipartitions of $n$.
For a bipartition $\lambda$, we denote by $j_\lambda$ the corresponding irreducible representation with the convention that $j_{(n;0)}=\1$ and $j_{(0;1^n)}=\St$. 
For type $B_2$, using the tactic outlined at the beginning of \cref{small_rk_comps} one gets the following $G^F$-isomorphisms
\begin{align*}
        i_{\{s_1\}}&\cong j_{(0;2)}\oplus j_{(1;1)}\oplus \1,\\
        i_{\{s_2\}}&\cong j_{(11;0)}\oplus j_{(1;1)}\oplus \1,\\
        H^1(X(c))&\cong j_{(1;1)}
\end{align*}
where $c\in W$ is a Coxeter element. 
Recall the unipotent cuspidal representation $\theta_{10}$ of a finite group of type $B_2$ (\cite[Section 13.7]{Carter1972}).
The same tactic as before gives us the following \cref{Pic_C2}.
To show the cycle map statement, we use \cref{rq_DMR_cycle_iso_on_BF_3} and the fact that the codomain is an irreducible $G^F$-representation \cite[Théorème 4.3.4.]{digne2006cohomologie}.

\begin{proposition}\label{Pic_C2}
    Let $G^F$ be of type $B_2$ and $w\in W\setminus \{w_0\}$ of full support and non-Coxeter (there are precisely $2$ such elements), then \begin{align*}
        \Pic(X(w))_{\Ql}\cong \begin{cases}
            j_{(11;0)}&\text{if }w=s_1s_2s_1\\
            j_{(0;2)}&\text{if }w=s_2s_1s_2.
        \end{cases}
    \end{align*}
    In both cases, the cycle map $\cl\colon\Pic(X(w))_{\Ql}\to H^2(X(w))_\pr\cong H^2(X(w))/\theta_{10}$ is an isomorphism. 
\end{proposition}

If $W$ is of type $B_3$, we have seen in \cref{sec:cyclic_shifts} that there is a terminal element $\varpi\in W\setminus \{w_0\}$ that is non-Coxeter.
To compute its rational Picard group (at least the principal series part), we need an additional vanishing result that is contained in \cref{H1_height0_nonCox}.
One can check (with GAP, for example) that $\Cyc(\varpi)\subset W$ has cardinality $4$ and is the subset of all such elements. 
From the following, we will only use the vanishing in degrees $8$ and $9$.

\begin{proposition}\label{H1_height0_nonCox}
    Let $G^F$ be of type $B_3$ and let $\varpi\in W\setminus \{w_0\}$ be terminal and non-Coxeter. 
    The principal series part of the cohomology of $X(\varpi)$ is given by \begin{align*}
        H^i_c(X(\varpi))_\pr\cong \begin{cases}
            j_{(0;111)}&\text{if }i=5 \\
            (j_{(0;21)}\oplus j_{(11;1)})(-2)&\text{if }i=6 \\
            (j_{(1;2)}\oplus j_{(21;0)})(-3)&\text{if }i=7 \\
            \1(-5)&\text{if }i=10 \\
            0&\text{else}.
        \end{cases}
    \end{align*} 
\end{proposition}
\begin{proof}
    Without loss of generality, we may assume $\varpi=s_1s_2s_3s_2s_3$. 
    To compute the corresponding cohomology groups, we use \cite[Corollary 7.20]{Digne_2014}.
    We write $\varpi=vw$ for $w\coloneq s_1s_2s_3s_2$ and $v\coloneq s_3$, and set $I\coloneq \{s_3\}$.
    Recall from \cref{rem:parab_covers} the parabolic Deligne--Lusztig variety cover $Y_I(\dot w)\to X_I(w)$.
    As $\prescript{w}{}I=I$, we get an isomorphism of $G^F$-varieties \begin{align*}
        X(\varpi)\cong Y_I(\dot{w})\times_{L_I^{\dot w F}}X_{L_I}(v).
    \end{align*}
    Note that $L_I$ is of type $A_1$ and hence $\dot w F$ is the standard Frobenius on $L_I$, which means $H^*_c(X_{L_I}(v))\cong \St_{L_I^{\dot wF}}[-1]\oplus \1_{L_I^{\dot wF}}(-1)[-2]$.
    As $\pi\colon Y_I(\dot w)\to X_I(w)$ is a $L_I^{\dot wF}$-torsor (see \cite[(7.9)]{Digne_2014}), the pushforward $\pi_! \Ql$ decomposes as the direct sum of two local systems parametrized by $\St_{L_I^{\dot wF}}$ and $\1_{L_I^{\dot wF}}$, which we denote by $\St$ and $\Ql$ respectively (see \cite[1.9]{originalpaper}).
    In particular, 
    \begin{align*}
        H^i_c(Y_I(\dot w),\Ql)&=H^i_c(X_I(w),\pi_!\Ql)\\
        &\cong H^i_c(X_I(w),\Ql)\oplus H^i_c(X_I(w),\St).
    \end{align*}
    By the Künneth formula and \cref{cohom_quot} \begin{align*}
        H^i_c(X(\varpi))&\cong H^{i-2}_c(Y_I(\dot w))\otimes_{\Ql[L_I^{\dot w F}]}\1(-1)\oplus H^{i-1}_c(Y_I(\dot w))\otimes_{\Ql[L_I^{\dot w F}]}\St\\
        &\cong H^{i-2}_c(X_I(w),\Ql)(-1)\oplus H^{i-1}_c(X_I(w),\St)
    \end{align*}
    as $\Mod(\Ql[G^F])$ is semisimple. 
    The principal series part of the cohomology of $X_I(w)$ has been computed in \cite[Proposition 3.13]{DUDAS20131} and is given by 
    \begin{align*}
            H^{3+k}_c(X_I(w),\Ql)_\pr\cong \begin{cases}
                j_{(0;21)}(-1)&\text{if }k=1\\
                j_{(1;2)}(-2)&\text{if }k=2\\
                \1(-4)&\text{if }k=5\\
                0&\text{else}
            \end{cases}    
    \end{align*}
    and 
    \begin{align*}
            H^{3+k}_c(X_I(w),\St)_\pr\cong \begin{cases}
                j_{(0;111)}&\text{if }k=1\\
                j_{(11;1)}(-2)&\text{if }k=2\\
                j_{(21;0)}(-3)&\text{if }k=3\\
                0&\text{else}.
            \end{cases}    
    \end{align*}
    Indeed, keeping in mind the following slight subtleties, our element $w$ corresponds precisely to $w_3$ defined in \cite[3.3.2]{DUDAS20131};
    \begin{enumerate}
        \item Dudas uses the same convention as GAP for the Dynkin diagram of $B_3$:\dynkin[backwards,labels={3,2,1},reverse arrows,label macro/.code={t_{#1}},
		edge length=.75cm]C{3}.
        \item There is a small typo in the proof of \cite[Proposition 3.13]{DUDAS20131}.
        Using the notions and notation of the preprint, the cohomology of the quotient of the piece $X_{x_2}$ should be shifted by $2$ instead of $1$, and a Tate twist $1$ should be added. 
        This is just the application of \cite[Remark 3.12]{DUDAS20131}.
        Moreover, Dudas assumes the convention $H^2_c(\A^1)=\Ql(1)$ for Tate twists (see the subsequent preprint \cite[Theorem 1.1]{DUDAS2013276} of Dudas for this clarification), which is precisely the opposite to ours. 
    \end{enumerate} 
    The claim follows.
\end{proof}

Using the tactic outlined at the beginning of \cref{small_rk_comps}, one gets the following $G^F$-isomorphisms
\begin{align*}
        i_{\{s_1,s_3\}}&\cong j_{(11;1)}\oplus j_{(21;0)}\oplus j_{(1;2)}\oplus j_{(2;1)}\oplus \1,\\
        i_{\{s_2,s_3\}}&\cong j_{(21;0)}\oplus j_{(2;1)}\oplus \1,\\
        i_{\{s_1,s_2\}}&\cong j_{(1;2)}\oplus j_{(2;1)}\oplus j_{(0;3)}\oplus \1,\\
        H^1(X(c))&\cong j_{(2;1)}.
\end{align*}
where $c\in W$ is a Coxeter element.
This gives the following.

\begin{proposition}\label{prop_type_B3}
    Let $G^F$ be of type $B_3$ and $\varpi\in W$ as in \cref{H1_height0_nonCox}.
    Let $w\in W\setminus \{w_0\}$ of full support and non-Coxeter.
    Up to choice of cyclic shift representative, the representation $\Pic(X(w))_{\Ql}$ is given in \cref{table:Sp6} where $\Pic(X(\varpi))_{\Ql}^{B^F}=0$.
    \begin{table}[h!]
    \centering
    \begin{tabular}{ m{8em} | m{5em} | m{5em} | m{15em}} 
    \makecell{$w$} & \makecell{$l(w)$} & \makecell{$|\Cyc(w)|$} & \makecell{$\Pic(X(w))_{\Ql}$} \\ [0.5ex] 
    \hline\hline
    \makecell{$s_3s_1s_2s_1$} & \makecell{$4$} & \makecell{$4$} & \makecell{$j_{(21;0)}$} \\
    \hline 
    \makecell{$s_2s_1s_3s_2$} & \makecell{$4$} & \makecell{$1$} & \makecell{$j_{(11;1)}\oplus j_{(21;0)}\oplus j_{(1;2)}$}  \\
    \hline 
    \makecell{$s_3s_1s_2s_3$} & \makecell{$4$} & \makecell{$2$} & \makecell{$j_{(1;2)}\oplus j_{(0;3)}$} \\
    \hline 
    \makecell{$s_3s_2s_1s_2s_3$} & \makecell{$5$} & \makecell{$1$} & \makecell{$j_{(21;0)}\oplus j_{(1;2)}\oplus j_{(2;1)}\oplus j_{(0;3)}$}   \\
    \hline 
    \makecell{$s_1s_2s_3s_2s_1$} & \makecell{$5$} & \makecell{$1$}& \makecell{$j_{(21;0)}^2\oplus j_{(2;1)}$} \\
    \hline 
    \makecell{$s_2s_1s_2s_3s_2$} & \makecell{$5$} & \makecell{$2$} & \makecell{$j_{(21;0)}$} \\
    \hline 
    \makecell{$\varpi$} & \makecell{$5$} & \makecell{$4$} & \makecell{$\Pic(X(\varpi))_{\Ql}$}  \\
    \hline 
    \makecell{$s_1s_2s_1s_3s_2s_1$} & \makecell{$6$} & \makecell{$1$} & \makecell{$j_{(21;0)}$} \\
    \hline 
    \makecell{$s_2s_3s_2s_3s_1s_2$} & \makecell{$6$} &\makecell{$4$} & \makecell{$\Pic(X(\varpi))_{\Ql}$} \\
    \hline 
    \makecell{$s_1s_3s_2s_3s_2s_1$} & \makecell{$6$} & \makecell{$1$} &\makecell{$\Pic(X(\varpi))_{\Ql}$} \\
    \hline 
    \makecell{$s_3s_2s_3s_1s_2s_3$} & \makecell{$6$} & \makecell{$1$} &\makecell{$\Pic(X(\varpi))_{\Ql}$} \\
    \hline 
    \makecell{$s_2s_3s_2s_1s_2s_3s_2$} & \makecell{$7$} & \makecell{$1$} &\makecell{$\Pic(X(\varpi))_{\Ql}$}
    \\
    \hline 
    \makecell{$s_2s_3s_1s_2s_3s_2s_1$} & \makecell{$7$} & \makecell{$4$} &\makecell{$\Pic(X(\varpi))_{\Ql}$}\\
    \hline 
    \makecell{$s_iw_0, i=1,2,3$} & \makecell{$8$} & \makecell{$1$} &\makecell{$\Pic(X(\varpi))_{\Ql}$}
    \end{tabular}
    \caption{Rational Picard group of $X(w)$ for $G^F$ of type $B_3$, where $w\neq w_0$ is of full support and non-Coxeter.}
    \label{table:Sp6}
    \end{table}
\end{proposition}

\subsubsection{Type \texorpdfstring{$G_2$}{G2}}\label{subsec:typeG2}

Let $G^F$ be of type $G_2$.
The group $G^F$ has $10$ unipotent characters, $6$ of which lie in the principal series. 
They are denoted by $\phi_{(1,0)}, \phi_{(2,1)}, \phi_{(2,2)}, \phi_{(1,3)}',\phi_{(1,3)}'',\phi_{(1,6)}$ with the convention that $\phi_{(1,0)}=\1$ and $\phi_{(1,6)}=\St$ (see \cite[Section 13.9]{Carter1972}).
Note that any Weyl group element length $4$ is terminal.
At least to the author, the cohomology ($\ell$-adic or motivic) of the corresponding Deligne--Lusztig varieties is unknown (see also \cref{rem:G2_comp_with_cohom}).
Using the tactic outlined at the beginning of \cref{small_rk_comps}, one gets the following $G^F$-isomorphisms
\begin{align*}
        i_{\{s_1\}}&\cong \phi_{(1,3)}'\oplus \phi_{(2,1)}\oplus \phi_{(2,2)}\oplus \1,\\
        i_{\{s_2\}}&\cong \phi_{(1,3)}''\oplus \phi_{(2,1)}\oplus \phi_{(2,2)}\oplus \1,\\
        H^3_c(X(c))&\cong \phi_{(2,1)}.
\end{align*}
where $c\in W$ is a Coxeter element.

\begin{proposition}\label{Pic_G2}
    Let $w\in W$ be a cyclic shift representative of length $3$. 
    Then 
    \begin{align*}
        \Pic(X(w))_{\Ql,\pr}&\cong \begin{cases}
            \phi_{(2,2)}+\phi_{(1,3)}''&\text{if }w=s_1s_2s_1\\
            \phi_{(2,2)}+\phi_{(1,3)}'&\text{if }w=s_2s_1s_2.
        \end{cases}
    \end{align*}
    In particular, $\Pic(X(w))_{\Ql,\pr}=\Pic(X(w))_{\Ql,\upr}$.
\end{proposition}

\begin{rem}\label{rem:G2_comp_with_cohom}
    The $\ell$-adic cohomology of the corresponding Deligne--Lusztig varieties has been partly computed in \cite[Théorème 4.4.3]{digne2006cohomologie}, in the sense that their table does not tell us in which degrees two particular unipotent representations, namely $\phi_{(2,2)}$ and one cuspidal, can appear. 
    \Cref{Pic_G2} together with the injection $\Pic(X(w))_{\Ql,\pr}\hookrightarrow H^2(X(w))_\pr$ of \cref{rq_DMR_cycle_iso_on_BF_3} tells us that $\phi_{(2,2)}$ does in fact appear in $H^4_c(X(w))\cong H^2(X(w))$ with nonzero multiplicity for $w\in W$ of length $3$.
\end{rem}

\bibliographystyle{alpha}
\bibliography{references}

\end{document}